\documentclass[10pt]{amsart}
\usepackage[T1]{fontenc}
\usepackage[utf8]{inputenc}
\usepackage[english]{babel}

\usepackage[a4paper]{geometry}%

\usepackage{amsmath, latexsym, amsfonts, amssymb}
\usepackage{amsthm}
\usepackage{mathrsfs,enumerate}
\usepackage{graphics,epsfig,color}

\usepackage{enumitem}
\usepackage[svgnames]{xcolor}
\usepackage{bbm}
\usepackage{comment}

\usepackage{hyperref}
\hypersetup{
    colorlinks=true,
    linkcolor=Blue,
    citecolor=Goldenrod,
    urlcolor=Blue,
    pdfborder={0 0 0}
}
\definecolor{fond}{rgb}{0.05,0.05,0.25}

\DeclareMathOperator{\mindist}{mindist}

\DeclareMathSymbol{\eset}{\mathalpha}{AMSb}{"3F}

\usepackage{enumitem}
\setlist[itemize,1]{nosep}
\setlist[enumerate,1]{itemsep=0pt,label=(\alph*)}

\numberwithin{equation}{section}

\theoremstyle{definition}
\newtheorem{thm}{Theorem}[section]

\newtheorem{example}[thm]{Example}

\newtheorem{lemma}[thm]{Lemma}
\newtheorem{corollary}[thm]{Corollary}
\newtheorem{prop}[thm]{Proposition}
\newtheorem{defn}[thm]{Definition}

\newtheorem{remark}[thm]{Remark}

\newcommand{\cA}{{\ensuremath{\mathcal A}} }
\newcommand{\cB}{{\ensuremath{\mathcal B}} }

\newcommand{\cE}{{\ensuremath{\mathcal E}} }
\newcommand{\cF}{{\ensuremath{\mathcal F}} }

\newcommand{\cL}{{\ensuremath{\mathcal L}} }

\newcommand{\cP}{{\ensuremath{\mathcal P}} }
\newcommand{\cR}{{\ensuremath{\mathcal R}} }

\newcommand{\bbE}{{\ensuremath{\mathbb E}} }
\newcommand{\E}{{\ensuremath{\mathbb E}} } 
\renewcommand{\P}{{\ensuremath{\mathbb P}} } 

\newcommand{\bbH}{{\ensuremath{\mathbb H}} }

\newcommand{\bbP}{{\ensuremath{\mathbb P}} }

\newcommand{\bbR}{{\ensuremath{\mathbb R}} }

\newcommand{\bbZ}{{\ensuremath{\mathbb Z}} }

\newcommand{\abs}[1]{\left\lvert #1 \right\rvert}

\newcommand{\ag}{\left\{ } 

\newcommand{\ad}{\right\} }

\newcommand{\cg}{\left[}
\newcommand{\cd}{\right]}
\newcommand{\pg}{\left(} 
\newcommand{\pd}{\right)}
\newcommand{\bg}{\left|}
\newcommand{\bd}{\right|}

\newcommand{\du}{{\ensuremath{\;:\;}}} 

\newcommand{\summ}[2]{\sum_{#1}^{#2}}

\newcommand{\intt}[2]{\int_{#1}^{#2}}

\newcommand{\esppi}[1]{\mathbb{E}_\pi \left[ #1 \right] }
\newcommand{\bbPpi}[1]{\mathbb{P}_\pi \left( #1 \right) }
\newcommand{\eps}{\varepsilon}

\newcommand{\report}[1]{\textcolor{black}{#1}}

\def\Z{\mathbb{Z}}

\def\Z{\mathbb{Z}}

\DeclareMathOperator{\tmix}{\textit{t}_{mix}}
\DeclareMathOperator{\tmixun}{\textit{t}_{mix}^{(1)}}
\DeclareMathOperator{\tmixinfty}{\textit{t}_{mix}^{(\infty)}}
\DeclareMathOperator{\tmixp}{\textit{t}_{mix}^{(p)}}
\DeclareMathOperator{\tmixq}{\textit{t}_{mix}^{(q)}}

\DeclareMathOperator{\trel}{\textit{t}_{rel}}

\DeclareMathOperator{\taucov}{\tau_{cov}}
\DeclareMathOperator{\tcov}{\textit{t}_{cov}}
\DeclareMathOperator{\thit}{\textit{t}_{hit}}
\DeclareMathOperator{\thitav}{\textit{t}_{\langle{hit}\rangle}}

\DeclareMathOperator{\Unif}{\text{Unif}}
\DeclareMathOperator{\dtv}{\text{d}_\text{TV}}
\DeclareMathOperator{\duncov}{\text{d}_\text{uncov}}
\DeclareMathOperator{\Aut}{\text{Aut}}
\DeclareMathOperator{\Diam}{\text{Diam}}
\DeclareMathOperator{\Cay}{\text{Cay}}

\DeclareMathOperator{\tsmax}{\textit{t}_{\textit{s}}}

\newcommand{\tsav}{{\ensuremath{t_{ {\scriptscriptstyle \langle} \hspace{-0.03cm}s {\scriptscriptstyle \rangle} } }}}

\newcommand{\tsprimav}{{\ensuremath{t_{ {\scriptscriptstyle \langle} \hspace{-0.03cm}s' {\scriptscriptstyle \rangle} } }}}

\newcommand{\tsavminuslittleoofone}{{\ensuremath{t_{ {\scriptscriptstyle \langle} \hspace{-0.03cm}s -o(1){\scriptscriptstyle \rangle} } }}}

\DeclareMathOperator{\tsexp}{\textit{t}^*_{\textit{s}}}

\usepackage[normalem]{ulem}

\numberwithin{equation}{section}

\begin{document}

\title[Universality of fluctuations for the cover time]{On the universality of fluctuations for the cover time}

\author{Nathanaël Berestycki}
\address{Universität Wien}
\curraddr{}
\email{nathanael.berestycki@univie.ac.at}
\thanks{N. B.'s research is supported by FWF grant P33083, ``Scaling limits in random conformal geometry''.}

\author{Jonathan Hermon}
\address{University of British Columbia}
\email{jhermon@math.ubc.ca}
\thanks{J. H.'s research is supported by an NSERC grant.}

\author{Lucas Teyssier}
\address{Université de Lorraine}
\email{lucas.teyssier@univ-lorraine.fr}
\thanks{L. T.'s research is supported by FWF grant P33083, ``Scaling limits in random conformal geometry''.}

\subjclass[2020]{}

\date{\today}

\dedicatory{}


\begin{abstract}
We consider random walks on finite vertex-transitive graphs $\Gamma$ of bounded degree. We find a simple geometric condition which characterises the cover time fluctuations: the renormalised cover time $\tfrac{\taucov}{\thit} - \log |\Gamma|$ converges to a standard Gumbel variable if and only if $\mathrm{Diam}(\Gamma)^2 = o(n/\log n)$, where $n = |\Gamma|$. We prove that this condition is furthermore equivalent to the decorrelation of the uncovered set.
The arguments rely on recent breakthroughs by Tessera and Tointon on finitary versions of Gromov's theorem on groups of polynomial growth, which we leverage into strong heat kernel bounds, and refined quantitative estimates on Aldous and Brown's exponential approximation of hitting times.
\end{abstract}

\maketitle

\tableofcontents

\section{Introduction}\label{S:intro}

\subsection{Context}
\label{SS:mainresults}

Let $\Gamma$ be a finite (connected) vertex-transitive graph. Let $X= (X_t)_{t \ge 0}$ denote the simple random walk in continuous time, which at constant rate 1 jumps to a randomly chosen neighbour, starting from a designated vertex called the root and denoted by $o$. Let us denote by $T_x = \inf \ag t\geq 0\du X_t = x \ad$ the hitting time of the vertex $x$.
The \textbf{cover time} variable of $\Gamma$ by $X$ is the first time that the walk has visited every vertex:
\begin{equation*}
\taucov = \max_{x\in \Gamma} T_x,
\end{equation*}
where with a slight abuse of notation we have used $\Gamma$ to also denote the vertex set of the graph $\Gamma$. 
In this article we are concerned with obtaining general fluctuation results for the cover time $\taucov$ of random walks on vertex-transitive graphs,  as the vertex set size $|\Gamma|$ tends to infinity. As we will see this question is deeply intertwined with the study of the structure of the uncovered set $U(t) = \{ x\in \Gamma: T_x > t\}$,
for $t$ close to the (expected) cover time.

\medskip Obtaining quantitative estimates on the cover time is a natural problem which parallels the intensively studied question of mixing time and cutoff for random walks on graphs (i.e.\ understanding how far the law of $X_t$ is from the stationary distribution, see Section \ref{S:relaxmixhittimes} for some definitions, and \cite{LivreLevinPeres2019MarkovChainsAndMixingTimesSecondEdition} and the references therein for an introduction). In common with  much of the literature on the subject, we will focus in this paper on vertex-transitive graphs. The restriction to this class of graphs is natural in order to avoid pathological examples, which can be arbitrarily badly behaved. At the same time it allows  
for a very rich range of behaviours, as we are about to discuss.

\medskip Results in this direction go back at least to the seminal work of Aldous \cite{Aldous1983cover} who proved that for random walks on finite groups and under mild geometric conditions, the cover time $\taucov$ is concentrated around its mean $\tcov = \E( \taucov)$ (since the distribution of $\taucov$ does not depend on the starting point, we may write $\bbE$ and $\bbP$ without specifying the starting measure) i.e.\ that $\tfrac{\taucov}{\tcov} \to 1$ in probability, and obtained the leading order behaviour of $\tcov$.

This was complemented a few years later by Matthews \cite{Matthews1988cover} who obtained a general upper bound (valid for \emph{any} graph) on the expected cover time.  The problem of the correlation structure of the uncovered set was raised in the physics literature through the work of Brummelhuis and Hilhorst \cite{BH91, BH92}.

\medskip It is worth noting that even for basic graphs such as the $d$-dimensional torus of side-length $n$, i.e.\ when $\Gamma = (\Z/ n\Z)^d$ (with $d$ fixed and $n \to \infty$), the problem of describing the cover time is highly nontrivial. {For more general Cayley graphs, even very concrete ones such as of the Heisenberg group (see Example \ref{ex:H} for details), the problem was until now completely open, even though the parallel problem of convergence to equilibrium and mixing of random walk on such groups of moderate growth has been studied at least since the seminal works of Diaconis and Saloff-Coste \cite{DiaconisSaloff-Coste1993, DSC}.}

Let us first summarise what is known in the case of the $d$-dimensional torus. It was only in 2004 that the first order of the cover time for the two-dimensional torus $\pg \bbZ/n\bbZ\pd^2$ was found by Dembo, Peres, Rosen and Zeitouni in \cite{DemboPeresRosenZeitouni2004cover}. This result was later refined, see \cite{BeliusKistler2017, Abe2021cover, BeliusRosenZeitouni2020}, but obtaining a convergence in distribution for the fluctuations of the cover time in two dimensions remains an open problem to this day (it is widely believed however that the fluctuations will be in any case different from the higher dimensional regime described below). 

In dimensions $d\ge 3$, it was proved by Belius \cite{Belius2013} in 2013 (solving an open question of Aldous and Fill \cite{AldousFill}), building on Sznitman's random interlacement model \cite{Sznitman2010} (which
 was in fact motivated by the work of Brummelhuis and Hilhorst \cite{BH91, BH92}), that the fluctuations of the cover time are asymptotically distributed according to a standard Gumbel law:
\begin{equation*}
\bbP\pg\frac{\taucov}{\thit} - \log |\Gamma| \leq s\pd \to e^{-e^{-s}},
\end{equation*}
where $\thit := \max_{x,y\in\Gamma} \bbE_x \cg T_y \cd$ is the maximal expected hitting time. (See also the results in Prata's thesis \cite{Prata} for partial results valid for more general graphs but under some restrictive assumptions.)

The Gumbel law is significant because it describes the asymptotic maximum of i.i.d.\ random variables, subject to some conditions on their common distribution.
The Gumbel fluctuations in the result above therefore suggest that the law of uncovered set at time $\tsmax := \thit(\log(|\Gamma|)+s)$ might  asymptotically be close to a product measure, where each vertex is uncovered with probability $e^{-s}/|\Gamma|$ independently of other vertices.
More formally, if $\mu_s$ is a Bernoulli variable of parameter $e^{-s}/|\Gamma|$, we might expect that as $|\Gamma| \to \infty$,
\begin{equation}\label{e:decor}
    \duncov(\tsmax) := \dtv\pg \cL(U(\tsmax)), \mu_s^{\otimes \Gamma} \pd \to 0,
\end{equation}
where $\cL(U(\tsmax))$ denotes the law of the uncovered set at time $\tsmax$ \report{(for the simple random walk on $\Gamma$ started at $o$)}, and the total variation distance between two probability measures $\mu$ and $\nu$ on a finite space $S$ is given by
\begin{equation*}
    \dtv(\mu,\nu) = \max_{A\subset S} |\mu(A) - \nu(A)|.
\end{equation*}
Part of the goal of this paper is to study the uncovered set and in particular prove \eqref{e:decor} in a general framework. 

Questions concerning the geometry of the uncovered set (even when the cover time is well understood) have recently become prominent.
Even on the $d$-dimensional torus, a number of basic problems remain open. For instance, it was proved in \cite{MillerSousi} and \cite{OleskerTaylorSousi2020} that the uncovered set \report{at time $a t^*$ (with $t^* = (1+ o(1))\tcov$)} is decorrelated in the above sense if [$a>7/8$ for all $d\geq 3$] or if [$a>3/4$ for $d$ large enough depending on $a$] and correlated (in the sense that \eqref{e:decor} does not hold) if $a< a_d = \tfrac{1+p_d}{2}$, where $p_d = \bbP_o(T_o^+ < \infty)$ for the simple random walk on $\bbZ^d$. 
Recently a substantial progress was made by Prévost, Rodriguez, and Sousi \cite{PrevostRodriguezSousi}, who established a phase transition at $a_d$ for a relaxed version of the problem (involving two-sided stochastic domination and sprinkling). We refer to \cite[Remark 7.3.7]{PrevostRodriguezSousi} for a recent detailed account on this problem.
It nevertheless still remains an open question whether there actually is a phase transition in total variation, which is expected to occur at $a=a_d$. In other words, the problem is to show that $\duncov(a\tcov) \to 0$ for $a>a_d$.

More generally, the study of the uncovered set fits into the theme of exceptional points for random walks. In two dimensions the structure of those exceptional points has recently been proved to be linked with Liouville quantum gravity, see \cite{AbeBiskup2022}, and, away from the cover time, to an even more singular and in some way intriguing object called Brownian multiplicative chaos, see \cite{Jego2020}.

Finally, we note that the expected cover time $\tcov$ of a graph is also closely related to the typical value of the maximum of the associated Gaussian free field, see for instance \cite{DingLeePeres2012coverGFF}, \cite{Ding2014coverGFF} and \cite{Zhai2018}.

\subsection{Main results}

Let $\Gamma$ be a finite (connected) bounded degree vertex-transitive graph. Write $\Gamma$ also for the vertex set of $\Gamma$ as above. We denote by $d(x,y)$ the graph distance between the vertices $x$ and $y$, by $D=\Diam(\Gamma) := \max_{x,y \in \Gamma} d(x,y)$ the diameter of $\Gamma$, and by $\pi$ its stationary distribution.

Our first main result shows that Gumbel fluctuations are universal. Perhaps even more surprisingly we obtain a sharp (necessary and sufficient) geometric condition for this universality.

\begin{thm}
 \label{T:main} Let $(\Gamma)$ be a collection of finite (connected) vertex-transitive graphs of fixed degree, and let $\chi$ be a standard Gumbel variable, i.e.\ $\P( \chi \le s) = e^{- e^{-s}}$ for $s \in \bbR$. Then
 \begin{equation}\label{E:main}
 \frac{\taucov}{\thit} - \log |\Gamma|  \xrightarrow[|\Gamma|\to \infty]{d} \chi
 \end{equation}
 if and only if
 \begin{equation}
\label{e:DC}
\tag{DC}
\frac{D^2 \log |\Gamma|}{|\Gamma|} \xrightarrow[|\Gamma| \to \infty]{} 0.
\end{equation}
\end{thm}
For future reference we note that, trivially, if we let $n = |\Gamma|$ denote the number of vertices of $\Gamma$ (which tends to infinity by assumption), then \eqref{e:DC} is equivalent to $D^2 = o ( n /\log n)$, and \eqref{E:main} is equivalent to the condition that for all $s\in \bbR$, as $n \to \infty$, 
\begin{equation*}
\P( \taucov \le \thit ( \log n + s) ) \to \exp ( - e^{-s} ).
\end{equation*}
We complement this result with a refined statement on the structure of the uncovered set. 

\begin{thm}\label{T:uncovered_main}
Let $(\Gamma)$ be a collection of finite (connected) vertex-transitive graphs of fixed degree. For $s\in\bbR$, let $\tsmax := \thit(\log(|\Gamma|)+s)$ and let
$\mu_s^{\otimes \Gamma}$ denote the product over all the vertices of the graph of the Bernoulli law $\mu_s$ with parameter $e^{-s}/|\Gamma|$. Then, \report{under $\P_o$ or $\P_\pi$},
\begin{equation}\label{e:tv conv law uncov set}
    \dtv ( \cL(U(\tsmax)), \mu_s^{\otimes \Gamma} ) \xrightarrow[|\Gamma|\to\infty]{} 0
\end{equation}
  for every $s\in\bbR$ if and only if the diameter condition \eqref{e:DC} holds.
 \end{thm}

This result will be strengthened in various ways under the assumption \eqref{e:DC} in Section \ref{S:uncoveredset}. For instance, in Theorem \ref{T:convergence of the uncovered set} we discuss the law of the uncovered set at the first time that exactly $k$ points remain to be covered. 

\medskip In particular, we recover that \eqref{e:decor} holds for the tori $(\bbZ/m\bbZ)^d$, for $d\geq 3$ fixed and as $m\to \infty$ (this is in fact already mentioned in the thesis of Prata \cite{Prata}, although the result was never published). We also point out that despite considerable effort over the last 50 years, in the parallel problem of the mixing time mentioned above,
a simple characterisation of the cutoff phenomenon (and even more so of the limiting profile) is currently completely out of reach.

\subsection{Application to Cayley graphs}
Cayley graphs have very diverse geometric behaviours. For example, even in graphs of polynomial growth, the size of balls can grow slowly (in tori), or initially very fast (e.g.\ in the product of a cycle and a Ramanujan Cayley graph, as we will detail in Section \ref{s:examples}). Understanding cover times with specific techniques could be very challenging, and was open even for thin tori. However, as we will see, in many natural examples, it is simple to check the diameter condition \eqref{e:DC}, and hence to deduce the cover time behaviour.
\begin{example}[Thin tori]
The first and most natural example to consider is the ``thin'' torus, which is nothing but a box (with periodic boundary conditions) of sidelengths $m$, $m$ and (height) $h = h_m$ (where $1\le h_m \le m$). Formally this is the Cayley graph $\Gamma = \Cay(G,S)$ of the group $G = (\Z/ m\Z)^2 \times (\Z/ h \Z)$  with respect to the generating set $S = \ag (\pm 1,0,0), (0, \pm 1, 0), (0,0,\pm 1) \ad$.

Thus the extreme cases $h_m =1$ and $h_m = m$ correspond to the familiar two- and three-dimensional cases respectively, while very little is known in general about intermediate cases. The diameter of $\Gamma$ is $D \asymp m$ (where $a_n \asymp b_n$ means that $a_n=O(b_n)$ and $b_n = O(a_n)$), and the volume is $n = m^2 h_m$, so the condition \eqref{e:DC} holds if and only if $h_m / \log m \to \infty$.
\end{example}
\begin{remark}
    The same analysis applies to more general tori $\prod_{i=1}^{r} (\bbZ/m_i \bbZ)$ (where $r$ is fixed and $m_r \leq ... \leq m_2 \leq m_1 \to \infty$). Then the diameter is $\asymp m_1$, so the diameter condition can be rewritten as $m_1^2 \log m_1 =o\pg \prod_{i=1}^r m_i\pd$.
\end{remark}

\begin{example}[Heisenberg group]\label{ex:H}
    For $m\geq 1$, the Heisenberg group $G = \bbH_3(m)$ is the group of upper triangular matrices with coefficients in $\bbZ/m\bbZ$ and unit diagonal. A natural set of generators is 
    \begin{equation*}
        S = \ag 
    \begin{pmatrix}
1 & \pm 1 & 0 \\
0 & 1 & 0\\
0 & 0 & 1
\end{pmatrix},
\begin{pmatrix}
1 & 0 & 0 \\
0 & 1 & \pm 1\\
0 & 0 & 1
\end{pmatrix} 
\ad.
    \end{equation*}
Let $\Gamma = \Cay(\bbH_3(m),S)$ be the corresponding Cayley graph.  It is relatively easy to check (see for instance \cite[Example 1]{DiaconisSaloff-Coste1993}),  that $\Diam(\Gamma) \asymp m$, so, since $\bg \bbH_3(m) \bd = m^3$, the diameter condition is satisfied as $m\to \infty$. A similar analysis can be carried out for generalisations of the Heisenberg group such as upper triangular matrices with unit diagonal, see \cite[Example 3]{DiaconisSaloff-Coste1993}.
\end{example}

\begin{example}
This example has a more arithmetic flavour.
Consider, for $p\geq 3$ prime, the group $G =M_3(p)$ obtained as the semi-direct product of $\bbZ/(p\bbZ)$ and $\bbZ/(p^2\bbZ)$ with multiplication law given by 
\begin{equation*}
    (a,b)(c,d) = (a+c, cb+d) \quad \quad \text{for $(a,b), (c,d) \in \bbZ/(p\bbZ) \times \bbZ/(p^2\bbZ)$}. 
\end{equation*}
A natural set of generators is $S := \ag (\pm1,0), (0, \pm1)\ad$. As shown in \cite[Example 2]{DiaconisSaloff-Coste1993}, the Cayley graph $\Gamma = \Cay(G,S)$ has diameter $\asymp p$ as $p\to \infty$, so, since it has $p^3$ vertices, it satisfies \eqref{e:DC}.
\end{example}

\begin{example}[Product of a cycle and a Ramanujan Cayley graph]
    Consider a product of a cycle of length $m$ and a Ramanujan Cayley graph of size $\ell = \ell_m$, where $\ell_m$ is a sequence which may be chosen arbitrarily (see Section 
\ref{s:examples} for precise definitions). If $\ell_m$ is chosen to grow at most polynomially in $m$, i.e.\ if $\ell_m = m^{O(1)}$, then the graph is \textit{macroscopically} of
polynomial growth (in the sense that the volume is bounded by a power of its diameter). Nevertheless its geometric behaviour is markedly different from the above examples: essentially, the graph is locally ``very recurrent'' in one direction but ``very transient'' in another. In particular, the volume growth of balls is initially exponential. 
Despite this difference with the above examples, it is straightforward to check that there are instances when the diameter condition \eqref{e:DC} is satisfied: more precisely, this condition is satisfied if and only if $\ell_m /(m\log m) \to \infty$. We will return to this example in detail in Section \ref{s:examples}, as we will see that this gives us a source of examples of graphs which are locally transient in a precise sense (more precisely in the sense of ``strong uniform transience'', which we introduce below in Section \ref{SS:localtransience}) but to which the diameter condition does not apply. In particular, the cover time does not have Gumbel fluctuations or, equivalently (see the next subsection) the uncovered set has nontrivial correlations within the scaling window.
\end{example}

\subsection{Decorrelation in the uncovered set} 

The reader might find it surprising that the condition \eqref{e:DC} determines the behaviour \eqref{E:main} (and \eqref{e:tv conv law uncov set}). Theorem \ref{T:uncovered_main} tells us that decorrelation in the uncovered set close to the cover time depends only on the relation of the diameter of the graph to its size, a completely global geometric condition. 
The local structure of the graph only plays a role in determining the leading order behaviour of $\thit$.

In the proof that \eqref{e:DC} is necessary for \eqref{e:tv conv law uncov set}, we will see that even the first moment of $U(\tsmax)$ is different from that of $\mu_s^{\otimes \Gamma}$. In that sense the theorem above does not capture the correlations of the uncovered set close to the cover time. However, there are other natural ways to adjust the time scaling, for instance by considering (for any $s \in \mathbb{R}$) the uncovered set at the time $\tsexp$ such that $\E(|U(\tsexp)|) = e^{-s}$. 
Such changes of normalisation do not affect the above results under \eqref{e:DC}. That is, under \eqref{e:DC}, both \eqref{E:main} and \eqref{e:tv conv law uncov set} hold with $\tsmax$ replaced by $\tsexp$. Indeed, we will prove in Proposition \ref{P:thitmaxav} and Proposition \ref{P:tdoubleprime s} that under the assumption \eqref{e:DC}, we have $\tsmax = \tsexp + o(D^2)$, where $D = \Diam(\Gamma)$ is the diameter of $\Gamma$. 

It is therefore natural to ask whether the law of $U(\tsexp)$ is close to a product measure. For this as well, we will show that the diameter condition \eqref{e:DC} is the sharp criterion for decorrelation. However, in order to state this stronger result we will need to make an additional geometric assumption. 

\begin{thm}\label{T:uncovered intro t hit av}
    Let $(\Gamma)$ be a collection of finite (connected) vertex-transitive graphs of fixed degree such that $\Diam(\Gamma)^2/|\Gamma| \to 0$ as $|\Gamma|\to \infty$. 
    For $s\in\bbR$, let
$\mu_s^{\otimes \Gamma}$ denote the product over all the vertices of the graph of the Bernoulli law $\mu_s$ with parameter $e^{-s}/|\Gamma|$. Fix $s\in \bbR$. The following are equivalent:
\begin{enumerate}[label=(\roman*)]
\item 
  $
  \dtv ( \cL(U(\tsexp)), \mu_s^{\otimes \Gamma} ) \xrightarrow[|\Gamma|\to\infty]{} 0,
  $

\item The diameter condition \eqref{e:DC} holds.
\end{enumerate}
\end{thm}

\begin{remark}
    All our statements for collections $(\Gamma)$ of graphs of \textit{fixed degree} are equivalent to statements for graphs of \textit{uniformly bounded degree}. Either this is immediate by inclusion or this follows from taking a maximum for the other direction.
\end{remark}

The \textit{relaxation time} of the simple random walk $X = (X_t)_{t\geq 0}$ on $\Gamma$ is defined by $\trel = 1/(1-\lambda_2)$, where $\lambda_2$ is the second largest eigenvalue of the transition matrix of $X$. If $\mu$ is a measure on (the vertex-set of) $\Gamma$, we also write $\bbP_\mu(\cdot) = \sum_{x\in \Gamma}\mu(x)\bbP_x(\cdot)$ and $\bbE_\mu(\cdot) = \sum_{x\in \Gamma}\mu(x)\bbE_x(\cdot)$.
\begin{remark}
    In fact, we will prove a result stronger than Theorem \ref{T:uncovered intro t hit av}: namely, the conclusion of Theorem \ref{T:uncovered intro t hit av} is valid as soon as $\trel = o(\thit)$, where $\trel$ is the \emph{relaxation time}, or inverse spectral gap. We will also show that these conditions are equivalent to a third one, where instead of $\tsexp $ we use the time $\tsav:=\thitav  ( \log n + s)$, where $\thitav := \bbE_\pi T_o$; this is in fact an integral part of our proof. The same statement also holds with $\tsmax$ instead of $\tsexp$. In other words, if $\trel = o(\thit)$ holds, we can strengthen Theorem \ref{T:uncovered_main}, replacing “for every $s\in \bbR$” by “for some $s\in\bbR$”.
    The theorem will be proved in these forms in Section \ref{SS:diam_nec_uncovered_av}, where we will also briefly explain why the assumption $D^2/ n \to 0$ implies $\trel = o(\thit)$. 

    We conjecture in fact that the theorem is valid without any assumption on $\trel$ or $\thit$.
    \end{remark}
    
\subsection{Diameter condition and local transience}
\label{SS:localtransience}

We have already mentioned that the sharpness of the diameter condition \eqref{e:DC} is a little surprising. Initially (and this was in fact our own belief when we began this work), one might have suspected that Gumbel fluctuations are perhaps more naturally linked with the following notions of local transience which we now define. For this it will be useful to recall the definition of the mixing time $\tmix (\varepsilon)$ at level $0<\varepsilon<1$ for a positive recurrent Markov chain $(X_{t})_{t \in \mathbb{R}_+} $ on some state space $\Gamma$ with invariant distribution $\pi$:
\begin{equation}\label{e:def tv mixing time at level eps}
   \tmix(\varepsilon) = \inf \{ t \ge 0: \sup_{x\in \Gamma} \dtv (p_t(x, \cdot) , \pi(\cdot)) \le \varepsilon \}, 
\end{equation}
where $p_t(x,\cdot)$ denotes the law of the Markov chain at time $t$ starting from $x\in \Gamma$.

\begin{defn}
\label{def:UT} We say that a sequence of finite Markov chains with state spaces $\Gamma_n$ and stationary distributions $\pi=\pi_n$ satisfying $\lim_{n \to \infty} \max_{x \in \Gamma_n}\pi_n(x)=0 $  is \emph{weakly uniformly transient} (\textbf{WUT}), or \emph{uniformly locally transient}, if
\begin{equation}
\label{e:WUT}
\max_{o \in \Gamma_n} \E_o ( L_o(t)) = O(1),
\end{equation}
where $t = \tmix (1/4)$, and $L_x(t) := \int_0^t  1_{\{ X_s = x \}} \mathrm{d}s$ is the \emph{local time} of the walk at $x$ up to time $t$. Again, writing  $t= \tmix (1/4)$,
we say that the sequence is (\textbf{SUT}) \emph{strongly uniformly transient}, or  \emph{uniformly globally transient},   if
\begin{equation}
\label{e:SUT}
\lim_{s \to \infty} \limsup_{n \to \infty} \max_{o \in \Gamma_n}  \E_o ( L_o(t)- L_o(s)) =0. \end{equation}
We say that a sequence of graphs $(\Gamma_n)$ is WUT (resp.\ SUT) if the sequence of simple random walks on $\Gamma_n$ is  WUT (resp.\ SUT).
\end{defn}

Clearly, \eqref{e:SUT} implies \eqref{e:WUT}. The reason why one might suspect that this notion might be related to Gumbel fluctuations is that uniform transience, especially strong uniform transience, should prevent clusterisation and hence lead to decorrelation in the uncovered set.

\medskip For example, a torus in dimension $d\ge 3$ is strongly uniformly transient, but a torus of dimension $d =1,2$ is not even weakly uniformly transient. More generally, let us return to the example of the thin torus $\Gamma = (\Z/m\Z)^2 \times ( \Z/ h\Z)$ discussed above, where $h = h_m$. For which values of $h$ is this weakly/strongly uniformly transient? Since a two-dimensional random walk returns to the origin roughly $\log m$ times by time $t =  \tmix (1/4) \asymp m^2$, it is not hard to see that the thin torus is strongly uniformly transient if and only if $h_m/ \log m \to \infty$ (and weakly uniformly transient if and only if $ (\log m)/h_m =O(1)$). This condition coincides with \eqref{e:DC}, as already observed.

\medskip This immediately raises the following question: could it be the case that the diameter condition is equivalent to strong (or weak) uniform transience? We will see as part of our analysis that any sequence of vertex-transitive graphs $(\Gamma_n)$  satisfying the diameter condition \eqref{e:DC} satisfies that
\begin{equation*}
    \lim_{r \to \infty} \limsup_{n \to \infty} \max_{y \in \Gamma_n :d(o,y) \ge r}  \E_y ( L_o(2\tmix (1/4))) =0, 
\end{equation*}
where $d(o,y)$ is the graph distance between $y$ and $o$ 
(see, e.g., Lemma \ref{proppxygrand}). Strong uniform transience can relatively easily be deduced from this  (in Proposition \ref{P: finite time spent at the origin before time n} we give a more direct argument). Thus
\begin{equation}\label{DCSUT}
\text{\eqref{e:DC} implies (SUT). }
\end{equation}
However the converse is, perhaps surprisingly again, not true. In Section \ref{s:examples} we will construct a sequence of finite vertex-transitive graphs $\Gamma_n$ of uniformly bounded degree which satisfy SUT but not \eqref{e:DC}. In particular, as a consequence of Theorem \ref{T:uncovered_main}, despite being locally transient in this strong sense, the uncovered set will display nontrivial correlations. This example will be constructed by considering the product of a Ramanujan graph and a cycle of suitable sizes. (Essentially, in this product, one direction is very recurrent, but the other is very transient). See Section \ref{s:examples} for the precise definition.

\medskip The implication \eqref{DCSUT} is closely related to a conjecture of Benjamini and Kozma \cite{BenjaminiKozma}. 
This conjecture states that for vertex-transitive graphs of bounded degree, if we assume $D^2 =O(n/\log n)$ then the effective resistance $\cR_{x \leftrightarrow y}$ between vertices of the graph is uniformly bounded above by some constant. This conjecture was recently solved by Tessera and Tointon \cite{TesseraTointonfinitary}. Furthermore, using the tools developed in their paper, it is not hard to show that the effective resistance is uniformly bounded if and only if (WUT) holds. The implication \eqref{DCSUT} is therefore the direct ``strong'' analogue of the Benjamini--Kozma conjecture (the left hand side replaces the assumption $D^2 =O(n/\log n)$ by $D^2 = o(n/\log n)$, and the weak uniform transience in the right hand side by the strong uniform transience).

\medskip \medskip We end this discussion with an instructive comparison with the results in a recent paper of Dembo, Ding, Miller and Peres \cite{DemboDingMillerPeres2019lamplighter}. This describes the first order (cutoff) behaviour of the mixing time of the lamplighter walk on the thin torus 
$\Gamma_m = (\Z/m\Z)^2 \times (\Z/h\Z)$ with $h = a \log m$. The (total variation) mixing time of the lamplighter group on the base graph $\Gamma_m$ is known to be closely related to the cover time of the simple random walk on that base graph $\Gamma_m$. This led the authors of \cite{DemboDingMillerPeres2019lamplighter} to use the ratio between the mixing time of the lamplighter on $\Gamma_n$ and the cover time on $\Gamma_n$ as an indicator of low- vs. high-dimensional behaviour. Indeed, this ratio is asymptotically equal to $1$ for a completely flat, two-dimensional torus by results of \cite{DemboPeresRosenZeitouni2004cover} and \cite{PeresRevelle2004}, whereas it is asymptotically $1/2$ for a three-dimensional torus, by results of Miller and Peres \cite{MillerPeres2012} (as well as on the complete graph, where the lamplighter walk reduces to the well known walk on the hypercube). Surprisingly, the authors in \cite{DemboDingMillerPeres2019lamplighter} show that on a thin torus, this ratio is strictly contained in the interval $(1/2, 1)$ for $0 < a < a_*$ and becomes asymptotically equal to 1/2 (as in the three-dimensional case) for $a \ge  a_*$, for some explicit $a_*$. They interpret this as a phase transition between low-dimensional (``recurrent'') and high-dimensional (``transient'') behaviour; see the discussion after Theorem 1.3 in that paper. By contrast, our results show that, at the level of fluctuations, high-dimensional behaviour only kicks in when $a =a_m \to \infty$, arbitrarily slowly, rather than for $a\ge  a_*$. 

\subsection{Discussion of proof ideas and organisation of paper}\label{s: discussion of proof ideas}

Our starting point for this paper is the remarkable series of papers by Tessera and Tointon \cite{TesseraTointonfinitary, TesseraTointonIsop}, which gives a \textbf{quantitative form of Gromov's theorem} on groups of polynomial growth. Recall that, since the graph is vertex-transitive, by results of Trofimov \cite{Trofimov1984} it is roughly isometric to the Cayley graph of a finite group.
Recall also that, for infinite groups, Gromov's theorem \cite{Gromov1981} shows that if the volume of balls of radius $r$ grows polynomially in the radius $r$ then the growth exponent $\alpha$ must be integer and the group is then (in the Gromov--Hausdorff sense) close to a $d$-dimensional Euclidean lattice.

The diameter condition \eqref{e:DC}, which says that the graph is slightly more than two-dimensional, combined with the results of Tessera and Tointon, therefore implies that at least for relatively moderate distances the volume growth is at least \textbf{three-dimensional}. In combination with isoperimetric profile bounds (coming for instance from the theory of \textbf{evolving sets} of Morris and Peres \cite{MorrisPeres2005}), this translates into very good decay for the heat kernel at small times and so gives excellent control on the number of returns of the random walk to its starting point in this time. Later visits to this point are controlled using two-dimensional estimates. As we will see now, it turns out this is at the root of the decorrelation in the uncovered set.

\report{These heat kernel bounds are used as follows. Suppose that \eqref{e:DC} holds and we wish to prove \eqref{E:main}.
It suffices to show that for each $k$, the $k$-th cumulant of the size $Z_s$ of the uncovered set at time $\tsav = \thitav ( \log n +s) $ converges to the $k$-th cumulant of a Poisson random variable with parameter $e^{-s}$, as $|\Gamma| \to \infty$. (As we will see $\tsav$ is very close to $\tsmax$, under \eqref{e:DC}.) Now, this cumulant can be expressed as a sum over all sets $A$ of size $k$ of the probability that $A$ was not visited by time $\tsav$, and by transitivity we can assume without loss of generality that the starting distribution is the uniform distribution $\pi$. The quantity we want to estimate therefore becomes $\P_\pi (T_A > \tsav)$.}

\report{The Aldous--Brown approximation \cite{AldousBrown1992} shows that this hitting time $\P_\pi (T_A > \tsav)$ is approximately an exponential random variable with mean approximately $\bbE_\pi[T_A]$, i.e.\ that $\P_\pi (T_A > \tsav)$ is close to 
\begin{equation}
    \exp\pg -\frac{\tsav}{\bbE_\pi[T_A]}\pd =  \exp\pg - \frac{\bbE_\pi[T_o]}{\bbE_\pi[T_A]} (\log n +s)\pd,
\end{equation}
and improved bounds from \cite{BerestyckiHermonTeyssier2025AB} make this approximation sufficiently precise to carry out the analysis under \eqref{e:DC}. This brings us to study,  for all subsets $A$ of $\Gamma$ of a given size $k$, the ratio 
\begin{equation}\label{qA}
    q_A: = \frac{\E_\pi[T_o]}{ \E_\pi [T_A]} 
\end{equation}
of expectations of hitting time of $A$ compared to that of a single point; see Proposition \ref{propapproxpita}. The quantity $q_A$ is analogous to the notion of capacity for infinite graphs.} For sets $A$ such that the points in $A$ are well separated, we typically expect $q_A \approx k$, because the hitting time of $A$ is close to the minimum of $k$ independent exponential random variables.
The challenge is to quantify this approximation and show that the contribution coming from sets where the points are not so well separated is negligible. It is here that the heat kernel bounds are very useful: indeed, the on-diagonal decay of $p_t(x,y)$ translates into a \textbf{strong off-diagonal decay} (using a subgaussian estimate, namely a recent variant due to Folz \cite{Folz2011}) and implies good quantitative bounds of the desired form for $q_A$.
In the most delicate case where $n $ is barely larger than $D^2 \log n$, the analysis uses a somewhat elaborate induction over scales which requires strong quantitative bounds.

In the opposite direction, when \eqref{e:DC} fails, the key task (both for Theorems \ref{T:main} and \ref{T:uncovered_main}) is to show that the expected number of uncovered points at time $\tsav$ is strictly smaller (by a factor $c$ bounded away from 1) than $e^{-s}$. Indeed this immediately implies Theorem \ref{T:uncovered_main}, and implies Theorem \ref{T:main} by considering the tail at $+\infty$ of $\taucov$ and a simple union bound. 

The proof of Theorem \ref{T:uncovered intro t hit av} is more complicated. 
It might initially be tempting to show that the uncovered set is ``too'' clustered, i.e.\ the probability that two relatively nearby points are uncovered is larger than it should be in an independent scenario. However, this turns out to be very difficult to control (moments of order at least two of the size of the uncovered set can explode when \eqref{e:DC} does not hold, precisely because of the contribution coming from nearby points). Instead we show that the uncovered set is negatively correlated at large distances. This requires controlling the macroscopic variations of the Green function, a task which occupies a good part of Section \ref{S:subcritical cases}.

\medskip \textbf{Organisation of the paper.} We start in Section \ref{S:prelim} with the preliminaries in which we set up the notation and obtain the heat kernel bounds (Corollary \ref{Cor:all bounds in one} for the on-diagonal term, and Section \ref{SS:offdiagHK} for the off-diagonal terms). We also obtain approximations of the capacities $q_A$ in Section \ref{SS:AldousBrown}, see in particular Corollary \ref{corapproxq_A}.

Section \ref{S:main} proves that the diameter condition implies Gumbel fluctuations (i.e.\ that \eqref{e:DC} implies \eqref{E:main} in Theorem \ref{T:main}), and discusses the cumulants of the size of the uncovered set. The strategy is provided in Section \ref{SS:strategy}.

Section \ref{S:uncoveredset} proves that the diameter condition implies the decorrelation of the uncovered set (i.e.\ that \eqref{e:DC} implies \eqref{e:tv conv law uncov set} in Theorem \ref{T:uncovered_main}), and  refinements on the structure of the uncovered set. 

Section \ref{S:subcritical cases} studies the cover time and uncovered set of graphs which do not satisfy \eqref{e:DC} and completes the proofs of Theorems \ref{T:main}, \ref{T:uncovered_main}, and \ref{T:uncovered intro t hit av}.

Section \ref{s:examples} contains the construction of an explicit example of graphs that are (SUT) but do not satisfy \eqref{e:DC}.

\medskip \textbf{Acknowledgements.} 
We are grateful to the anonymous referee for their helpful suggestions, which significantly improved the structure and clarity of the paper.
We thank Persi Diaconis, Roberto Imbuzeiro Oliveira, Perla Sousi and Matt Tointon for some useful discussions. In particular we thank Matt Tointon for asking us about the relation between Definition \ref{def:UT} and the definition in \cite{TesseraTointonIsop}. This work started when the first two authors were at the University of Cambridge, during which time they were supported by EPSRC grant  EP/L018896/1. A first version of this paper was finished while the first author was in residence at the Mathematical Sciences Research Institute in Berkeley, California, during
the Spring 2022 semester on \emph{Analysis and Geometry of Random Spaces}, which was supported by the National Science Foundation under Grant No. DMS-1928930.

N.B. gratefully acknowledges the support from the Austrian Science Fund (FWF) grants P33083 on ``Scaling limits in random conformal geometry'', 10.55776/F1002 on ``Discrete random structures: enumeration and scaling limits" and 10.55776/PAT1878824 on ``Random Conformal Fields''.
J.H.'s research is supported by NSERC grants. L. T.'s research was supported by FWF grant P33083, ``Scaling limits in random conformal geometry''.

\section{Preliminaries}\label{S:prelim}

\subsection{A priori bounds on the heat kernel}

From now on, and for the rest of the article, to ease notation, we will write with a small abuse of notation $\Gamma$ also for the vertex set of the graph $\Gamma$. We will denote the size of the vertex set of $\Gamma$ by $|\Gamma|$ or $n$. Our constants may implicitly depend on the degree $d$, which we will not always recall. This also applies to certain notation such as $a_n = o(b_n)$, $a_n = O(b_n)$. We also sometimes use $a_n \ll b_n$ to mean that $a_n = o(b_n)$, i.e.\ that $a_n / b_n \to 0$ and $a_n \lesssim b_n$ to mean that $a_n = O(b_n)$, i.e.\ that $a_n /b_n$ is bounded. We warn the reader 
that the conventions in geometric group theory are a little different:
while using $\ll$ for $o(\cdot)$ is common in probability, it is often used for $O(\cdot)$ in geometric group theory, including in the work of Tessera and Tointon \cite{TesseraTointonIsop,TesseraTointonfinitary}.

It will be at times convenient to allow the graph $\Gamma$ to not be vertex-transitive. Let $P(\cdot,\cdot)$ be the transition kernel of our walk. Define the conductance $\Phi(S)$ of a (non-empty) set $S \subset \Gamma$ by
\begin{equation}
    \Phi(S) = \frac{Q(S, S^c)}{\pi(S)},
\end{equation}
where for $A,B \subset \Gamma$,
\begin{equation}\label{e:Q}
    Q(A,B) = \sum_{a\in A, b\in B} \pi(a) P(a,b).
\end{equation}
The conductance profile is the function $\phi$ defined, for $\min_{x\in\Gamma} \pi(x) \leq u \leq 1/2$, by
\begin{equation}\label{e:conduct_def}
    \phi(u) = \inf \ag \Phi(S) \du \pi(S) \leq u \ad.
\end{equation}
We first recall here a crucial consequence of the theory of evolving sets which allows us to get bounds on the heat kernel given a conductance profile.
Let $\varepsilon>0$. Then we know by \cite[Theorem 13]{MorrisPeres2005}, considering only diagonal transitions, that for every $x\in \Gamma$ and every $\varepsilon$ such that $\pi(x) \leq 1/\varepsilon \leq 1/8$,
\begin{equation}\label{E:isoperimetric time}
t \ge     \intt{4\pi(x)}{4/\varepsilon}\frac{8\mathrm{d}u}{u\phi(u)^2} \implies p_t(x,x) \le (1+ \varepsilon) \pi(x).
\end{equation}
Consider now the vertex-transitive case where $\pi(x) = 1/n$ for every $x\in \Gamma$, where $n=|\Gamma|$. When applying this inequality it is essential to note that the condition on $\varepsilon$ is $8 \le \varepsilon  \le n$. (It may be slightly perverse to name $\varepsilon$ a quantity which by assumption is greater or equal to 8, but this notation is by now relatively well established. Here $w = n/ \varepsilon$ plays at large times a role similar to the volume in the graph, namely we get $p_t(o,o) \asymp 1/w$, and $\varepsilon$ will typically be small compared to the volume $n$.)

Our first task will be to transform the condition on the conductance profile into a condition on the volume growth, or equivalently, on isoperimetry. (Combined with the work of Tessera and Tointon on a quantitative form of Gromov's theorem giving strong control on the volume growth, this will give us excellent control on the return probabilities.)
For $x\in \bbR$, we denote by $\lfloor x \rfloor = \max \ag y \in \bbZ \du y\leq x \ad$ its integer part and by $\ag x \ad = x - \lfloor x \rfloor$ its fractional part.
\report{We first recall for ease of use later on a result of Tessera and Tointon \cite{TesseraTointonIsop}, which we specialise to the case where the radius $r$ (in their notation) is the diameter $D$.}
\report{\begin{lemma}[{\cite[Theorem 7.2]{TesseraTointonIsop}}]\label{lem:tesseratointonisop}
   Let $m\geq 1$ be an integer. There exists a constant $c_0(m)>0$ such that the following holds. Let $q\geq m$ be a real number.  Let $\Gamma$ be a (finite connected) locally finite vertex-transitive graph such that $n = |\Gamma| \geq D^q$. Then for any subset $S\subset \Gamma$ such that $1\leq |S| \leq n/2$, we have
   \begin{equation}
       |\partial S| \geq c_0(m) \min \pg |S|^{\frac{m}{m+1}} , D^{\frac{\ag q \ad}{m}} |S|^{\frac{m-1}{m}}  \pd,
   \end{equation}
   where $\partial S$ is the external vertex boundary of $S$, i.e.\ $\partial S = \ag x \in \Gamma \backslash S \du x\sim y \text{ for some } y\in S \ad$.
\end{lemma}}
\begin{remark}
\report{Tessera and Tointon \cite{TesseraTointonIsop} prove results in various equivalent forms, in terms of volume growth, effective resistances, and isoperimetric inequalities. We note that the connection between isoperimetry and volume growth can also be found in \cite[Lemma 10.46]{LyonsPeresBook}, and for the infinite case in \cite[Lemma 7.2]{LyonsMorrisSchramm2008}.} 
\end{remark}

\begin{prop}\label{P: bounds on return probabilities}
 Let $m\geq 1$ and $d\geq 2$ be integers. There exists a constant $C_0(m,d)>0$ such that the following holds. Let $q\geq m$ be a real number. Let $\Gamma$ be a (finite connected) vertex-transitive graph of degree $d$ such that $n = |\Gamma| \geq D^q$. Then for every $0< t\leq D^2$, we have
    \begin{equation}\label{HKbound_diag}
p_t(o,o) \leq C_0 \max \pg \frac{1}{t^{(m+1)/2}}, \frac{1}{Rt^{m/2}} \pd,
\end{equation}
where $R = D^{\ag q \ad}$.
\end{prop}

\begin{remark}
  For a given $m = \lfloor q \rfloor$ this upper bound is sharp, as is easily shown by considering ``flat'' tori of the form $\Gamma = (\bbZ/ L \bbZ)^m \times (\bbZ / R \bbZ)$ with $1 \ll R \le L$. In this case, the heat kernel initially decays like in dimension $(m+1)$; but eventually (for $t \gtrsim R^2$) the decay becomes only $m$-dimensional.
\end{remark}

\begin{remark}\label{rem:dimension at least 5}
    In the statements of Lemma \ref{lem:tesseratointonisop} and Proposition \ref{P: bounds on return probabilities}, it is not assumed that $\lfloor q \rfloor = m$. In particular, applying  $m=5$ Proposition \ref{P: bounds on return probabilities} with $m=q=5$ shows that uniformly for all finite transitive graphs that satisfy $|\Gamma| \geq D^5$, the bound $p_t(o,o) \leq C_0(5,d)/t^{5/2}$ holds for all $0<t\leq D^2$. 
\end{remark}

\begin{remark}
We will only apply this result with $1\leq m \leq 5$, so we do not need to track the dependence on $m$.
\end{remark}

\begin{remark}\label{R:discrete}
  By using \cite[Theorem 1]{MorrisPeres2005} instead of \cite[Theorem 13]{MorrisPeres2005}, the same estimate holds for the discrete time transition probabilities.
\end{remark}
\begin{proof}
For every $S \subset \Gamma$ such that $0<\pi(S) \leq 1/2$,
we have 
\begin{equation}
    Q(S,S^c) = \sum_{a\in S, b \in \partial S} \pi(a)P(a,b) = \sum_{a\in S, b \in \partial S} \frac{1}{nd} \geq  \frac{|\partial S|}{nd},
\end{equation}
and therefore by Lemma \ref{lem:tesseratointonisop} 
\begin{equation}
    d\cdot \Phi(S) \geq \frac{1}{\pi(S)}\frac{|\partial S|}{n} = \frac{|\partial S|}{\pi(S)} \geq c_0(m) \min \pg |S|^{-1/(m+1)}, R^{1/m} |S|^{-1/m}\pd.
\end{equation}
For $0<v \le n/2$, considering sets $S$ such that $ v / n \leq \pi(S) \leq 1/2$, we get immediately that
\begin{equation} \label{E:conductance profile original chain}
\phi(v/n) \geq \frac{c_0(m)}{d} \min \pg v^{-1/(m+1)}, R^{1/m} v^{-1/m}\pd.
\end{equation}
We deduce, making the change of variables $v = nu$ and setting  $w = 8n/ \varepsilon$ with $ 8\leq \varepsilon <n$,
\begin{align*}
\intt{4/n}{4/\varepsilon}\frac{8\mathrm{d}u}{u\phi(u)^2} = \intt{4}{4n/\varepsilon}\frac{8\mathrm{d}v}{v\phi(v/n)^2} & \lesssim \intt{4}{4n/\varepsilon}  \max\pg v^{2/(m+1)}, \frac{v^{2/m}}{R^{2/m}}\pd \frac{\mathrm{d}v}{v} \\
& \lesssim \max\pg w^{2/(m+1)}, \frac{w^{2/m}}{R^{2/m}}\pd.
\end{align*}
Consequently, there exists $C\geq 1$ such that if $t = C \max\pg w^{2/(m+1)}, \frac{w^{2/m}}{R^{2/m}}\pd$, we may apply \eqref{E:isoperimetric time} and obtain (recall that $w = 8n / \varepsilon$),
\begin{equation*}
p_t(o,o) \leq \frac{1+\varepsilon}{n} = \frac{1}{n} + \frac{8}{w} \leq \frac{9}{w}.
\end{equation*}
As $t = C \max\pg w^{2/(m+1)}, \frac{w^{2/m}}{R^{2/m}}\pd$, we have in particular \\
$w \gtrsim \min\pg t^{(m+1)/2},Rt^{m/2} \pd$, and the bound on $p_t(o,o)$ can be rewritten as
\begin{equation}\label{E: some bound}
p_t(o,o) \lesssim \max\pg \frac{1}{t^{(m+1)/2}},\frac{1}{Rt^{m/2}}  \pd.
\end{equation}
Finally, as we defined first $1<w\leq n$, and then $t = C \max\pg w^{2/(m+1)}, \frac{w^{2/m}}{R^{2/m}}\pd$, \eqref{E: some bound} holds for all $t$ such that $C < t\leq CD^2$, and hence, as we took $C \geq 1$ and up to changing the implicit constant because of the values of $t$ smaller than $C$, for all $t\leq D^2$, as desired.
\end{proof}

To avoid separating the proof into too many cases and cutting the integrals into several pieces, the following corollary will be helpful. It provides upper bounds on the diagonal return probabilities which are true for all $t$. For finite connected vertex-transitive graphs $\Gamma$, we recall that we write $n=|\Gamma|$ and $D = \Diam(\Gamma)$. 
\report{
We also set for graphs $\Gamma$ such that $n\geq D^2$,
\begin{equation}
     f(\Gamma) = \frac{|\Gamma|}{\Diam(\Gamma)^2} = \frac{n}{D^2}.
\end{equation}
\begin{corollary}\label{Cor:all bounds in one} Let $d\geq 2$.
\begin{enumerate}
    \item Uniformly over all finite vertex-transitive graphs of degree less or equal to $d$ such that $n\geq D^2$, we have
\begin{equation}
\label{e:ondiag1}
p_t(o,o) \lesssim \frac{1}{t^{3/2}} + \frac{1}{f(\Gamma)t}, \quad \quad \text{ for all } 1 \leq t \leq D^2,
\end{equation}
and
\begin{equation}\label{e:ondiag2}
p_t(o,o) \lesssim \frac{1}{D^3} + \frac{1}{n}, \quad \quad \text{ for all } t \geq D^2.
\end{equation}
\item Uniformly over all finite vertex-transitive graphs of degree less or equal to $d$  such that $n\geq D^5$, and $t>0$, we have
\begin{equation}
p_t(o,o) \lesssim \frac{1}{t^{5/2}} + \frac{1}{D^5}.
\end{equation}
\end{enumerate}
Moreover, the implicit constants in (a) and (b) depend only on $d$.
\end{corollary}}

\begin{remark}
  The bounds we give are not necessarily optimal if we fix the value of $m = \lfloor q\rfloor$, but have the advantage that they do not depend on $m$ and so can be used regardless. This leads to fewer cases to treat separately further down the proof and so makes the argument more unified.
  \end{remark}

\begin{proof}
Set $q$ such that $n = D^q$ and write $R = D^{\ag q\ad}$. If $2\leq q <3$, i.e.\ $1\leq f(\Gamma) < D$, we have $R = f(\Gamma)$, and applying Proposition \ref{P: bounds on return probabilities} gives \eqref{e:ondiag1}. If $q\geq 3$, i.e.\ $f(\Gamma) \geq D$, applying Proposition \ref{P: bounds on return probabilities} with $m=3$ gives that $p_t(o,o) \lesssim \frac{1}{t^{3/2}}$ uniformly for all $1\leq t\leq D^2$, which implies \eqref{e:ondiag1}.

In all cases, applying $\eqref{e:ondiag1}$ at time $t=D^2$ therefore gives $p_{D^2}(o,o) \lesssim \frac{1}{D^3} + \frac{1}{n}$. Finally, since $t\mapsto p_t(o,o)$ is decreasing, the bound $p_t(o,o) \leq p_{D^2}(o,o) \lesssim \frac{1}{D^3} + \frac{1}{n}$ holds uniformly for all $t\geq D^2$, which proves \eqref{e:ondiag2} and concludes the proof of (a).
(b) follows by Proposition \ref{P: bounds on return probabilities} with $q=5$, proceeding exactly as in the proof of (a).
\end{proof}

\subsection{Hitting times and volume bounds} In this section we show that under \eqref{e:DC} the hitting times of small sets are of order $|\Gamma|$. We also show elementary bounds on the growth of balls in vertex-transitive graphs.

\begin{prop}\label{prop: thit at most order n under DC}
     Let $(\Gamma)$ be a collection of finite (connected) vertex-transitive graphs of fixed degree that satisfies \eqref{e:DC}. Recall that we denote $n=|\Gamma|$. We have 
     \begin{equation}
         \thit \asymp n.
     \end{equation}
\end{prop}
\begin{proof}
   Using transitivity and the commute-time identity (see \cite[Proposition 10.7]{LivreLevinPeres2019MarkovChainsAndMixingTimesSecondEdition}), we have for all $x,y \in \Gamma$
\begin{equation*}
2\E_x( T_y) = \E_x( T_y) + \E_y(T_x) = d\cR({x\leftrightarrow y}) n.
\end{equation*}
Moreover, writing 
$\cR_{\Gamma, 2} := \max_{x,y\in\Gamma} \cR(x \leftrightarrow y)$, and denoting the degree of $\Gamma$ by $d$ we have from \cite[Theorem 2.3]{TesseraTointonIsop}, that
\begin{equation}
    \cR_{\Gamma, 2} \lesssim \frac{1}{d} + \frac{D^2\log n}{n} = \frac{1}{d} + o(1) = O(1).
\end{equation}
This shows that $\thit \leq \frac{d}{2} \cR_{\Gamma, 2} n = O(n)$.

The other direction is simpler: by \cite[Proposition 1.19]{LivreLevinPeres2019MarkovChainsAndMixingTimesSecondEdition}, 
\begin{equation}\label{e:simple general lower bound on the maximal hitting time}
    \thit = \max_{x,y}\bbE_x T_y \geq \max_x \bbE_x T_x^+ -1 = \bbE_o T_o^+ - 1 = n-1.
\end{equation}
This concludes the proof.
\end{proof}

A crucial part of the argument will be to obtain good bounds on the hitting time of a finite arbitrary set $A$ of fixed cardinality $k \ge 1$. Problems usually arise when some points in $A$ are relatively close to one another. The following proposition shows that under \eqref{e:DC}, the expected hitting time of $A$ is always of order $n$.  

\begin{prop}\label{propepitathetan}
 Let
 $(\Gamma)$ be a collection of finite (connected) vertex-transitive graphs of fixed degree $d$, and let $k\geq 1$. Then, as $n=|\Gamma|\to\infty$, uniformly over all subsets $A$ of cardinality $k$ of $\Gamma$,
\begin{equation}
     \bbE_\pi T_A \gtrsim_{d,k} n.
\end{equation}
Moreover, if $(\Gamma)$ satisfies \eqref{e:DC},  as $n=|\Gamma|\to\infty$ we have, uniformly over all subsets $A$ of cardinality $k$ of $\Gamma$,
\begin{equation} 
    \bbE_\pi T_A\asymp_{d,k} n.
\end{equation}
\end{prop}

\begin{proof}
Let $A\subset \Gamma$ of cardinality $k$. Averaging \eqref{E: hitting probability and local time} (exceptionally using an independent result from the next section) with respect to $x$, and taking $y=o$, we have for every $t\geq 0$
\begin{equation}
 \bbP_\pi (T_o \leq t) \leq e\bbE_\pi L_o(t+1) = e(t+1)/n.
\end{equation}
It follows that for every $t \geq 1$, as $t+1 \leq 2t$, and taking
$t =  n/(4ek)$,
\begin{equation*}
\bbP_\pi(T_A \leq t) \leq k\bbP_\pi(T_o\leq t) \leq 2ekt/n\le 1/2.
\end{equation*}
By Markov's inequality, we deduce
\begin{equation*}
\bbE_\pi T_A \geq t \bbP_\pi(T_A > t) = t \pg 1 - \bbP_\pi(T_A \leq t) \pd \geq t/2 \gtrsim n,
\end{equation*}
which concludes the proof of the lower bound.
The upper bound is straightforward, since
 $\bbE_\pi T_A \leq \bbE_\pi T_o \leq \thit$, $\thit = O(n)$ under \eqref{e:DC} by Proposition \ref{prop: thit at most order n under DC}.
\end{proof}

Let us also collect two simple bounds on volumes. We denote the volume of a ball of radius $r$ in a vertex-transitive graph by $V(r)$.

\begin{lemma}\label{lem:volume bounds}
     Let $\Gamma$ be a finite connected vertex-transitive graph of degree $d\geq 2$. Denote $n=|\Gamma|$ and $D=\Diam(\Gamma)$.
\begin{enumerate}
    \item \report{We have $V(r)\leq (d+1)^r$ for any $1\leq r \leq D$. In particular, $n \leq (d+1)^D$, and \begin{equation}\label{eq: trivial bound diameter}
        D\geq \frac{\log n}{\log (d+1)}\gtrsim_d \log n.
    \end{equation}}
    \item For any $1\leq r \leq D$, we have
    \begin{equation}\label{volumegrowth}
V(r) \le \frac{6r}{D}n.
\end{equation}
\end{enumerate}    
\end{lemma}
\begin{proof}
Each vertex has $d$ neighbours. Therefore for $r\geq 0$ we have $V(r) \leq (d+1)^r$. The bound $n\leq (d+1)^D$ then follows, since $V(D) = n$, and is equivalent to $D\geq (\log n)/\log(d+1)$, which proves (a). Let us now prove (b).
Consider two points $x,y\in \Gamma$ such that $d(x,y) = D$. There exists a sequence of vertices $x= x_0, x_1, \ldots, x_D = y$ such that $\ag x_i, x_{i+1}\ad$ is an edge of $\Gamma$ for each $0\leq i\leq D-1$. Let $1\leq r \leq D$. The bound $V(r) \le \frac{6r}{D}n$ is trivial if $r\geq D/6$, so we may assume that $r<D/6$.
By the triangle inequality the $\lfloor D/(3r)\rfloor$ balls of radius $r$ centered at $x_{0}, x_{3r}, x_{2\cdot 3r}, \ldots, x_{\lfloor D/(3r)\rfloor\cdot 3r}$ are disjoint, and therefore $\lfloor D/(3r)\rfloor V(r) \leq n$. Finally, since $r<D/6$, we have $\lfloor D/(3r)\rfloor \geq D/(6r)$, and we conclude that $ V(r) \leq n/\lfloor D/(3r)\rfloor \leq 6nr/D$. 
\end{proof}

\subsection{Relaxation times and mixing times}\label{S:relaxmixhittimes}

Recall that the relaxation time of the chain is the inverse of the spectral gap:
\begin{equation}
    \trel = \frac{1}{1-\lambda_2},
\end{equation}
where $\lambda_2$ is the second largest eigenvalue of the chain.
By a classical argument (see e.g.\ \cite[Theorem 7.6]{ChungBook}) the following bound on the relaxation time is valid on every finite (connected) vertex-transitive graph:
\begin{equation}\label{treldiff}
\trel \leq dD^2.
\end{equation}
We now show that for graphs of polynomial growth, $\trel \asymp D^2$. 
\begin{lemma}\label{lem:trel asymp diam squared}
    Let $m\geq 1$ be an integer. There exists a  constant $c = c(m)$
    such that for every finite (connected) vertex-transitive graph of degree $d$ such that $n = |\Gamma| \leq D^m$,
\begin{equation}
    cD^2 \leq \trel \leq dD^2.
\end{equation}
\end{lemma}
\begin{proof}
We follow arguments from the proof of the lower bound of \cite[Theorem 3.1]{DiaconisSaloff-Coste1993}.

By the minimax characterisation of the spectral gap (see for instance \cite[Lemma 13.7 and Remark 13.8]{LivreLevinPeres2019MarkovChainsAndMixingTimesSecondEdition}) we have, setting $g = d(o,\cdot)$,
\begin{equation}
    1-\lambda_2 \leq \frac{\cE(g)}{\mathrm{Var}_\pi(g)}.
\end{equation}
Since $g$ is 1-Lipschitz, $\cE(g) \leq 1/2$. We hence have
\begin{equation}
    \trel \geq 2 \mathrm{Var}_\pi(g) = \frac{1}{n^2}\sum_{x,y \in \Gamma} (g(x)-g(y))^2.
\end{equation}
Let $z \in \Gamma$ such that $d(o,z) = D$. We deduce that
\begin{equation}
    \trel \geq \frac{1}{n^2} \sum_{x\in B(o,D/4), y \in B(z,D/4)} (D/2)^2 =\frac{D^2}{4} \pg \frac{V(D/4)}{n} \pd^2.
\end{equation}
Moreover, by \cite[Proposition 6.1]{TesseraTointonIsop}, there exists a constant $c_0 = c_0(m)$ such that for all finite (connected) vertex-transitive graphs satisfying $n \leq D^m$,
\begin{equation}
    V(D/4) \geq c_0 n,
\end{equation}
where we recall that $V(r)$ denotes the volume of a ball of radius $r$.
For those graphs, we hence have
\begin{equation} 
    (c_0^2/4)D^2 \leq \trel \leq dD^2. \qedhere
\end{equation}
\end{proof}
Let us now define mixing times and the distance to stationarity. We restrict here to our framework of simple random walks on vertex-transitive graphs, so the walk is in particular transitive and the stationary distribution is the uniform distribution $\pi$. See \cite[Chapter 4]{LivreLevinPeres2019MarkovChainsAndMixingTimesSecondEdition} for more general definitions.
For $1\leq p< \infty$, we define the $L^p$ distance to stationarity of the chain $(X_t)_{t\geq 0}$ at time $t$ by
\begin{equation}
    d^{(p)}(t) =\pg \frac{1}{n} \sum_{x\in\Gamma} \bg np_t(o,x) - 1 \bd^p \pd^{1/p}.
\end{equation}
We also define the $L^\infty$ distance to stationarity as (see \cite[Proposition 4.15]{LivreLevinPeres2019MarkovChainsAndMixingTimesSecondEdition})
\begin{equation}
    d^{(\infty)}(t) = np_t(o,o) - 1.
\end{equation}
We define the $L^p$ mixing time at level $\varepsilon$, for $\varepsilon<0$ and $1\leq p\leq \infty$, by
\begin{equation}
    \tmixp(\varepsilon) = \inf \ag t\geq 0 \du d^{(p)}(t) \leq \varepsilon \ad. 
\end{equation}
Note that the total variation distance is just half of the $L^1$ distance. 
For all $x,y\in \Gamma$ and $t\geq 0$, we set
\begin{equation}
    h_t(x,y) = p_t(x,y) - \frac{1}{n}.
\end{equation}

\begin{prop}\label{P:tmix asymp diam squared}
 Let $m\geq 1$ and $d\geq 2$ be integers. For $0<\varepsilon<1$, we have, uniformly over all finite (connected) vertex-transitive graphs of degree $d$ such that $n = |\Gamma| \leq D^m$,
\begin{equation}
\tmixinfty(\varepsilon) \asymp_{m,d,\varepsilon}\tmixun(\varepsilon) \asymp_{m,d,\varepsilon} \trel \asymp_{m,d} D^2.
\end{equation}
\end{prop}
\begin{proof}
We already proved in Lemma \ref{lem:trel asymp diam squared} that $\trel \asymp_d D^2$. By convexity, if $1\leq p\leq q \leq \infty$, and $t\geq 0$, $d^{(p)}(t)\leq d^{(q)}(t)$, and hence for $0<\varepsilon<1$, we have $\tmixp(\varepsilon) \leq \tmixq(\varepsilon)$. By \cite[Lemma 20.11]{LivreLevinPeres2019MarkovChainsAndMixingTimesSecondEdition}, we hence have $d^{(\infty)}(t) \geq d^{(1)}(t) \geq e^{-t/\trel}$ (since the total variation distance is exactly the half of the $L^1$ distance) and thus for $0<\varepsilon<1$,
\begin{equation}\label{eq:lower bound mixing time}
    \tmixinfty(\varepsilon) \geq \tmixun(\varepsilon) \geq \log(1/\varepsilon)\trel.
\end{equation}
Let us now prove that $\tmixinfty(\varepsilon) \lesssim \trel$.
By spectral estimates, we have for all $t,s\geq 0$, 
\begin{equation}\label{Eq:spec}
    h_{t+s}(o,o) \leq e^{-s/\trel}h_t(o,o).
\end{equation}
Moreover, by Proposition \ref{P: bounds on return probabilities}, there is a constant $C(m,d)$ such that uniformly over all finite connected vertex-transitive graphs of degree $d$ such that $n\leq D^m$, $h_{D^2}(o,o) \leq C(m,d)$.
Hence, at time $D^2 + i\trel$, where $i = \left \lceil \log(C(m,d)/\varepsilon )\right \rceil $,
\begin{equation}
    h_{D^2 + i\trel}(o,o) \leq e^{-i}h_{D^2}(o,o) \leq \varepsilon,
\end{equation}
so, recalling that $\trel \leq dD^2$,
\begin{equation}\label{eq:upper bound mixing time}
    \tmixinfty(\varepsilon) \leq D^2 + \left \lceil \log(C(m,d)/\varepsilon )\right \rceil\trel \leq D^2\pg 1 + d\left \lceil \log(C(m, d)/\varepsilon )\right \rceil\pd,
\end{equation}
concluding the proof.
\end{proof}

\begin{remark}
Proposition \ref{P:tmix asymp diam squared} shows that when $n\leq D^m$, the relaxation time and the mixing time are both of order $D^2$. The explicit bounds \eqref{eq:lower bound mixing time} and \eqref{eq:upper bound mixing time} show in particular that there is no \textit{cutoff}, i.e.\ that for some $0<\varepsilon<1$, $\tfrac{\tmixinfty(1-\varepsilon)}{\tmixinfty(\varepsilon)}$ does not converge to 1 as $n = |\Gamma| \to \infty$. See \cite[Chapter 18]{LivreLevinPeres2019MarkovChainsAndMixingTimesSecondEdition} for more details on the cutoff phenomenon.
\end{remark}

The following proposition is a classical result which can be traced back to \cite{Aldous1982someinequalities}.

\begin{prop}\label{prop:comparison average and max hitting times} Uniformly over all finite (connected) vertex-transitive graphs we have 
\begin{equation}
    \thit - \thitav \asymp \tmixun\pg (2e)^{-1} \pd.
\end{equation}
\end{prop}
\begin{proof} 
Let $\tau_i$ for $1\leq i\leq 4$ be as in \cite{Aldous1982someinequalities}, and set $\tau_5 := \thit - \thitav$. From \cite[Equation (17)]{Aldous1982someinequalities}, we have $\tau_5 \leq \tau_2$. On the other hand, we have 
\begin{equation}
    \tau_3 = \max_{x,y\in \Gamma}\sum_{z\in \Gamma}\pi(z)\bg \bbE_x T_z - \bbE_y T_z \bd = 2\max_{x,y\in \Gamma}\sum_{z\in \Gamma \du \bbE_x T_z \geq \bbE_y T_z}\pi(z)\pg \bbE_x T_z - \bbE_y T_z \pd.
\end{equation}
Moreover, for any $x,y,z\in \Gamma$, we have $\bbE_x T_z - \bbE_y T_z \leq \thit - \bbE_y T_z$, and (by definition of $\thit$) $\thit - \bbE_y T_z \geq 0$. We deduce from this that
\begin{equation}
    \tau_3 \leq 2 \max_{x,y\in \Gamma} \sum_{z\in\Gamma} \pi(z)(\thit - \bbE_y T_z) = 2\tau_5,
\end{equation}
so we have proved that $\tau_3/2 \leq \tau_5 \leq \tau_2$. We conclude using \cite[Theorem 5]{Aldous1982someinequalities} that
\begin{equation}\label{eq: prop:comparison average and max hitting times eq inter}
    \tau_5 \asymp \tau_1 = \tfrac{1}{2}\tmixun\pg (2e)^{-1} \pd. \qedhere
\end{equation}
\end{proof}
\begin{remark}\label{remark: QS will be used for comparison average and max hitting times}
Proposition \ref{prop:comparison QS and max hitting times} below proves a variant of Proposition \ref{prop:comparison average and max hitting times}, where $\thitav = \bbE_\pi T_o$ is replaced by the quasi-stationary distribution associated with the set $A = \{o\}$ (as defined below in Section \ref{SS:AldousBrown}).
\end{remark}

\begin{prop}\label{P:thitmaxav}
Let $(\Gamma)$ be a collection of finite (connected) vertex-transitive graphs of fixed degree that satisfies  \eqref{e:DC}. As $n= |\Gamma| \to \infty$ we have
\begin{equation}
    \thit = \thitav(1+o(1/\log n)).
\end{equation}
In particular, for every $s \in \bbR$, we have as $n=|\Gamma| \to \infty$,
\begin{equation}
    \thit((\log n) + s) = \thitav((\log n) + s + o(1)).
\end{equation}
\end{prop}
\begin{proof}
It is enough to prove that $\thit = \thitav(1+o(1/\log n))$, i.e.\ that $\frac{\thit}{\thitav} - 1 =o(1/\log n)$. Since $\thitav \leq \thit$ by definition, we only need to upper bound $\frac{\thit}{\thitav} - 1$. By Proposition \ref{propepitathetan} we have $\thitav = \bbE_\pi T_o \asymp n$. We deduce from Proposition \ref{prop:comparison average and max hitting times} that
\begin{equation}
    \frac{\thit}{\thitav} - 1\asymp \frac{\tmixun((2e)^{-1})}{n}.
\end{equation}
\report{If $n\leq D^5$ then $\tmixun((2e)^{-1}) \asymp D^2$ by Proposition \ref{P:tmix asymp diam squared} so 
\begin{equation}
    \frac{\tmixun((2e)^{-1})}{n} \asymp \frac{D^2}{n},
\end{equation}
and the result holds by the assumption \eqref{e:DC}. On the other hand, if $n> D^5$, then by \cite[Theorem 12.4]{LivreLevinPeres2019MarkovChainsAndMixingTimesSecondEdition}, recalling \eqref{treldiff}, and since $D\gtrsim \log n$ by Lemma \ref{lem:volume bounds} (a), we have
\begin{equation}
     \frac{\tmixun((2e)^{-1})}{n} \lesssim \frac{\trel \log n}{D^5} \lesssim \frac{\log n}{D^3} \lesssim \frac{1}{(\log n)^2} = o\pg \frac{1}{\log n} \pd. \qedhere
\end{equation}}
\end{proof}

\subsection{Off-diagonal heat kernel bounds}
\label{SS:offdiagHK}

Our first task will be to translate the on-diagonal bounds described in the previous section into off-diagonal bounds. A general upper bound can always be obtained with the Carne--Varopoulos inequality, but we will need a sharper recent result due to Folz, see \cite[Corollary 1.2]{Folz2011}.
This general result takes as an input a time-dependent bound on the on-diagonal heat kernel and deduces from this a general off-diagonal bound. The result is actually quite general but we will only use it in the most simple case where the Lipschitz function is taken to be the graph distance, and the ``volume'' is measured with respect to cardinality.

Set $1/g(t)$ to be the right hand side of \eqref{e:ondiag1}, so that $1/g(t)$ upper bounds $p_t(o,o)$; that is,
\begin{equation}\label{g}
1/ g(t)  =
\begin{cases}
\frac{C}{t^{3/2}}
& \text{ for } 1\leq t \leq f(\Gamma)^2 \\
 \frac{C}{t f(\Gamma)} & \text{ for } f(\Gamma)^2 \leq t \leq D^2.
 \end{cases}
 \end{equation}
(Recall that if $m =2$ then $R = f(\Gamma)$.) A consequence of his result is the following inequality. Assume that $f(\Gamma) < D$. Then there exists a constant $c>0$ such that for every $x \neq y\in \Gamma$ and $d(x,y)\leq t \leq D^2$,
\begin{equation}\label{Folz}
    p_t(x,y) \lesssim \frac{1}{g(t)}\exp\pg-c\frac{d(x,y)^2}{t}\pd,
\end{equation}
where the implicit constant in $\lesssim$ depends only on the degree bound.
For smaller values of $t$, we will simply bound the heat kernel via the Carne--Varopoulos bound, \begin{equation*}
p_t(x,y) \lesssim  \exp ( - c d(x,y)^2 / t) \le \exp ( - c d(x,y))
\end{equation*}
when $t \le d(x,y)$.

We will use \eqref{Folz} to get upper bounds on the off-diagonal heat kernel, and therefore by integrating, on the expected local time at a given vertex. In turn, this can be used to upper bound the probability to visit a vertex $y$ far away from $x$ in a relatively short time, via the following elementary lemma.

\begin{lemma}\label{L:unionbound}
{Let $\Gamma$ be a finite connected graph,} and let $x,y\in \Gamma$. For the (rate-1) continuous time simple random walk on $\Gamma$, we have
\begin{equation}\label{E: hitting probability and local time}
\bbE_x L_y(t+1) \geq \frac{1}{e} \bbP_x \pg T_y \leq t \pd.
\end{equation}
\end{lemma}
\begin{proof}
Let us define the event $E$ by
\begin{equation}
E := \ag \text{the walk stays for time at least 1 at } y \text{ just after } T_y \ad.
\end{equation}
Then, we have
\begin{equation*}
\bbE_x L_y(t+1) \geq \bbE_x \pg L_y(t+1) | T_y \leq t, E \pd \bbP_x \pg T_y \leq t, E \pd \geq \bbP_x \pg T_y \leq t \pd \bbP(E) \geq \bbP_x \pg T_y \leq t \pd e^{-1}.
\end{equation*}
This proves the lemma.
\end{proof}

We now combine the above ideas to get the following bounds:

\begin{prop}\label{P:bound on transition probabilities before some time}  Let
$(\Gamma)$ be a collection of finite (connected) vertex-transitive graphs of fixed degree. Recall that we denote $n=|\Gamma|$, $D = \Diam(\Gamma)$, and $f(\Gamma) = n/D^2$.
\begin{enumerate}
    \item Assume that $n \geq D^2$. Then, uniformly over all $\delta \in [1,D/2]$, all $x,y\in\Gamma$ such that $d(x,y) \geq \delta$, and all $t\geq 1$, we have 
    \begin{equation}
\bbE_x L_y(t) \lesssim \frac{1}{{\delta}} + \frac{\log (D/\delta)}{f(\Gamma)} + \frac{t}{n} + \frac{t}{D^3}.
\end{equation}
\item Assume that $n  \geq D^5$. Then, uniformly over all $\delta \in [1,D/2]$, all $x,y\in\Gamma$ such that $d(x,y) \geq \delta$, and all $t\geq 1$, we have \begin{equation}
\bbE_x L_y(t) \lesssim \frac{1}{\delta^{3}} + \frac{t}{D^5}.
\end{equation}
\end{enumerate}
Note that by Lemma \ref{L:unionbound}, the same bounds hold also for $\bbP_x(T_y \leq t)$.
\end{prop}
\begin{proof}
Let $x,y \in \Gamma$.
Since for all $s$, $p_s(x,y) \leq p_s(o,o)$ (see \cite[Lemma 20, in particular Equation (3.60)]{AldousFill})
\begin{equation}
\begin{split}
\bbE_x L_y(t) & =  \intt{0}{t}p_s(x,y)\mathrm{d}s \\ & \leq  \int_0^{d(x,y)} p_s(x,y) \mathrm{d}s  + \intt{d(x,y)}{d(x,y)^2}p_s(x,y)\mathrm{d}s + \mathrm{1}_{t\geq d(x,y)^2}\intt{d(x,y)^2}{t}p_s(o,o)\mathrm{d}s.    
\end{split}
\end{equation}
The first integral is by the Carne--Varopoulos bound smaller or equal to $\delta \exp( - c \delta) \lesssim 1/\delta$. Let us consider the second integral. By \eqref{Folz},
we have
\begin{align*}
\intt{d(x,y)}{d(x,y)^2}p_s(x,y)\mathrm{d}s & \lesssim \int_0^{d(x,y)^2} \left( s^{-3/2} + \frac1{f(\Gamma)s}\right) \exp (- c \delta^2/s) \mathrm{d}s.
\end{align*}
Now, studying the function $t \mapsto t^{-b} \exp (  -A/t)$ we see that this is maximised at $t = A/b$ and so is always $< (A/b)^{-b}$. Applying this with $b = 1$ and $b = 3/2$ as well as $A \asymp \delta^2$, we get
\begin{equation*}
\intt{0}{d(x,y)^2}p_s(x,y)\mathrm{d}s \lesssim \delta^2 \frac1{\delta^3} + \delta^2\frac{1}{ \delta^2 f(\Gamma)} \le \frac1{\delta} + \frac1{f(\Gamma)}.
\end{equation*}
For the third integral, using Corollary \ref{Cor:all bounds in one}, (a), we immediately get
\begin{align*}
\intt{d(x,y)^2}{t}p_s(o,o)\mathrm{d}s & \lesssim \intt{\delta^2}{D^2}\pg \frac{1}{s^{3/2}} + \frac{1}{f(\Gamma)s} \pd \mathrm{d}s +  \mathrm{1}_{t\geq D^2}\intt{D^2}{t}\pg \frac{1}{D^3} + \frac{1}{n} \pd \mathrm{d}s \\
& \lesssim \frac{1}{{\delta}} + \frac{\log (D/\delta)}{f(\Gamma)} + \frac{t}{n} + \frac{t}{D^3}.
\end{align*}
We finally have, as $1 \le \delta \le D/2$,
\begin{equation}\label{e:preciseboundlocaltimexy}
\bbE_x L_y(t) \lesssim \frac{1}{{\delta}} + \frac{\log (D/\delta)}{f(\Gamma)} + \frac{t}{n} + \frac{t}{D^3} + \frac{1}{f(\Gamma)}
\lesssim \frac{1}{{\delta}} + \frac{\log (D/\delta)}{f(\Gamma)} + \frac{t}{n} + \frac{t}{D^3},
\end{equation}
which proves (a). The second bound is proved the same way, using part (b) of Corollary \ref{Cor:all bounds in one} and taking $1/g(t) = C/t^{5/2}$ for $1\leq t \leq D^2$.
\end{proof}

\begin{corollary}\label{proppxygrand} Let
$\cF = (\Gamma)$ be a collection of finite (connected) vertex-transitive graphs of fixed degree satisfying \eqref{e:DC}. Let $(t_\Gamma)_{\Gamma \in \cF}$ and $(\omega_\Gamma)_{\Gamma \in \cF}$ such that 
$t_\Gamma =  o\pg{\min(|\Gamma|, \Diam(\Gamma)^3)}\pd$ and $\omega_\Gamma\to \infty$ as $|\Gamma| \to \infty$. Then, as $|\Gamma| \to \infty$, uniformly over $x,y\in \Gamma$ with $d(x,y) \ge \omega_\Gamma$,
\begin{equation}
    \bbE_x L_y(t_\Gamma) \to 0.
\end{equation}
In particular,
\begin{equation} \label{E:hitbeforen}
    \bbP_x(T_y \leq t_\Gamma) \to 0.
\end{equation}
\end{corollary}
\begin{proof}
It follows immediately from part (a) of Proposition \ref{P:bound on transition probabilities before some time}, Lemma \ref{L:unionbound}, and \eqref{e:DC}.
\end{proof}

When we consider the diagonal case we get a similar bound, but this time of order 1; in this case we can therefore afford to consider times that are as big as the volume. Such a bound is also a signature of our local transience condition. As mentioned in the introduction, Tessera and Tointon \cite{TesseraTointonIsop} proved \eqref{e:WUT2}. We present its proof for the sake of completeness. 

\begin{prop}\label{P: finite time spent at the origin before time n} Let $\cF = (\Gamma)$ be a collection of finite (connected) vertex-transitive graphs of fixed degree.
If $(\Gamma)$ satisfies $|\Gamma| \gtrsim \Diam(\Gamma)^2\log |\Gamma|$ as $|\Gamma| \to \infty$, then as $|\Gamma| \to \infty$, at time $t=|\Gamma|$ we have
\begin{equation}
\label{e:WUT2}
\bbE_o L_o(|\Gamma|) = O(1).
\end{equation}
If $(\Gamma)$ satisfies \eqref{e:DC}, then
\begin{equation}
\label{e:SUT2}
\lim_{s \to \infty} \limsup_{|\Gamma| \to \infty}  \E_o ( L_o(t_{\mathrm{mix}}(1/4))- L_o(s)) =0.
\end{equation}
\end{prop}

\begin{proof}
Let $\Gamma \in \cF$. Denote as usual $n=|\Gamma|$, $D = \Diam(\Gamma)$, and $f(\Gamma) = n/D^2$.
We first prove \eqref{e:WUT2}. We start as in the proof of Proposition \ref{P:bound on transition probabilities before some time}, cutting the integral at time 1 instead of time $d(x,y)^2$. At time $t = D^2$, we have
\begin{equation}
\begin{split}
\bbE_o L_o(D^2)  = \intt{0}{1}p_s(o,o)\mathrm{d}s + \intt{1}{D^2}p_s(o,o)\mathrm{d}s \lesssim 1 + \left(1 + \frac{\log D}{f(\Gamma)} + \frac{t}{n} + \frac{t}{D^3}\right) = O(1).    
\end{split}
\end{equation}
Recall from \eqref{treldiff} that $\trel \leq dD^2$.
Furthermore, we have, for every $t\geq D^2$, (see for instance \cite[Lemma 4.18 and Equation 20.18]{LivreLevinPeres2019MarkovChainsAndMixingTimesSecondEdition})
\begin{equation}\label{E:Poincaré}
p_t(o,o) -\frac{1}{n} \leq \exp\pg-\frac{t-D^2}{dD^2}\pd p_{D^2}(o,o).
\end{equation}
Hence,
\begin{equation}
\intt{D^2}{n}p_s(o,o)\mathrm{d}s = \intt{D^2}{n}\pg p_s(o,o) - \frac{1}{n}\pd \mathrm{d}s + \intt{D^2}{n}\frac{1}{n}\mathrm{d}s
 \lesssim \frac{1}{n}\intt{D^2}{\infty} \exp\pg \frac{-t}{dD^2} \pd \mathrm{d}s +1  = O(1),
\end{equation}
which concludes the proof of \eqref{e:WUT2}. We now prove \eqref{e:SUT2}. By Corollary \ref{Cor:all bounds in one} we have that $t_{\mathrm{mix}}(1/4) = o(n)$. The remainder of the proof is analogous to that of \eqref{e:WUT2} and hence omitted.
\end{proof}

With this, it is easy to deduce that any point (distinct from the starting point) can be avoided with positive probability up to a time slightly smaller than the volume.

\begin{prop} \label{proppxysmall}
Let $\cF = (\Gamma)$ be a collection of finite (connected) vertex-transitive graphs of fixed degree satisfying \eqref{e:DC}, and
let $(t_\Gamma)_{\Gamma\in \cF}$ such that $t_\Gamma = o(|\Gamma|)$ as $|\Gamma|\to \infty$.
There exists a constant $\beta>0$ such that as $|\Gamma|\to \infty$,  uniformly over all distinct $x,y \in \Gamma$, we have
\begin{equation}
    \bbP_x(T_y>t_\Gamma) \geq \beta +o(1).
\end{equation}
\end{prop}

\begin{proof}
This can be proved using \cite[Theorem 1.8]{TesseraTointonIsop}, but we give a simple, self-contained argument. We argue by contradiction and assume (without loss of generality, perhaps taking a subsequence if needed) that there exist $x,y\in \Gamma$ such that $\bbP_x(T_y\leq t_\Gamma) = 1-o(1)$. As the law of $T_x$ starting from $y$ is the same as the law of $T_y$ starting from $x$ by symmetry (see the proof of \cite[Proposition 2]{Aldous1989}), this implies, using the Markov property, that we also have
\begin{equation}
\bbP_x(T_x^+ \leq 2t_\Gamma) \geq \bbP_x(T_y\leq t_\Gamma)\bbP_y(T_x\leq t_\Gamma) =   1-o(1).
\end{equation}
By the same argument, for every integer $m \geq 1$,
\begin{equation}
    \bbP_x(\text{at least } m \text{ returns to } x \text{ before time } 2mt_\Gamma) \geq \bbP_x(T_x^+ \leq 2t_\Gamma)^m  = 1-o(1).
\end{equation}
Since $t_\Gamma = o(n)$, we also have $2m t_\Gamma = o(n)$ and hence we deduce that
$\E_o ( L_o(n)) \ge m (1- o(1))$ and so is unbounded. This contradicts Proposition \ref{P: finite time spent at the origin before time n}.
\end{proof}

\subsection{Exponential approximation of stationary hitting times}

\label{SS:AldousBrown}

The aim of this subsection is to approximate
the hitting time of a fixed set $A$ of cardinality $k \ge 1$ (starting from stationarity) by an exponential random variable with mean $\bbE_\pi T_A$. The key tool to do this will be to consider the quasi-stationary distribution $\alpha_A$, for which the hitting distribution of $A$ is \emph{exactly} exponential, that is, for $t\geq 0$,
\begin{equation}
\label{e:QS tail}
    \bbP_{\alpha_A}(T_A>t) = \exp\pg - \frac{t}{\bbE_{\alpha_A}\cg T_A \cd}\pd \, .
\end{equation} 
Such an idea goes back to the work of Aldous and Brown \cite{AldousBrown1992}.
The quasi-stationary distribution of $A$ is a distribution whose support is contained in $B:=A^c$, and of support $B$ when $\mathbb{P}_a[T_b<T_A] $ for all $a,b \in B$. 

In fact, it is in principle possible for the set $A$ to disconnect $\Gamma \setminus A$ into several connected components, in which case there are multiple quasi-stationary distributions.
However, under our diameter assumption (in fact, as soon as $n\geq D^2$, say), only one of these components may be of macroscopic size as $n \to \infty$ while $k \ge 1$ is fixed: for a fixed $k$ and all sufficiently large $n$, the rest of the components are of bounded sizes by isoperimetry (see \cite[Theorem 7.1]{TesseraTointonIsop}), where the bound on their sizes depends only on $k$ and the degree.
In such a case we therefore define $\alpha_A$ to be the quasi-stationary distribution associated with that macroscopic component, and we do so without further commenting on this case.

In a companion work \cite{BerestyckiHermonTeyssier2025AB} we develop tools to analyse this quasi-stationary distribution, which allow us to obtain explicit error bounds when we approximate the law of the hitting time of $A$ by this exponential distribution. The error terms improve on those obtained by Aldous and Brown \cite{AldousBrown1992}, and are particularly explicit in the case of vertex-transitive graphs. In this case, the error bound depends on a quantity $\beta(\Gamma)$ which is defined as follows. Suppose $\Gamma$ is a finite vertex-transitive graph. We can write its size in a unique way as $n = D^qR$, where $D = \Diam(\Gamma)$ is the diameter of $\Gamma$, $q$ is an integer and $1\leq R<D$. We then set
\begin{equation}\label{e:def of beta Gamma}
    \beta(\Gamma) := \begin{cases}
\frac{D^4}{(\mathbb{E}_{\pi}[T_o])^{2}} & \text{ if } q \in \ag 1,2\ad \\
\frac{D^4}{(\mathbb{E}_{\pi}[T_o])^{2}}\pg 1 + \frac{R\log R}{D}\pd & \text{ if } q=3 \\
\frac{D^4}{(\mathbb{E}_{\pi}[T_o])^{2}}\pg R + \log (\frac{D}{R})\pd & \text{ if } q=4 \\
1/|\Gamma| & \text{ if } q \ge 5  \\
\end{cases}.
\end{equation}
As we will see below (see, e.g., Theorem \ref{thm: improvement of AB for finite sets under DC}), $\beta(\Gamma)$ is a parameter which allows us to bound the difference between the quasi-stationary distribution $\alpha_A$ and the uniform distribution $\pi$. We also set, if $\Gamma$ is a vertex-transitive graph such that $n = |\Gamma| \geq D^2\log n$, (and with $|\Gamma| \geq 3$), recalling that $f(\Gamma) = n/D^2$,
\begin{equation}\label{eq: definition of beta Gamma}
   b(\Gamma) = \frac{f(\Gamma)}{\log n} = \frac{n}{D^2\log n},
\end{equation}
and
\begin{equation}\label{eq: definition of t Gamma star}
   t_\Gamma^* = D^2 \sqrt{1+ \log b(\Gamma)}.
\end{equation}

\begin{remark}\label{rem: t star lesssim D power 5/2}
By Lemma \ref{lem:volume bounds} (a), for vertex-transitive graphs of a given degree, if the number of vertices diverges then the diameter must also diverge. Moreover, if $n\geq D^2\log n$, \eqref{eq: trivial bound diameter} implies that  $\log b(\Gamma)\leq \log n \lesssim D$, and we therefore have
\begin{equation}
    t_\Gamma^* \lesssim D^{5/2}.
\end{equation}
\end{remark}

The following bound will be crucial.
\begin{lemma}\label{lem: bound on beta Gamma biased by time}
    Let $(\Gamma)$ be a collection of finite (connected) vertex-transitive graphs of fixed degree satisfying \eqref{e:DC}. As $n=|\Gamma|\to \infty$, we have
\begin{equation}
\max\pg \frac{t_\Gamma^*}{n}, \frac{n}{t_\Gamma^*}\beta(\Gamma)\pd = o\pg \frac{1}{\log n} \pd.
\end{equation}
and in particular,
\begin{equation}
        \beta(\Gamma) = o\pg \frac{1}{(\log n)^2} \pd.
\end{equation}
\end{lemma}
\begin{proof}
Set $\varepsilon_\Gamma := \frac{\log n}{f(\Gamma)}$. By \eqref{e:DC} we have $\varepsilon_\Gamma\to 0$. First we have
\begin{equation}
    \frac{t_\Gamma^*}{n} =\frac{\varepsilon_\Gamma \sqrt{\log 1/\varepsilon_\Gamma}}{\log n}  = \frac{o(1)}{\log n} = o\pg \frac{1}{\log n} \pd.
\end{equation}
Let us now bound $\frac{n}{t_\Gamma^*}\beta(\Gamma)$ as a function of the diameter. Since $t_\Gamma^*\gtrsim D^2$, it is enough to bound $\frac{n}{D^2}\beta(\Gamma)$. Recall also from Proposition \ref{propepitathetan} that under \eqref{e:DC} we have $\bbE_\pi T_o \asymp n$.

If $D^2 \log n \ll n<D^3$ then $\frac{n}{D^2}\beta(\Gamma) \asymp \frac{D^2}{n} = \frac{1}{f(\Gamma)} = o\pg \frac{1}{\log n} \pd$.

If $D^3\leq n <D^4$ then $\frac{n}{D^2}\beta(\Gamma) \lesssim \frac{D^2}{n}\log D \lesssim \frac{\log D}{D} = o\pg \frac{1}{\log n} \pd$.

If $D^4\leq n <D^5$ then $\frac{n}{D^2}\beta(\Gamma) \lesssim \frac{D^2}{n} D \leq 1/D = o\pg \frac{1}{\log n} \pd$.

Finally, if $n\geq D^5$ then $\frac{n}{D^2}\beta(\Gamma) = 1/D^2 \lesssim \pg\frac{1}{(\log n)^2}\pd = o\pg \frac{1}{\log n} \pd$. (Here we used that since $n\leq (d+1)^{D}$, where $d$ is the degree of the graphs, we have $D\gtrsim_d \log n$.)

The second display follows from the first one, since $\beta(\Gamma) = \frac{t_\Gamma^*}{n} \cdot \frac{n}{t_\Gamma^*}\beta(\Gamma) \leq \max\pg \frac{t_\Gamma^*}{n}, \frac{n}{t_\Gamma^*}\beta(\Gamma)\pd^2$.
\end{proof}

By applying \cite[Theorem 1.4]{BerestyckiHermonTeyssier2025AB} 
to the particular case of sets $A$ of a given size, we obtain the following result.

\begin{thm}
\label{thm: improvement of AB for finite sets under DC}
      Let $d\geq 2$ and $k\geq 1$. Let $(\Gamma)$ be a collection of finite (connected) vertex-transitive graphs of degree $d$ satisfying \eqref{e:DC}. As $|\Gamma| \to \infty$, uniformly over all non-empty subsets $A\subset \Gamma$ of cardinality at most $k$, we have, denoting $n= |\Gamma|$,
\begin{equation}\label{eq: improvement of AB for finite sets under DC probabilities}
\frac{\bbP_\pi(T_A>t)}{\bbP_{\alpha_{A}}(T_A>t)} = 1 + O( \beta(\Gamma)) = 1 + o\pg \frac{1}{(\log n)^2} \pd
\end{equation}
uniformly over all $t\geq 0$, and 
\begin{equation}\label{eq: improvement of AB for finite sets under DC expectations}
    \frac{\bbE_\pi \cg T_A \cd}{\bbE_{\alpha_{A}} \cg T_A \cd} = 1 + O( \beta(\Gamma)) = 1 + o\pg \frac{1}{(\log n)^2} \pd.
\end{equation}
\end{thm}
\begin{proof}
    In \eqref{eq: improvement of AB for finite sets under DC probabilities} and \eqref{eq: improvement of AB for finite sets under DC expectations}, the first equality is precisely \cite[Theorem 1.4]{BerestyckiHermonTeyssier2025AB}, with the explicit dependence on $k$ and $d$ absorbed into the $O(\cdot)$ term, and the second follows from Lemma \ref{lem: bound on beta Gamma biased by time}.
\end{proof}

\begin{remark}
    
\end{remark}

We recall from the introduction that $\tsav = \bbE_\pi T_o (\log(n) + s)$ for $s \in \bbR$. Recall also that for $\eset \ne A\subset \Gamma$ we write
\begin{equation*}
    q_A = \frac{\bbE_\pi T_o}{\bbE_\pi T_A}.
\end{equation*}
Theorem \ref{thm: improvement of AB for finite sets under DC} enables proving the following exponential approximation of $\bbP_\pi(T_A > \tsav)$, which will later be useful to rewrite the moments of the uncovered set at time $\tsav$.

\begin{prop}\label{propapproxpita} 
Let $d\geq 2$ and $k\geq 1$. Let $\cF = (\Gamma)$ be a collection of finite (connected) vertex-transitive graphs of degree $d$ satisfying \eqref{e:DC}. Let $s\in \bbR$. Let $(\varepsilon_\Gamma)_{\Gamma \in \cF}\in \bbR^\cF$ such that $\varepsilon_\Gamma = o(1/\log |\Gamma|)$ as $|\Gamma|\to \infty$. As $|\Gamma| \to \infty$, uniformly over all non-empty subsets $A\subset \Gamma$ of cardinality at most $k$, we have, denoting $n=|\Gamma|$,
\begin{equation}
    \bbP_\pi\pg T_A > \tsav(1+\varepsilon_\Gamma) \pd = n^{-q_A}e^{-q_A s} (1+o(1)).
\end{equation}
\end{prop}
\begin{proof}
Let $\Gamma \in \cF$, let $\eset \ne A\subset \Gamma$ of cardinality at most $k$, and denote $\alpha = \alpha_A$. 

By Theorem \ref{thm: improvement of AB for finite sets under DC}, and recalling that $\tsav = \bbE_\pi T_o (\log(n) + s)$, we have
\begin{equation}
    \frac{\tsav(1+\varepsilon_\Gamma)}{\bbE_\alpha T_A} =  \frac{\tsav(1+\varepsilon_\Gamma)}{\bbE_\pi T_A}(1+o(1/(\log n)^2)) = q_A((\log n)+s+o(1)).
\end{equation}
We deduce from the definition of quasi-stationarity (recall \eqref{e:QS tail}) that
\begin{equation}
    \bbP_\alpha\pg T_A > \tsav(1+\varepsilon_\Gamma) \pd = \exp\pg -\frac{\tsav(1+\varepsilon_\Gamma)}{\bbE_\alpha T_A}\pd =  n^{-q_A} e^{-q_A(s+o(1))}.
\end{equation}
Since $\frac{\bbE_\pi T_o}{\bbE_\pi T_A} = \Theta(1)$ by Proposition \ref{propepitathetan}, this leads to
\begin{equation}
   \bbP_\alpha\pg T_A > \tsav(1+\varepsilon_\Gamma) \pd =  n^{-q_A} e^{-q_As}(1+o(1)).
\end{equation}
Finally, by Theorem \ref{thm: improvement of AB for finite sets under DC} we have $ \bbP_\pi\pg T_A > \tsav(1+\varepsilon_\Gamma) \pd =  \bbP_\alpha\pg T_A > \tsav(1+\varepsilon_\Gamma) \pd(1+o(1))$. This concludes the proof.
\end{proof}

Using Taylor expansions at 0 of the exponential function, we obtain a relation between $\bbP_\pi\pg T_A > t \pd$ and $\bbE_\pi T_A$ for not too large values of $t$. As before and after, asymptotic notation such as $o(\cdot)$ and $O(\cdot)$ may depend on fixed parameters $k$ and $d$, which we do not always emphasise.

\begin{prop}\label{P:linear}
Let $d\geq 2$ and $k\geq 1$. Let $\cF = (\Gamma)$ be a collection of finite (connected) vertex-transitive graphs of degree $d$ satisfying \eqref{e:DC}. 
Then as $|\Gamma| \to \infty$, denoting $n=|\Gamma|$, we have uniformly over all non-empty subsets $A\subset \Gamma$ of cardinality at most $k$ and all $1\leq t\leq n$, 
\begin{equation}
\bbE_\pi T_A = \frac{t}{\bbPpi{T_A \leq t}}\pg 1+  O\left( \frac{t}{n}\right) + O\left( \frac{n}{t}\beta(\Gamma)\right) \pd.
\end{equation}
\end{prop}
\begin{proof}
Let $\Gamma\in \cF$, $\eset \ne A\subset \Gamma$ of cardinality at most $k$, and  $1\leq t \leq n$. Denote $\alpha = \alpha_A$ and set $x = t/\bbE_\alpha T_A$. By Theorem \ref{thm: improvement of AB for finite sets under DC} and expanding the exponential around 0 (we use that $\bbE_\alpha T_A \geq \bbE_\pi T_A \gtrsim n$, so that $x\lesssim 1$), we have
\begin{align*}
\P_\pi (T_A > t) = e^{-x}\pg 1+ O(\beta(\Gamma))\pd
&=\left( 1- x + O \left( x^2 \right) \right)\left(1+ O(\beta(\Gamma))\right)\\
& = 1 - x + O\left( x^2\right) +O(\beta(\Gamma)).
\end{align*}
It follows recalling again that $\E_\alpha (T_A) \asymp n$, that 
\begin{equation*}
\begin{split}
     \P_\pi (T_A \le t) = x + O\left( x^2\right) +O(\beta(\Gamma)) & = x \pg 1+  O\left( x\right) + O\left(x^{-1}\beta(\Gamma)\right) \pd \\
     & = \frac{t}{\bbE_\alpha T_A} \pg 1+  O\left( \frac{t}{n}\right) + O\left( \frac{n}{t}\beta(\Gamma)\right) \pd.
\end{split}
\end{equation*}
This concludes the proof since by Theorem \ref{thm: improvement of AB for finite sets under DC}, we have
\begin{equation}
    \frac{\bbE_\pi T_A}{\bbE_\alpha T_A} = 1 + O\pg \beta(\Gamma)\pd = 1 + O\left( \frac{n}{t}\beta(\Gamma)\right). \qedhere
\end{equation}
\end{proof}

The following corollary gives a particularly useful approximation of $q_A$ in terms of hitting probabilities.
\begin{corollary}\label{corapproxq_A}
Let $d\geq 2$ and $k\geq 1$. Let $\cF = (\Gamma)$ be a collection of finite (connected) vertex-transitive graphs of degree $d$ satisfying \eqref{e:DC}. 
Then as $|\Gamma| \to \infty$, denoting $n=|\Gamma|$, we have uniformly over all non-empty subsets $A\subset \Gamma$ of cardinality at most $k$ and all $1\leq t\leq n$, 
\begin{equation}
   q_A =  \frac{\bbE_\pi T_o}{\bbE_\pi T_A} = \frac{\bbPpi{T_A \leq t}}{\bbPpi{T_o \leq t}}\pg 1+  O\left( \frac{t}{n}\right) + O\left( \frac{n}{t}\beta(\Gamma)\right) \pd,
\end{equation}
and in particular, at time $t_\Gamma^*$, we have
\begin{equation}\label{eq:approx q A with t star good error term}
    q_A = \frac{\bbPpi{T_A \leq t_\Gamma^*}}{\bbPpi{T_o \leq t_\Gamma^*}}\pg 1 + o\pg \frac{1}{\log n} \pd \pd.
\end{equation}
\end{corollary}
\begin{proof}
    The first point follows from applying Proposition \ref{P:linear} to approximate both the numerator and the denominator of $\frac{\bbE_\pi T_o}{\bbE_\pi T_A}$. The second point follows from the first one and Lemma \ref{lem: bound on beta Gamma biased by time}.
\end{proof}

\begin{remark}\label{R:expapprox and limit of AB}
The ratio $q_A : = \frac{\bbE_\pi T_o}{\bbE_\pi T_A}$ will play a crucial role in our analysis. The approximation in Corollary \ref{corapproxq_A} will be very useful, as the hitting probabilities $\bbPpi{T_A \leq t}$ are easier to estimate than the expectations $\bbE_\pi T_A$. In what follows we will want to apply \eqref{eq:approx q A with t star good error term},
and we will soon see that an error term of $o(1/\log n)$ is indeed sufficient for our purpose.
\end{remark}

We conclude this section showing that $\tsexp$ is close to $\tsav$. 
\begin{prop}\label{P:tdoubleprime s}
    Let $(\Gamma)$ be a collection of finite (connected) vertex-transitive graphs of fixed degree that satisfies  \eqref{e:DC}. Recall that for $s\in \bbR$,  $\tsav = \thitav((\log n) + s)$, and $\tsexp$ is such that $\bbE|U(\tsexp)| = e^{-s}$. Let $s\in \bbR$. As $n=|\Gamma| \to \infty$, we have
    \begin{equation}
       \tsexp = \tsav + o(D^2).
    \end{equation}
\end{prop}
\begin{proof}
First observe that we have for every $t\geq 0$, setting $\alpha = \alpha_o$ the quasi-stationary distribution associated with $A = \{o\}$,
\begin{equation}
    \bbE|U(t)| = n\bbP_\pi(T_o >t) = n\frac{\bbP_\pi(T_o >t)}{\bbP_\alpha(T_o >t)} \bbP_\alpha(T_o >t) = n\frac{\bbP_\pi(T_o >t)}{\bbP_\alpha(T_o >t)}e^{-t/\bbE_\alpha T_o}.
\end{equation}
It follows, since $t\mapsto \bbE|U(t)|$ is decreasing, that $\tsexp$ is unique and that $\tsexp$ satisfies
\begin{equation}
    \tsexp = \bbE_\alpha T_o \pg (\log n) + s - \log \pg \frac{\bbP_\alpha(T_o >\tsexp)}{\bbP_\pi(T_o >\tsexp)}\pd \pd.
\end{equation}
Since \eqref{e:DC} holds, by Theorem \ref{thm: improvement of AB for finite sets under DC} we have $\frac{\bbE_\alpha T_o}{ \bbE_\pi T_o } = 1 + O(\beta(\Gamma))$ and $\log \pg \tfrac{\bbP_\alpha(T_o >\tsexp)}{\bbP_\pi(T_o >\tsexp)}\pd = O(\beta(\Gamma))$. 
This implies that, since $\thitav \leq \thit$, and  $\thit = O (n)$ by Proposition \ref{prop: thit at most order n under DC},
\begin{equation*}
\tsexp = \tsav + O \pg (\thitav \log n)\beta(\Gamma) \pd = \tsav + O \pg (n \log n)\beta(\Gamma) \pd. 
\end{equation*}
Finally, a simple case by case analysis shows that under \eqref{e:DC} we have $(n \log n)\beta(\Gamma) = o(D^2)$. (For the case $q\geq 5$, this uses that $D\gtrsim \log n$ which holds since the degree is bounded.)
This concludes the proof.
\end{proof}

\section{Convergence of the uncovered set}
\label{S:main}

We now start the proof of the main result of this paper. In this section we show that \eqref{e:DC} is a sufficient condition for Gumbel fluctuations, i.e.\ that \eqref{e:DC} implies \eqref{E:main} in Theorem \ref{T:main}. We recall this statement as the following proposition.

\begin{prop}\label{prop: DC implies Gumbel fluctuations}
    Let $(\Gamma)$ be a collection of finite vertex-transitive graphs of fixed degree satisfying \eqref{e:DC}. Let $(\varepsilon_\Gamma)_{\Gamma \in \cF}\in \bbR^\cF$ such that $\varepsilon_\Gamma = o(1/\log |\Gamma|)$ as $|\Gamma|\to \infty$. As $|\Gamma|\to \infty$, $\frac{\taucov}{\thitav(1+\varepsilon_\Gamma)} - \log |\Gamma|$ converges in distribution to a standard Gumbel variable. In particular, as $|\Gamma|\to \infty$, $\frac{\taucov}{\thit} - \log |\Gamma|$ converges in distribution to a standard Gumbel variable.
\end{prop}

\subsection{Strategy}
\label{SS:strategy}

As discussed in Section \ref{s: discussion of proof ideas}, it will be important to estimate the capacities $q_A = \frac{\bbE_\pi T_o}{\bbE_\pi T_A}$ for finite sets $A$.
Initially it is not clear whether one needs to worry about the whole spectrum of possibilities for the mutual distances between points of $A$ or if a cruder bound on the \emph{minimum distance} between points of $A$ is sufficient to distinguish between the good and the bad cases. As it turns out, this cruder strategy is sufficient. We therefore introduce the following quantities.
\begin{defn}
    Let $\Gamma$ be a finite connected vertex-transitive graph. Recall that the graph distance of $\Gamma$ is denoted by $d(\cdot, \cdot)$. For $A\subset \Gamma$ of cardinality at least 2, we set
    \begin{equation}
    \mindist(A) := \min_{x,y\in A, \; x\ne y} d(x,y).
\end{equation}
For $k\geq 2$, let $\cA = \cA(\Gamma, k) $ be the set of all (ordered) $k$-tuples of pairwise distinct elements of $\Gamma$. Here and below, we make another small abuse of notation, where we identify the ordered sets $A\in\cA$ with unordered subsets of $\Gamma$.
For $k\geq 2$ and $1\leq \delta \leq D$, we define
\begin{equation}
    \cA_{\delta} = \cA_\delta(\Gamma, k) :=  \ag A \in \cA(\Gamma,k) \du \mindist(A) \leq \delta \ad.
\end{equation}
Finally, for $k\geq 2$, $1\leq \delta \leq D$, and $s\in \bbR$, we set
\begin{equation}
S_\Gamma(\delta, k, s) :=  \sum_{A \in \cA_{\delta}(\Gamma, k)} n^{-q_A}e^{-q_As}.
\end{equation}
\end{defn}

In order to motivate our proofs, we start by reducing Proposition \ref{prop: DC implies Gumbel fluctuations} to a simple statement about the existence of an intermediate distance satisfying three conditions. In Lemma \ref{lem: good sequence delta implies Gumbel fluctuations} below, Condition (c) ensures that all sets $A$ of a given size $k$ such that $\mindist(A)>\delta_\Gamma$ have almost the same hitting distribution, and hence almost the same contribution to the $k$-th moment; Condition (a) ensures that most sets satisfy $\mindist(A)>\delta_\Gamma$; and Condition (b) ensures that the contribution of other sets to the moments is negligible.

\begin{lemma}\label{lem: good sequence delta implies Gumbel fluctuations}
    Let $\cF = (\Gamma)$ be a collection of finite vertex-transitive graphs of fixed degree satisfying \eqref{e:DC}. Let $(\varepsilon_\Gamma)_{\Gamma \in \cF}\in \bbR^\cF$ such that $\varepsilon_\Gamma = o(1/\log |\Gamma|)$ as $|\Gamma|\to \infty$. Assume that there exists $(\delta_\Gamma)_{\Gamma \in \cF}$ such that as $n= |\Gamma|\to \infty$ the following holds:
    \begin{enumerate}
        \item $\delta_\Gamma = o(\Diam(\Gamma))$;
        \item for each $k\geq 2$, $S_\Gamma(\delta_\Gamma, k, 0) = o(1)$;
        \item for each $k\geq 2$,
        \begin{equation}
            \min_{A \in \cA(\Gamma, k) \backslash \cA_{\delta_\Gamma}(\Gamma, k)} q_A \geq  k-o\pg \frac{1}{\log n} \pd.
        \end{equation}
    \end{enumerate}
Then as $|\Gamma|\to \infty$, $\frac{\taucov}{\thitav(1+\varepsilon_\Gamma)} - \log |\Gamma|$ converges in distribution to a standard Gumbel variable. In particular, as $|\Gamma|\to \infty$, $\frac{\taucov}{\thit} - \log |\Gamma|$ converges in distribution to a standard Gumbel variable.\end{lemma}
\begin{proof}
Let $k\geq 1$ and $s\in \bbR$. Recall that $\tsav = \bbE_\pi T_o (\log(n) + s)$. Denote for this proof $\Tilde{t}_s := \tsav(1+\varepsilon_\Gamma)$ and let $\Tilde{Z}_s = \sum_{x \in \Gamma}\mathbbm{1}_{T_x > \Tilde{t}_s}$ be the size of the uncovered set at time $\Tilde{t}_s$. Identifying the ordered sets $A\in\cA(\Gamma, k)$ with subsets of $\Gamma$, we can rewrite the factorial moments of $\Tilde{Z}_s$ as
\begin{equation}
    \esppi{(\Tilde{Z}_s)^{\downarrow k}} = \sum_{A \in \cA(\Gamma, k)} \bbPpi{T_A > \Tilde{t}_s},
\end{equation}
where $z^{\downarrow k} = z(z-1) \cdots (z-k+1)$. We deduce from Proposition \ref{propapproxpita} that
\begin{equation}\label{approxepizk}
    \esppi{(\Tilde{Z}_s)^{\downarrow k}}
    =  (1+o(1))\sum_{A \in \cA(\Gamma, k)} n^{-q_A} e^{-q_As} = (1+o(1))S_\Gamma(D, k, s).
\end{equation}
Since $\cA(\Gamma, 1)$ is identified with $\Gamma$ and $q_A = 1$ if $A$ is a singleton, we deduce that as $n\to \infty$,
\begin{equation}
    \bbE_\pi[\Tilde{Z}_s] \to e^{-s}.
\end{equation}
From now on we may therefore assume that $k\geq 2$. 
Let $(\delta_\Gamma)$ be a sequence satisfying the three assumptions of the lemma. For any subset $A\ne \eset$ we have $q_A \geq q_{\ag o\ad} = 1$. Moreover by Corollary \ref{corapproxq_A} and a union bound $\P_\pi(T_A < t_\Gamma^*) \le k \P_\pi (T_o < t_\Gamma^*)$ on the numerator, we also always trivially have as $n\to \infty$
\begin{equation}\label{qAub}
  \max_{A\in \cA(\Gamma, k)} q_A \le k + o\pg \frac{1}{\log n}\pd = k + o(1).
\end{equation}
 Such an \emph{upper bound} on $q_A$ is useful since $q_A$ appears not only in the term $n^{-q_A}$ but also in the term $e^{-q_A s}$, while $s \in \bbR$ can be negative. Therefore $S_\Gamma(\delta_\Gamma, k, s) \lesssim_{k,s} S_\Gamma(\delta_\Gamma, k, 0)$, and using Assumption (b) we obtain as $n\to\infty$
 \begin{equation}\label{eq: good sequence delta implies Gumbel fluctuations inter 1}
S_\Gamma(\delta_\Gamma, k, s) = o(1).
 \end{equation}
Denote $\cB := \cA(\Gamma, k) \backslash \cA_{\delta_\Gamma}(\Gamma, k)$. By Assumption (a) and Lemma \ref{lem:volume bounds} (b), we have $V(\delta_\Gamma) = o(n)$. We deduce that $\left|\cA_{\delta_\Gamma}(\Gamma, k) \right| = o(n^k)$, and since $\left|\cA(\Gamma, k) \right| = n^{\downarrow k} =n^k(1-o(1))$, we obtain
\begin{equation}\label{eq: good sequence delta implies Gumbel fluctuations inter 2}
   |\cB| = n^k(1-o(1)).
\end{equation}  
Moreover, by Assumption (c) and \eqref{qAub}, as $n\to \infty$ we have $q_A = k + o(1/\log n)$ uniformly over all $A\in \cB$, and therefore, using again Proposition \ref{propapproxpita} for the first equality,
\begin{equation}
    S_\Gamma(D, k,s) - S_\Gamma(\delta_\Gamma, k, s) = (1+o(1)) \sum_{A\in \cB} n^{-q_A} e^{-q_A s} = (1+o(1)) \frac{|\cB|}{n^k} e^{-ks}.
\end{equation}
It then follows from \eqref{eq: good sequence delta implies Gumbel fluctuations inter 1} and \eqref{eq: good sequence delta implies Gumbel fluctuations inter 2} that as $n\to \infty$, $ S_\Gamma(D, k,s) \to e^{-ks}$, and therefore by \eqref{approxepizk} that 
\begin{equation}
   \esppi{(\Tilde{Z}_s)^{\downarrow k}} \to e^{-ks}.
\end{equation}
We have proved that the factorial moments (for each $k\geq 1$) of $\Tilde{Z}_s$ converge to those of a Poisson variable with parameter $e^{-s}$. A standard application of the method of moments then shows that $
\Tilde{Z}_s$ converges in distribution to a Poisson random variable with parameter $e^{-s}$: that is, for all $k\geq 0$ we have
\begin{equation}
  \label{eq:poisson_conv}
  \P( \Tilde{Z}_s = k) \to \exp( - e^{-s}) \frac{e^{-ks}}{k!}.
\end{equation}
Therefore as $|\Gamma|\to \infty$ we have
\begin{equation}
\P( \taucov \leq \Tilde{t}_s ) = \P( \Tilde{Z}_s = 0) \to \exp( - e^{-s}).
\end{equation}
This shows that $\frac{\taucov}{\thitav(1+\varepsilon_\Gamma)}- \log |\Gamma|$ converges in distribution to a standard Gumbel variable, and the same holds for $\frac{\taucov}{\thit} - \log |\Gamma|$ since $\thit = \thitav(1+o(1/\log n))$ by Proposition \ref{P:thitmaxav}.
\end{proof}

The goal of this section is therefore to prove that under \eqref{e:DC}, there exists a sequence $(\delta_\Gamma)$ as in Lemma \ref{lem: good sequence delta implies Gumbel fluctuations}.
Given a finite connected vertex-transitive graph $\Gamma$  such that $n\geq D^2 \log n$, where $n=|\Gamma|$ and $D = \Diam(\Gamma)$, we recall that  
\begin{equation}
    b(\Gamma) = \frac{n}{D^2\log n} = \frac{f(\Gamma)}{\log n},
\end{equation}
and we set
\begin{equation}\label{eq: definition of delta star}
    \delta^*_\Gamma := \begin{cases}
       \sqrt{\log n} & \text{ if } n> D^5, \\
       \sqrt{D} & \text{ if } D^2(\log n)^{9/4} < n \leq D^5,\\
       D/e^{\sqrt{b(\Gamma)}} &  \text{ if } D^2\log n < n \leq D^2(\log n)^{9/4}.
    \end{cases}
\end{equation}
We will prove that $(\delta^*_\Gamma)$ satisfies the assumptions of Lemma \ref{lem: good sequence delta implies Gumbel fluctuations}. Since $(\delta^*_\Gamma)$ trivially satisfies Assumption (a), we have to prove that it also satisfies Assumptions (b) and (c).

\subsection{Skeleton of a set}
\label{SS:skeleton}

We will need to group the points of $A \in \cA_\delta$ into subsets of points which are close to one another. Let us define an equivalence relation which will allow us to realise such partitions in a convenient way.
\begin{defn}
Let $A$ be a subset of cardinality at least 2 of a finite connected vertex-transitive graph $\Gamma$, and let $1 \leq \delta \leq \Diam(\Gamma)$. We say that two points $x,y\in A$ are $(A,\delta)$-linked if there exist an integer $r\geq 2$ and a sequence $x=x_1,...,x_r = y$ of points in $A$ such that for all $1\leq i \leq r-1$, $d(x_i,x_{i+1}) \leq \delta$. We will write $A$ under the form
\begin{equation}
    A = A_1 \sqcup ... \sqcup A_\ell,
\end{equation}
where the $A_i$ are the $(A,\delta)$-connected components of $A$, and $\ell_{\delta}(A)$ is the number of such components.
\end{defn}

\begin{remark}
It is worth noting that some points can be closer to points in different connected components than to some points in their own components. For example on $\bbZ$, take $\delta = 2$ and $A = \ag 0,2,4,7 \ad$. The partition is then $\ag\ag 0,2,4\ad,\ag 7\ad\ad$, though 4 is closer to 7 than to 0.
\end{remark}

Let $1\leq \delta \leq D$ and $A = A_1 \sqcup \ldots \sqcup A_\ell \in\cA_\delta$. Since at least two points of $A$ are in the same component, we assume without loss of generality that $\abs{A_{1}}\geq 2$, and we have $\ell_\delta(A) \leq k - 1$. Recall that by Corollary \ref{corapproxq_A} we have (under \ref{e:DC})
\begin{equation}\label{E:o(1) approximation of q_A}
    q_A = \frac{\bbPpi{T_A \leq t_\Gamma^*}}{\bbPpi{T_o \leq t_\Gamma^*}} \pg 1 + o\pg \frac{1}{\log n} \pd\pd = \frac{\bbPpi{\bigcup_{i=1}^\ell \ag T_{A_i} \leq t_\Gamma^*\ad}}{\bbPpi{T_o \leq t_\Gamma^*}} \pg 1 + o\pg \frac{1}{\log n} \pd\pd,
\end{equation}
where $t_\Gamma^*$ was defined in \eqref{eq: definition of t Gamma star}.
We want to show that when the sets $A_i$ are quite far from one another, \report{the events $\ag T_{A_i} \leq t_\Gamma^*\ad$ are essentially independent.}

\begin{prop}\label{propeiej}
Let $d\geq 2$. Let $\cF = (\Gamma)$ be a collection of finite (connected) vertex-transitive graphs of degree $d$ satisfying \eqref{e:DC}. Let $k\geq 2$. There exists a constant $C = C(k)>0$ such that the following holds. Let $\Gamma \in \cF$, $\delta \in [1, \Diam(\Gamma)]$, $A = A_1 \sqcup ... \sqcup A_\ell \in \cA_\delta(\Gamma, k)$, where each $A_i$ is $(A, \delta)$-connected and $|A_1|\geq 2$. For any $1\leq i\ne j \leq \ell$, we have
\begin{equation}
\frac{\bbP_\pi(\max(T_{A_i}, T_{A_j})<t_\Gamma^*)}{ \bbPpi{T_o < t_\Gamma^*}} \leq \report{C \max_{\underset{d(x,y) \geq \delta}{x,y\in \Gamma}} \bbP_x(T_y < t_\Gamma^*).}
\end{equation}
\end{prop}
\begin{proof}
Denote $t^* = t_\Gamma^*$ and let $x,y \in \Gamma$. By the Markov property we have
\begin{equation}
\begin{split}
     \bbP_\pi(T_x < T_y <t^*)  = \bbP_\pi(T_x <t^*)\bbP_\pi(T_x < T_y < t^* \mid T_x < t^*) & \le \bbP_\pi(T_x<t^*) \bbP_x( T_y < t^*)\\
    & = \bbP_\pi(T_o<t^*) \bbP_x( T_y < t^*).
\end{split}
\end{equation}
  Moreover $\bbP_\pi(T_x < T_y <t^*)=\bbP_\pi(T_y < T_x <t^*)$ and $\bbP_x(T_y < t^*) = \bbP_y(T_x< t^*)$ by symmetry, so altogether,
\begin{equation}
    \bbP_\pi(\max(T_x,T_y)<t^*) = 2\bbP_\pi(T_x < T_y <t^*) \leq 2\bbP_\pi(T_o<t^*)\bbP_x( T_y < t^*).
\end{equation}
Therefore for any $1\leq i\ne j \leq \ell$, we have
\begin{equation}
\begin{split}
    \bbP_\pi(\max(T_{A_i}, T_{A_j})<t^*) & = \bbP_\pi \pg \bigcup_{x\in A_i, y \in A_j} \ag \max(T_x,T_y)<t^* \ad \pd \\ & \leq \abs{
A_i}\abs{A_j} \max_{\underset{d(x,y) \ge \delta}{x,y\in \Gamma}} \bbP_\pi(\max(T_x,T_y)<t^*) \\
& \leq 2k^2 \bbP_\pi(T_o<t^*) \max_{\underset{d(x,y) \ge \delta}{x,y\in \Gamma}} \bbP_x(T_y < t^*),
\end{split}    
\end{equation}
which concludes the proof.
\end{proof}
\begin{remark}
    With the notation of Proposition \ref{propeiej}, by the strong Markov property we have
\begin{equation*}
    \max_{\underset{d(x,y) \ge  \delta}{x,y\in \Gamma}} \bbP_x(T_y < t^*)=\max_{\underset{d(x,y) = \lceil \delta \rceil}{x,y\in \Gamma}} \bbP_x(T_y < t^*). \qedhere
\end{equation*}  
\end{remark}

Let us end this subsection with a useful lower bound on $q_A$.
\begin{prop}\label{P: lower bound on q_A}
Let $d\geq 2$. Let $\cF = (\Gamma)$ be a collection of finite (connected) vertex-transitive graphs of degree $d$ satisfying \eqref{e:DC}. Let $k\geq 2$. As $n=|\Gamma| \to \infty$, we have uniformly over all $\delta \in [1, \Diam(\Gamma)]$ and all $A = A_1 \sqcup ... \sqcup A_\ell \in \cA_\delta(\Gamma, k)$, where each $A_i$ is $(A, \delta)$-connected and $|A_1|\geq 2$,
\begin{equation}
q_A - \summ{i=1}{\ell_{\delta}(A)} \frac{\bbP_\pi(T_{A_i}<t_\Gamma^*)}{\bbP_\pi(T_o < t_\Gamma^*)} \gtrsim - \pg \max_{\underset{d(x,y) \geq \delta}{x,y\in \Gamma}} \bbP_x(T_y < t_\Gamma^*) + o\pg \frac{1}{\log n} \pd\pd;
\end{equation}
and in particular, if $\delta < \mindist(A)$, we have
\begin{equation}\label{qAlb_spec}
q_A - k \gtrsim -\max_{\underset{d(x,y) \geq \delta}{x,y\in \Gamma}} \bbP_x(T_y < t_\Gamma^*) + o\pg \frac{1}{\log n} \pd.
\end{equation}
\end{prop}
\begin{proof}
Write again $t^* = t_\Gamma^*$.
From the Bonferroni inequality, we have
\begin{equation}
\begin{split}
     \bbPpi{T_A < t^*} & = \bbP_\pi\pg \bigcup_{i=1}^{\ell_{\delta}(A)} \ag T_{A_i}<t^* \ad \pd \\
     & \geq \summ{i=1}{\ell_{\delta}(A)} \bbP_\pi(T_{A_i}<t^* ) - \sum_{1\leq i < j \leq \ell_{\delta}(A)} \bbP_\pi(T_{A_i}<t^* , T_{A_j}<t^* ).
\end{split}
\end{equation}
Dividing both sides by $\bbPpi{T_o < t^*}$, we obtain
\begin{equation}
\frac{\bbP_\pi(T_A < t^*)}{\bbP_\pi(T_o < t^*)} - \summ{i=1}{\ell_{\delta}(A)} \frac{\bbP_\pi(T_{A_i}<t^*)}{\bbP_\pi(T_o < t^*)} \gtrsim - \sum_{1\leq i < j \leq \ell_{\delta}(A)} \frac{\bbP(T_{A_i}<t^* , T_{A_j}<t^*)}{\bbP_\pi(T_o < t^*)}.
\end{equation}
The result then follows from Corollary \ref{corapproxq_A} and Proposition \ref{propeiej}.
\end{proof}

\subsection{Bounds on capacities for sets with a large minimal distance}
The aim of this subsection is to prove that the sequence $(\delta^*_\Gamma)$ satisfies Assumption (c) from Lemma \ref{lem: good sequence delta implies Gumbel fluctuations}. This is contained in the following proposition.
\begin{prop}\label{prop: lower bound capacities high mindist assumption c}
     Let $\cF = (\Gamma)$ be a collection of finite vertex-transitive graphs of fixed degree satisfying \eqref{e:DC}. Let $(\delta^*_\Gamma)$ as defined in \eqref{eq: definition of delta star}. Let $k\geq 2$. As $n=|\Gamma| \to \infty$, uniformly over all $A \in \cA(\Gamma, k) \backslash \cA_{\delta_\Gamma^*}(\Gamma, k)$ we have
        \begin{equation*}
             q_A = k + o\pg \frac{1}{\log n} \pd.
        \end{equation*}
\end{prop}
\begin{proof}
    By \eqref{qAub} we only need to prove that $q_A \geq k - o\pg \frac{1}{\log n} \pd$. Denote $t^* = t^*_\Gamma$ and $\delta^* = \delta^*_\Gamma$. In regard of Corollary \ref{corapproxq_A}, Proposition \ref{P: lower bound on q_A}, and Lemma \ref{L:unionbound}, it is enough to show that 
    \begin{equation}
        \max_{\underset{d(x,y) \geq  \delta^*}{x,y\in \Gamma}} \bbE_x[L_y(t^*)] = o\pg \frac{1}{\log n} \pd.
    \end{equation}
Recall from Lemma \ref{lem:volume bounds} and Remark \ref{rem: t star lesssim D power 5/2} that $D\gtrsim \log n$ and $t^* \lesssim D^{5/2}$. We split the proof depending on the diameter, as in the definition of $\delta^*$.

First assume that $n> D^5$. Then $\delta^* = \sqrt{\log n}$, $(\delta^*)^{-3} = (\log n)^{-3/2} = o(1/\log n)$, and $D^{5/2}/D^5 = D^{-5/2} \lesssim (\log n)^{-5/2} = o(1/\log n)$. Therefore by Proposition \ref{P:bound on transition probabilities before some time} (b) we have
\begin{equation}
    \max_{\underset{d(x,y) \geq  \delta^*}{x,y\in \Gamma}}\bbE_x L_y(t^*) \lesssim \frac{1}{(\delta^*)^3} + \frac{t^*}{D^5} \lesssim \frac{1}{(\delta^*)^3} + \frac{D^{5/2}}{D^{5}} = o\pg \frac{1}{\log n} \pd.
\end{equation}

Now assume that $D^2(\log n)^{9/4} < n \leq D^5$. Then $\delta^* = \sqrt{D}$ by definition so $1/\delta^* = o(1/\log n)$ since $D\geq n^{1/5}$. Moreover, $$\log(D/\delta^*) = \log D^{1/2} \asymp \log n$$ so $$\log(D/\delta^*)/f(\Gamma) \lesssim (\log n)^{-{5/4}} = o(1/\log n).$$ 
We also have $$t^*/D^3 \lesssim D^{5/2}/D^{3} \leq n^{-1/10} =o(1/\log n),$$ and similarly, 
$$t^*/n \lesssim D^2 \sqrt{b(\Gamma)}/(D^2 b(\Gamma) \log n) = b(\Gamma)^{-1/2}/\log n = o(1/\log n).$$ Therefore by Proposition \ref{P:bound on transition probabilities before some time} (b),
\begin{equation}
   \max_{\underset{d(x,y) \geq  \delta^*}{x,y\in \Gamma}} \bbE_x L_y(t^*) \lesssim  \frac{1}{{\delta^*}} + \frac{\log (D/\delta^*)}{f(\Gamma)} + \frac{t^*}{n} + \frac{t^*}{D^3} =o\pg \frac{1}{\log n} \pd.
\end{equation}
Assume finally that $D^2\log n \ll n \leq D^2(\log n)^{9/4}$, so $1\ll b (\Gamma) \le (\log n)^{5/4}$. Then $\delta^* = D/e^{\sqrt{b(\Gamma)}} \ge \sqrt{D} \to \infty$. The terms $1/\delta^*$, $t^*/n$ and $t^*/D^3$ are bounded in the same way as in the previous case. The difficulty is to bound $\log(D/\delta^*)/f(\Gamma)$. But by our choice of $\delta^*$, we have 
\begin{equation}
    \frac{\log (D/\delta^*)}{f(\Gamma)} = \frac{ \sqrt{b(\Gamma)}}{b(\Gamma)\log n} =\frac{1}{\sqrt{b(\Gamma)}\log n} = o\pg \frac{1}{\log n}\pd.
\end{equation}
Applying Proposition \ref{P:bound on transition probabilities before some time} (b) again concludes the proof.
\end{proof}

We have seen that $(\delta^*_\Gamma)$ satisfies Assumptions (a) and (c) of Lemma \ref{lem: good sequence delta implies Gumbel fluctuations}. The rest of this section is dedicated to show that it also satisfies Assumption (b).

\subsection{Bounds on capacities and volumes for sets with a small minimal distance}

The aim of this subsection is to prove lower bounds on the capacities $q_A$ that will be sufficient to show that $(\delta^*_\Gamma)$ satisfies Assumption (b) of Lemma \ref{lem: good sequence delta implies Gumbel fluctuations}. We start with sets $A$ which have two very close points.

\begin{prop}\label{propmicroqa}
Let $d\geq 2$. Let $\cF = (\Gamma)$ be a collection of finite (connected) vertex-transitive graphs of degree $d$ satisfying \eqref{e:DC}. For $\Gamma \in \cF$ set $\delta_\text{micro} = \delta_\text{micro}(\Gamma):= \sqrt{\log |\Gamma|}$. Let $k\geq 2$. There exists a constant $\beta' = \beta'(k,d)>0$ such that as $n=|\Gamma| \to \infty$, uniformly over all $A \in \cA_{\delta_\text{micro}(\Gamma)}(\Gamma, k)$, 
\begin{equation}
    q_A \geq \ell_{\delta_\text{micro}}(A) + \beta' +o(1).
\end{equation}
\end{prop}
\begin{proof}
Let $\Gamma \in \cF$ and $A\in \cA_{\delta_{\text{micro}}}(\Gamma, k)$. To lighten the notation, we write $\delta$ for $\delta_{\text{micro}}$ and $t^*$ for $t^*_\Gamma$ in this proof. Write $A = A_1 \sqcup ... \sqcup A_\ell \in \cA_\delta(\Gamma, k)$, where each $A_i$ is $(A, \delta)$-connected and $|A_1|\geq 2$.
Applying the lower bound on $q_A$ proven in Proposition \ref{P: lower bound on q_A}, together with Corollary \ref{proppxygrand}, we have,
\begin{equation}
q_A - \summ{i=1}{\ell_{\delta}(A)} \frac{\bbP_\pi(T_{A_i}<t^*)}{\bbP_\pi(T_o < t^*)} \geq o(1).
\end{equation}
Let $x,y$ be distinct points of $A_1$ such that $d(x,y) \leq \delta$.
By symmetry, we have
\begin{equation}\label{evite1}
         \bbP_\pi(T_{\ag x,y\ad} < t^*) + \bbP_\pi (T_x<t^*, T_y<t^*) = \bbP_\pi (T_x<t^*) + \bbP_\pi (T_y<t^*) = 2 \bbP_\pi (T_o<t^*).
\end{equation}
Now, by the strong Markov property and symmetry, and since $t^* = t^*_\Gamma = o(|\Gamma|)$ by definition,
\begin{equation}\label{evite2}
     \P_\pi( T_x < t^*, T_y  <t^*)  \le \P_\pi ( T_{\{x,y\}} < t^*)\bbP_x(T_y< t^*) \leq \P_\pi ( T_{\{x,y\}}<t^*) (1- \beta +o(1)),
\end{equation}
where $\beta$ is the constant from Proposition \ref{proppxysmall}. Combining \eqref{evite1} and \eqref{evite2} together we get
\begin{equation*}
\bbP_\pi(T_{\ag x,y \ad}<t^*) (2- \beta +o(1)) \geq 2 \bbP_\pi(T_o<t^*),
\end{equation*}
so that, setting $\beta':= \frac{2}{2-\beta}-1$ and using that $\bbP_\pi(T_{A_1}<t^*) \geq \bbP_\pi(T_{\ag x,y \ad}<t^*)$,
\begin{equation*}
\frac{\bbP_\pi(T_{A_1}<t^*)}{\bbP_\pi(T_{o}<t^*)}\geq 1+ \beta'+ o(1).
\end{equation*}
We can bound the probability of hitting the $\ell_{\delta}(A) - 1$ other $A_i$'s before time $t^*$ very crudely: as the sets $A_i$ are non-empty, we have
\begin{equation}
\frac{\bbP_\pi(T_{A_i}<t^*)}{\bbP_\pi(T_o < t^*)} \geq 1.
\end{equation}
Finally, putting everything together, we obtain
\begin{equation}
    q_A \geq 1+ \beta'  + (\ell_{\delta}(A) - 1) + o(1) = \ell_{\delta}(A) + \beta' +o(1). \qedhere
\end{equation}
\end{proof}

The bound on $q_A$ we obtained depends on the number of components in the $\delta$-skeleton of $A$. Given a finite connected vertex-transitive graph $\Gamma$, an integer $k\geq 2$, and $\delta \in [1, \Diam(\Gamma)]$, we set for $1\leq \ell \leq k$,
\begin{equation}
    N_\Gamma(\delta, k,\ell) := \left|\ag A\in \cA_{\delta}(\Gamma, k) \du \ell_{\delta}(A) = \ell \ad\right|.
\end{equation}
\begin{lemma}\label{lem: bound N delta Gamma k ell}
    Let $\Gamma$ be a finite connected vertex-transitive graph.
    Let $k\geq 2$,  $1\leq \ell \leq k-1$, and $\delta \in [1, \Diam(\Gamma)]$. Then
    \begin{equation}
        N_{\Gamma}(\delta, k,\ell) \leq k^k n^{\ell}V(\delta)^{k-\ell} \leq k^k n^{k-1}V(\delta). \qedhere
    \end{equation}
\end{lemma}
\begin{proof}
    Sets $A$ of size $k$ with $\ell$ $\delta$-components can all be obtained as follows. First we pick $\ell$ vertices in $\Gamma$: there are at most $n^\ell$ possibilities, where $n=|\Gamma|$. Then add iteratively $k-\ell$ vertices, which are at distance at most $\delta$ from one of the vertices that we already have, i.e.\ are in the ball of radius $\delta$ around one of at most $k$ vertices: at each of these $k-\ell$ steps there are at most $kV(\delta)$ possibilities. It follows, using that $V(\delta)\leq n$ for the last inequality, that
\begin{equation}
    N_\Gamma(\delta, k,\ell) \leq n^\ell (kV(\delta))^{k-\ell} \leq k^k n^\ell V(\delta)^{k-\ell} \leq k^kn^{k-1}V(\delta). \qedhere
\end{equation}
\end{proof}

We can already show that Assumption (b) of Lemma \ref{lem: good sequence delta implies Gumbel fluctuations} is satisfied for graphs with a small diameter (i.e.\ if $n > D^5$).

\begin{prop}\label{prop: bound S micro case}
     Let $\cF = (\Gamma)$ be a collection of finite vertex-transitive graphs of fixed degree satisfying \eqref{e:DC}, and let $k\geq 2$. As $n=|\Gamma| \to \infty$, we have $S_\Gamma(\sqrt{\log |\Gamma|}, k, 0) = o(1)$.
     In particular, if $|\Gamma| > \Diam(\Gamma)^5$ (for $|\Gamma|$ large enough)
     $S_\Gamma(\delta_\Gamma^*, k, 0) = o(1)$.
\end{prop}
\begin{proof}
   Since $\sqrt{\log|\Gamma|} = o(\Diam(\Gamma))$, for each $1\leq \ell \leq k-1$ by Lemma \ref{lem:volume bounds} (a) we have that $V(\sqrt{\log|\Gamma|}) = n^{o(1)}$, and therefore by Lemma \ref{lem: bound N delta Gamma k ell} that $N_\Gamma(\sqrt{\log|\Gamma|},k,\ell) = n^{\ell + o(1)}$. Combining this with Proposition \ref{propmicroqa}, we obtain
\begin{equation}
    S_\Gamma(\delta, k, 0) \leq \summ{\ell=1}{k-1} N_\Gamma(\sqrt{\log|\Gamma|},k,\ell) n^{-(\ell + \beta' +o(1))} \lesssim n^{-\beta' +o(1)} = o(1). \qedhere
\end{equation}
\end{proof}

Let us now show that Assumption (b) of Lemma \ref{lem: good sequence delta implies Gumbel fluctuations}  is also satisfied for graphs with an intermediate diameter (i.e.\ if $D^2(\log n)^{9/4} < n \leq D^5$).

\begin{lemma}\label{lem: bound q A meso case}
Let $\cF = (\Gamma)$ be a collection of finite vertex-transitive graphs of fixed degree satisfying \eqref{e:DC}. Let $(\delta_\Gamma)_{\Gamma \in \cF}$ such that $\delta_\Gamma \to \infty$ as $|\Gamma|\to \infty$. Let $k\geq 2$. As $n=|\Gamma|\to \infty$, uniformly over all $A\in \cA(\Gamma, k) \backslash \cA_{\delta_\Gamma}(\Gamma, k)$, we have
 \begin{equation}
            q_A =  k + o(1).
        \end{equation}
In particular, this holds for $\delta_\Gamma = \delta_{\text{micro}} = \sqrt{\log |\Gamma|}$.
\end{lemma}
\begin{proof}
    Since $\delta_{\Gamma}\to \infty$ as $|\Gamma|\to \infty$, the lower bound follows immediately from Proposition \ref{P: lower bound on q_A} (see in particular \eqref{qAlb_spec}) and Corollary \ref{proppxygrand}, and the upper bound was already written as \eqref{qAub}.
\end{proof}

\begin{prop}\label{prop: bound S meso case}
     Let $(\Gamma)$ be a collection of finite vertex-transitive graphs of fixed degree satisfying \eqref{e:DC} and $n \leq D^5$ for $n=|\Gamma|$ large enough, where $D=\Diam(\Gamma)$.
     Let $k\geq 2$. As $|\Gamma| \to \infty$, we have $S_\Gamma(\sqrt{D}, k, 0) = o(1)$.
     In particular, if $D^2(\log n)^{9/4} < n \leq D^5$ (for $n$ large enough), we have
    $S_\Gamma(\delta_\Gamma^*, k, 0) = o(1)$.
\end{prop}
\begin{proof}
    Write $\delta_{\text{micro}} = \sqrt{\log n}$ and $\delta_{\text{meso}} = \sqrt{D}$. Since $n\leq D^5$ by assumption, for $n$ large enough we have $\delta_{\text{micro}}< \delta_{\text{meso}}$. 
By Lemma \ref{lem:volume bounds} (b), observing that $n\leq D^5$ is equivalent to $D^{1/2}\geq n^{1/10}$, we have
\begin{equation}
V(\delta_{\text{meso}}) \lesssim \frac{n \delta_{\text{meso}}}{D} \leq n^{9/10}.
\end{equation}
Plugging this into Lemma \ref{lem: bound N delta Gamma k ell}, we obtain that for each $1\leq \ell \leq k-1$,
\begin{equation}
     N_\Gamma(\delta_{\text{meso}}, k,\ell) \lesssim n^{k-1}V(\delta_{\text{meso}}) \lesssim n^{k-1/10}.
\end{equation}
We deduce, since $S_\Gamma(\delta_{\text{micro}} = o(1)$ by Proposition \ref{prop: bound S micro case}, that 
\begin{equation}
    S_\Gamma(\delta_{\text{meso}}, k, 0) \lesssim S_\Gamma(\delta_{\text{micro}}, k, 0) + \sum_{\ell=1}^{k-1} n^{k-1/10} n^{-k+o(1)} \lesssim o(1) + n^{-1/10 + o(1)} = o(1),
\end{equation}
as desired.
\end{proof}

We delay the lower bound on $q_A$ for graphs with a larger diameter to Proposition \ref{P:bootstrapbasis}.

\subsection{Bootstrap argument}
\label{SS:bootstrap}
To show that $(\delta_\Gamma^*)$ satisfies Assumption (b) when $n > D^2(\log n)^{9/4}$, we had to consider only one scale ($\delta \leq \delta_{\text{micro}} = \sqrt{\log n}$) and had to to bootstrap this once if  $D^2(\log n)^{9/4} < n \leq D^5$, to extend the control up to $\delta_{\text{meso}}$. For graphs that merely satisfy the diameter condition, we need to bootstrap more than once, and a diverging number of times if $n=D^2(\log n)^{1+o(1)}$.
The basis of the bootstrap argument is the following estimate. 

\begin{prop}\label{P:bootstrapbasis}
Let $\cF = (\Gamma)$ be a collection of finite vertex-transitive graphs of fixed degree satisfying \eqref{e:DC} and $n \leq D^2(\log n)^{9/4}$, where $n = |\Gamma|$ and $D = \Diam(\Gamma)$. Let $k\geq 2$. There exists a constant $K$ that depends on $k$ and the degree of the graphs such that as $|\Gamma| \to \infty$, uniformly over all $D^{1/2}\leq \delta \leq D/2$ and $A\in \cA(\Gamma, k)$ such  that $\delta \leq \mindist(A)$, we have
\begin{equation}
    q_A \geq k - K \frac{\log(D/\delta)}{f(\Gamma)} + o\pg \frac{1}{\log n}\pd.
\end{equation}
\end{prop}
\begin{proof}
Let $\Gamma \in \cF$. Let $t^* = t^*_\Gamma$ as defined in \eqref{eq: definition of t Gamma star}. Let $D^{1/2}\leq \delta \leq D/2$ and $A\in \cA(\Gamma, k)$ such that $\delta \leq \mindist(A)$.
As $|\Gamma|\to \infty$ we have $1/\delta \leq D^{-1/2} = o(1/\log n)$ and
\begin{equation}
    \frac{t^*}{D^3} \leq \frac{t^*}{n} \leq  \frac{\log b(\Gamma)}{f(\Gamma)} \leq \frac{\sqrt{b(\Gamma)}}{b(\Gamma)\log n} = o\pg \frac{1}{\log n}\pd.
\end{equation} We deduce from Lemma \ref{L:unionbound} and Proposition \ref{P:bound on transition probabilities before some time} (a) that 
\begin{equation}
   \max_{\underset{d(x,y) \geq \delta}{x,y\in \Gamma}} \bbP_x(T_y < t^*) \lesssim \bbE_x L_y(t^*) \lesssim \frac{1}{{\delta}} + \frac{\log (D/\delta)}{f(\Gamma)} + \frac{t^*}{n} + \frac{t^*}{D^3} = \frac{\log (D/\delta)}{f(\Gamma)} + o\pg \frac{1}{\log n}\pd.
\end{equation}
Applying Proposition \ref{P: lower bound on q_A} concludes the proof.
\end{proof}

\begin{lemma}\label{L:iteration}
Let $(\Gamma)$ be a collection of finite vertex-transitive graphs of fixed degree satisfying \eqref{e:DC} and $n \leq D^2(\log n)^{9/4}$, where $n = |\Gamma|$ and $D = \Diam(\Gamma)$. Let $k\geq 2$. As $n=|\Gamma| \to \infty$, uniformly over all $\delta, \delta'$ such that $D^{1/2} \leq \delta \leq \delta' \leq D/2$, we have
\begin{equation}
     S_\Gamma(\delta', k, 0) - S_\Gamma(\delta, k, 0) \lesssim \frac{\delta'}{D} \pg \frac{D}{\delta}\pd^{K/b(\Gamma)},
\end{equation}
where $K$ is the constant from Proposition \ref{P:bootstrapbasis}. 
\end{lemma}
\begin{proof}
Let $\sqrt{D} \leq \delta \leq \delta' \leq D/2$. Write $\cA_{\delta}$ and $ \cA_{\delta'}$ for $\cA_{\delta}(\Gamma, k)$ and $ \cA_{\delta'}(\Gamma, k)$.
From Proposition \ref{P:bootstrapbasis}, we have for $A \in \cA_{\delta'}\backslash \cA_{\delta}$,
\begin{equation}
    n^{-q_A}n^{-o(1/\log n)} \leq n^{-k} n^{K\log(D/\delta)/f(\Gamma)} = n^{-k} (D/\delta)^{K(\log n)/f(\Gamma)} = n^{-k} (D/\delta)^{K/b(\Gamma)},
\end{equation}
so for $n$ large enough we have 
\begin{equation}
    n^{-q_A} \leq 2 n^{-k} (D/\delta)^{K/b(\Gamma)}.
\end{equation}
Moreover, by the volume bound from Lemma \ref{lem:volume bounds} (b), we have
\begin{equation}
   \bg \cA_{\delta'}\backslash \cA_{\delta} \bd \leq  \bg \cA_{\delta'} \bd \leq  \binom{k}{2}V(\delta')n^{k-1} \leq \binom{k}{2}\frac{3\delta'}{D}n n^{k-1} \leq \frac{3k^2}{2}\frac{\delta'}{D}n^{k}
\end{equation}
We conclude that for $n$ sufficiently large we have
\begin{equation}
   S_\Gamma(\delta', k, 0) - S_\Gamma(\delta, k, 0) = \sum_{A \in \cA_{\delta'}\backslash \cA_{\delta}} n^{-q_A} \leq 3k^2 \frac{\delta'}{D} \pg \frac{D}{\delta}\pd^{K/b(\Gamma)}. \qedhere
\end{equation}
\end{proof}

Lemma \ref{L:iteration} allows us to increase the value of $\delta$ iteratively in such a way that $S_\Gamma(\delta, k, 0) = o(1)$ (except for macroscopic scales).

\begin{defn}
    Let $\Gamma$ be a finite connected vertex-transitive graph such that $D^2\log n <n \leq D^2 (\log n)^{9/4}$, where $n=|\Gamma|$ and $D = \Diam(\Gamma)$. We set
    \begin{equation}
J = J(\Gamma) := \left\lfloor 4 \frac{\log\log n}{\log(b(\Gamma))} -1\right\rfloor.
\end{equation}
(Note that $J\ge 2$ since $b(\Gamma) \le (\log n)^{5/4}$.) Furthermore, we set for $1\leq j \leq J$,
\begin{equation}
    \delta_j = \delta_j(\Gamma) := D \exp\left(-\frac{\log n}{b(\Gamma)^{j/4}}\right) = D n^{-1/b(\Gamma)^{j/4}}.
\end{equation}
This defines a sequence of scales $(\delta_j)_{1\le j \le J}$.
\end{defn}
We can now prove that Assumption (b) of Lemma \ref{lem: good sequence delta implies Gumbel fluctuations} is satisfied also for graphs that barely satisfy the diameter condition.
\begin{prop}\label{P:bootstrap}
Let $(\Gamma)$ be a collection of finite vertex-transitive graphs of fixed degree satisfying \eqref{e:DC} and $n \leq D^2(\log n)^{9/4}$ for $n=|\Gamma|$ large enough, where $D=\Diam(\Gamma)$. 
     Let $k\geq 2$. As $|\Gamma| \to \infty$, we have
    $S_\Gamma(\delta_\Gamma^*, k, 0) = o(1)$.
    \end{prop}
\begin{proof}
In this proof we write $S(\delta)$ for $S_\Gamma(\delta, k,0)$. By Lemma \ref{L:iteration}, we also have $S(\delta_1) - S(\sqrt{D}) \to 0$, and therefore using also Proposition \ref{prop: bound S meso case} we obtain that $S(\delta_1) = o(1)$. Furthermore, by definition we have
\begin{equation}\label{eq:J}
4\frac{\log\log n}{\log(b(\Gamma))} - 2 \leq J  \leq  4\frac{\log\log n}{\log(b(\Gamma))} - 1,
\end{equation}
and it follows from the lower bound of \eqref{eq:J} and the definition of $\delta_\Gamma^*$ that
\begin{equation}
    \delta_{J} \geq D n^{-\sqrt{b(\Gamma)}/\log n} = D e^{-\sqrt{b(\Gamma)}} = \delta_\Gamma^*.
\end{equation}
It is therefore enough to prove that
\begin{equation}
    S(\delta_J) - S(\delta_1) = o(1).
\end{equation}
Now observe that for $1 \leq j \leq J-1$, from Lemma \ref{L:iteration} and some elementary computations,
\begin{equation}
    S(\delta_{j+1}) - S(\delta_j) \lesssim \frac{\delta_{j+1}}{D} \pg \frac{D}{\delta_j}\pd^{K/b(\Gamma)} = n^{-\phi(j)},
\end{equation}
where as before the implicit constant is uniform in $j$, and
\begin{equation*}
    \phi(j) := \frac{1}{b(\Gamma)^{\frac{j+1}{4}}} - \frac{K}{b(\Gamma)^{\frac{j}{4}+1}}.
\end{equation*}
Hence, for $n$ large enough (which we assume in the following), we have
$
    \phi(j) \geq \tfrac{1}{2b(\Gamma)^{({j+1})/{4}}}
$
for every $1\leq j \leq J-1$.
Then,
\begin{equation}
    S(\delta_J) - S(\delta_1) = \summ{j=1}{J-1} (S(\delta_{j+1}) - S(\delta_j)) \lesssim \summ{j=1}{J-1} n^{-1/\pg 2b(\Gamma)^{(j+1)/4}\pd }.
\end{equation}
Note that from \eqref{eq:J}, we have
\begin{equation}
b(\Gamma)^{J/4} \leq (\log n)/b(\Gamma)^{1/4}.
\end{equation}
Making the change of variables $i = J - j$, we therefore have
\begin{equation}
    \summ{j=1}{J-1} n^{-1/\pg 2b(\Gamma)^{(j+1)/4}\pd } = \summ{i=1}{J-1} n^{-1/\pg 2b(\Gamma)^{J/4 -(i-1)/4}\pd }\leq \summ{i=1}{J-1} e^{-b(\Gamma)^{i/4}/2} = \summ{i=1}{J-1} u_i,
\end{equation}
with $u_i := \exp({-b(\Gamma)^{i/4}/2})$.
Note also that for $n$ sufficiently large, we have for all $1 \le i \le J-1$ that
\begin{equation}
\frac{u_{i+1}}{u_i} \leq 1/e.
\end{equation}
The sum $\summ{i=1}{J-1} u_i$ being sub-geometric, it is of the same order of magnitude as its first term, $u_1$, which tends to 0 since $b(\Gamma) \to \infty$.
This concludes the proof.
\end{proof}

The combination of Propositions \ref{prop: bound S micro case}, \ref{prop: bound S meso case}, and \ref{P:bootstrap} shows that $(\delta_\Gamma^*)$ satisfies Assumption (b) of Lemma \ref{lem: good sequence delta implies Gumbel fluctuations}, and can be rewritten as follows.

\begin{prop}\label{P: bound S delta star k 0 combined}
Let $(\Gamma)$ be a collection of finite vertex-transitive graphs of fixed degree satisfying \eqref{e:DC}. 
     Let $k\geq 2$. As $|\Gamma| \to \infty$, we have $S_\Gamma(\delta_\Gamma^*, k, 0) = o(1)$.
\end{prop}

We can now prove Proposition \ref{prop: DC implies Gumbel fluctuations}.

\begin{proof}[Proof of Proposition \ref{prop: DC implies Gumbel fluctuations}]
As $n = |\Gamma|\to \infty$ the following holds. By definition of $\delta_\Gamma^*$,  we have $\delta_\Gamma^* = o(\Diam(\Gamma))$. Moreover, for each $k\geq 2$, we have $S(\delta_\Gamma^*, k,0) = o(1)$ by Lemma \ref{P: bound S delta star k 0 combined}, and 
\begin{equation}
         \min_{A \in \cA(\Gamma, k) \backslash \cA_{\delta_\Gamma^*}(\Gamma, k)} q_A \geq  k-o\pg \frac{1}{\log n} \pd
\end{equation} 
by Proposition \ref{prop: lower bound capacities high mindist assumption c}.
Therefore $(\delta_\Gamma^*)$ satisfies the three assumptions of Lemma \ref{lem: good sequence delta implies Gumbel fluctuations}. This concludes the proof.
\end{proof}

\section{Law of the uncovered set}
\label{S:uncoveredset}

The goal of this section is to show that the information on the law of the cover time can be supplemented by a precise description of the law of the uncovered sets before the cover time. From the results in Section \ref{S:main} and in particular from \eqref{eq:poisson_conv}, we know that at a time $\tsav = \bbE_\pi T_o (s + \log n)$, the \emph{size} of the uncovered set converges to a Poisson random variable with parameter $e^{-s}$.  We will then turn to describe the geometry of this uncovered set: roughly, we aim to show that the uncovered points are approximately uniformly chosen from the vertex set of the graph.

There are different ways to express this idea. The first one is to consider the uncovered set at time $\tsav(1+o(1/\log n))$ (under \eqref{e:DC} this range includes $\tsav$ naturally, but also $\tsmax$ by Proposition \ref{P:thitmaxav} and $\tsexp$ by Proposition \ref{P:tdoubleprime s}) and show that its distribution approximately follows a product structure, as stated in Theorems \ref{T:uncovered_main} and \ref{T:uncovered intro t hit av}. The main goal of this section is to prove the first part of Theorems \ref{T:uncovered_main} and \ref{T:uncovered intro t hit av}: that \eqref{e:DC} implies a product structure for the uncovered set. We restate it as the following proposition.
\begin{prop}\label{prop: theorem intro uncovered set first part: DC implies decorrelation}
    Let $\cF = (\Gamma)$ be a collection of finite (connected) vertex-transitive graphs of fixed degree satisfying \eqref{e:DC}. Let $s\in\bbR$. Let $(\varepsilon_\Gamma)_{\Gamma \in \cF}\in \bbR^\cF$ such that $\varepsilon_\Gamma = o(1/\log |\Gamma|)$ as $|\Gamma|\to \infty$. Recall that $\tsav = \bbE_\pi T_o (\log (|\Gamma|) + s)$. Denote the product over all the vertices of the graph of the Bernoulli law $\mu_s$ with parameter $e^{-s}/|\Gamma|$ by
$\mu_s^{\otimes \Gamma}$. Then, for the simple random walk on $\Gamma$ (under $\P_o$ or $\P_\pi$),
\begin{equation}
    \dtv ( \cL(U(\tsav(1+\varepsilon_\Gamma))), \mu_s^{\otimes \Gamma} ) \xrightarrow[|\Gamma|\to\infty]{} 0.
\end{equation}
In particular, recalling that $\tsmax = \thit(\log(|\Gamma|)+s)$ and that $\tsexp$ is such that $\E(|U(\tsexp)|) = e^{-s}$, we also have $\dtv ( \cL(U(\tsmax)), \mu_s^{\otimes \Gamma} ) \xrightarrow[|\Gamma|\to\infty]{} 0$ and $\dtv ( \cL(U(\tsexp)), \mu_s^{\otimes \Gamma} ) \xrightarrow[|\Gamma|\to\infty]{} 0$.
\end{prop}

Another way is to still look at time $\tsav$, and condition on the size $Z_s = |U(\tsav)| = k$ of the uncovered set at this time. Then we would want to show that the law of the set $U(\tsav)$ itself is close (in the total variation sense) to a uniformly chosen set of size $k$. The proofs of these two variants are both carried out in Section \ref{SS:uncovered_givensize}.

A third way is to consider the stopping time $\tau_k$ which is the first time at which the size of the uncovered set is equal to $k$ (so $\tau_0 = \taucov$), and prove the same approximate uniformity. This is carried out in Section \ref{SS:lastkpoints}.

As we will see, another form of convergence (namely, convergence in total variation) for the cover time itself will follow relatively quickly from these results. This will be explained in Section \ref{SS:tv_conv}.

\medskip

In this whole section, $\cF = (\Gamma)$ is a collection of finite connected vertex-transitive graphs of same degree $d$ that satisfies \eqref{e:DC}.
As before, we write $n = |\Gamma|$ and $D = \Diam(\Gamma)$. We also recall from \eqref{eq: definition of delta star} the definition of $\delta_\Gamma^*$. For convenience, let $\cA^* = \cA^*(\Gamma,k)$ denote the set of all \textit{unordered} subsets $A \subset \Gamma$ (i.e.\ $A$ is a set of vertices of $\Gamma$) of size $k\geq 1$ such that $\mindist (A) > \delta_\Gamma^*$, and
\begin{equation}
    \cA^*_o(\Gamma, k) := \ag A \in \cA^*(\Gamma, k) \mid d(o,A) > \delta_\Gamma^* \ad.
\end{equation}
We emphasize that while $\cA^*(\Gamma,k)$ corresponds to the same subsets as $\cA(\Gamma, k) \backslash \cA_{\delta_{\Gamma}^*}(\Gamma, k)$ from Section \ref{S:main}, its sets are viewed as being \textit{unordered}.

We also recall from Proposition \ref{prop: lower bound capacities high mindist assumption c} that (since $\cF$ satisfies \eqref{e:DC}) for each $k\geq 2$, we have $q_A \geq k - o\pg \frac{1}{\log n}\pd$ uniformly over sets $A\in \cA^*(\Gamma, k)$, i.e.\ sets $A$ of size $k$ such that $\mindist(A) > \delta_\Gamma^*$.
Finally, all asymptotic statements are as $n = |\Gamma| \to \infty$.

\subsection{Convergence to a product measure}
\label{SS:uncovered_givensize}

In this subsection, we prove the first form of uniformity for the uncovered set mentioned above. We start by introducing a second intermediate time.
Given a vertex-transitive graph $\Gamma$ such that $n = |\Gamma| \geq D^2\log n$, 
(and with degree $\geq 2$), we recall that $b(\Gamma) = n/(D^2\log n)$ and set
\begin{equation}\label{eq: definition of t Gamma double star}
   t_\Gamma^{**} := n/\sqrt{b(\Gamma)} = D^2 (\log n)\sqrt{b(\Gamma)}.
\end{equation}

\begin{lemma}\label{lem: good properties of t double star}
    Let $(\Gamma)$ be a collection of finite connected vertex-transitive graphs of fixed degree $d\geq 2$ satisfying \eqref{e:DC}. As $n=|\Gamma| \to \infty$,
   \begin{equation}
      \tmix(1/e) (\log n) b(\Gamma)^{1/9} \lesssim t_\Gamma^{**} \ll n.
   \end{equation}
\end{lemma}
\begin{proof}
    We have $b(\Gamma) \to \infty$ since \eqref{e:DC} holds, so  $t_\Gamma^{**} = n / \sqrt{b(\Gamma)} =  o(n)$. Recall (for instance from the proof of Proposition \ref{P:thitmaxav}) that there exists a constant $C = C(d)$ such that if $n\leq D^8$ then $\tmix(1/e)\leq CD^2$, and if $n>D^5$ then $\tmix(1/e) \leq CD^3$. 
    In the first case we immediately obtain
    \begin{equation}
        t_\Gamma^{**} \gtrsim \tmix(1/e) (\log n) \sqrt{b(\Gamma)} \geq \tmix(1/e) (\log n) b(\Gamma)^{1/9},
    \end{equation}
    and in the second case we have
\begin{equation}
    n/\sqrt{b(\Gamma)} \gtrsim \sqrt{n} \geq D^3 n^{1/8} \gtrsim  D^3 (\log n) n^{1/9} \gtrsim \tmix(1/e) (\log n) b(\Gamma)^{1/9}.
\end{equation}
This concludes the proof.
\end{proof}
\begin{lemma}\label{lem: alpha x close to pi x for vertices x far away from A}
    Let $\cF = (\Gamma)$ be a collection of finite connected vertex-transitive graphs of fixed degree satisfying \eqref{e:DC}. Let $k\geq 1$. As $n=|\Gamma| \to \infty$, uniformly over all $A \in \cA^*_o(\Gamma, k)$, we have
    \begin{equation}
        \frac{\alpha_A(o)}{\pi(o)} = 1+o(1).
    \end{equation}
\end{lemma}
\begin{proof}
    Let $\Gamma\in \cF$ and $A\in \cA^*_o(\Gamma, k)$. Write $\alpha = \alpha_A$. By \cite[Lemma 4.1]{BerestyckiHermonTeyssier2025AB}, we have the following spectral decomposition: for any $t\geq 0$,
\begin{equation}
    \P_o[T_{A}>t]=\sum_{i=1}^m c_{i}f_i(o)e^{-\lambda_it},
\end{equation}
where $m = |A^c|$, $\lambda_1 = 1/\bbE_\alpha T_A$, $c_1^{-2} = \| \alpha/\pi\|_2^2 = \sum_{x\in \Gamma}\alpha(x)^2/\pi(x)$, $f_1 = \frac{\alpha/\pi}{\| \alpha/\pi\|_2}$, $\sum_{i=1}^m c_i^2 = \bbP_\pi(T_A>0) \leq 1$, and $\lambda_i \geq 1/\trel$ for $i\geq 2$.
Recall also from \eqref{treldiff} that $\trel \leq dD^2$, so that by \cite[Lemma 10 (b)]{AldousBrown1992} we have 
\begin{equation}\label{eq: c one asymp form}
    c_1^2 = \frac{1}{\| \alpha/\pi\|_2^2} = 1 + O\pg \frac{\trel}{\bbE_{\alpha} T_A} \pd= 1+ O\pg \frac{D^2}{n}\pd = 1+o(1).
\end{equation}
Moreover, by Corollary \ref{proppxygrand} and a union bound, we have $\P_o[T_{A}>t_\Gamma^{**}] = 1-o(1)$. 
Also, $\pi(o)f_i(o)^2 \leq \sum_{x\in \Gamma} \pi(x) f_i(x)^2$, so $|f_i(o)| \leq \sqrt{n}$ for every $1\leq i\leq m$. It follows, since $t^{**}_\Gamma \gg \tmix(1/e) \log n$ by Lemma \ref{lem: good properties of t double star}, that  
\begin{equation}
    1-o(1) = \P_o[T_{A}>t_\Gamma^{**}]=\sum_{i=1}^m c_{i}f_i(o)e^{-\lambda_it} = f_1(o)(1+o(1)) + n^{-7k} = f_1(o)(1+o(1)), 
\end{equation}
where, in the penultimate step, we used that $\lambda_i \geq 1/\trel$ for $i\geq 2$.
We conclude that
\begin{equation}
   \alpha(o)/\pi(o) = f_1(o) \| \alpha/\pi\|_2 = 1 + o(1). \qedhere
\end{equation}
\end{proof}

Let us now show that starting from a vertex $o$ or from the uniform distribution $\pi$ does not have a significant impact on hitting probabilities of \textit{nice} sets.
\begin{lemma}\label{lem: P o P pi equivalent k s}
Let $\cF = (\Gamma)$ be a collection of finite connected vertex-transitive graphs of fixed degree satisfying \eqref{e:DC}. Let $k\geq 1$ and $s\in\bbR$. Let $(\varepsilon_\Gamma)_{\Gamma \in \cF}\in \bbR^\cF$ such that $\varepsilon_\Gamma = o(1/\log |\Gamma|)$ as $|\Gamma|\to \infty$. 
As $n=|\Gamma| \to \infty$, the following holds.
\begin{enumerate}
    \item Uniformly over all $A \in \cA^*(\Gamma, k)$, we have
    \begin{equation}
        \bbP_\pi(T_A > \tsav(1+\varepsilon_\Gamma)) =  n^{-k}e^{-ks}(1+o(1)).
    \end{equation}
    \item Uniformly over all $A \in \cA^*_o(\Gamma, k)$, we have
\begin{equation}
    \bbP_o(T_A > \tsav(1+\varepsilon_\Gamma)) =  n^{-k}e^{-ks}(1+o(1)).
\end{equation}
\end{enumerate}
\end{lemma}
\begin{proof}
    Denote for this proof $\Tilde{t}_s = \tsav(1+\varepsilon_\Gamma)$. The approximation $\bbP_\pi\pg T_A > \Tilde{t}_s \pd = n^{-k}e^{-ks} (1+o(1))$ follows from Propositions \ref{propapproxpita} and \ref{prop: lower bound capacities high mindist assumption c}. This proves (a). To prove (b) we use the spectral decomposition as in Lemma \ref{lem: alpha x close to pi x for vertices x far away from A}, and follow the notation from its proof. By Theorem \ref{thm: improvement of AB for finite sets under DC} we have $\frac{\bbE_\pi T_A}{\bbE_\alpha T_A} = 1+o(1/\log n)$. Moreover by Proposition \ref{prop: lower bound capacities high mindist assumption c} we have $q_A =\frac{\bbE_\pi T_o}{\bbE_\pi T_A} = k+ o(1/\log n)$. It follows that
\begin{equation}
    \lambda_1 \bbE_\pi T_o = \frac{\bbE_\pi T_A}{\bbE_\alpha T_A} \frac{\bbE_\pi T_o}{\bbE_\pi T_A} = k+o\pg \frac{1}{\log n}\pd.
\end{equation}
We deduce that $\lambda_1 \Tilde{t}_s = k(\log n + s + o(1))$ and conclude, using also Lemma \ref{lem: alpha x close to pi x for vertices x far away from A} (and recalling also \eqref{eq: c one asymp form}), that
\begin{equation}
    \P_o[T_{A}>\Tilde{t}_s] = c_1 f_1(o) n^{-k} e^{-ks}(1+o(1)) + O\pg n^{-7k}\pd = n^{-k} e^{-ks}(1+o(1)). \qedhere
\end{equation}
\end{proof}
In what follows, given $\Gamma$ and an integer $k\geq 1$, we write
\begin{equation}
\cB_k = \ag A\subset \Gamma \du |A| = k, A \notin \cA^*_o(\Gamma, k)\ad.
\end{equation}
\begin{lemma}\label{L:non-macro sets negligible} Let $(\Gamma)$ be a collection of finite connected vertex-transitive graphs of fixed degree satisfying \eqref{e:DC}. Let $k\geq 1$ and $s\in\bbR$. Let $(\varepsilon_\Gamma)_{\Gamma \in \cF}\in \bbR^\cF$ such that $\varepsilon_\Gamma = o(1/\log |\Gamma|)$ as $|\Gamma|\to \infty$. 
As $n=|\Gamma| \to \infty$, we have
\begin{equation}
    \sum_{A\in \cB_k}\bbP_o(U(\tsav(1+\varepsilon_\Gamma)) = A) = o(1).
\end{equation}
\end{lemma}
\begin{proof}
Denote again $\Tilde{t}_s = \tsav(1+\varepsilon_\Gamma)$. Let $A \in \cB_k$.
First we have the trivial bound $\bbP_o(U(\Tilde{t}_s) = A) \leq \bbP_o(T_A > \Tilde{t}_s)$. Denote by $\mu$ the distribution of the walk by time $t_\Gamma^{**}$, conditioned on  starting at $o$. By Lemma \ref{lem: good properties of t double star}, we have $\dtv(\mu, \pi) \leq n^{-7k}$ (for $n$ large enough). Therefore by the Markov property we have (for $n$ large enough)
\begin{equation}\label{eq: L:non-macro sets negligible eq inter}
    \bbP_o(T_A > \Tilde{t}_s) \leq n^{-7k} +\bbP_\pi(T_A>\Tilde{t}_s - t_\Gamma^{**}).
\end{equation}
But since $t_\Gamma^{**} = o(n)$, we have $\Tilde{t}_s - t_\Gamma^{**} = \tsav(1+\varepsilon_\Gamma')$ where $\varepsilon_\Gamma' = \varepsilon_\Gamma - \frac{t_\Gamma^{**}}{\tsav} = o(1/\log n)$.
Summing $\bbP_\pi(T_A>\tsav - t_\Gamma^{**})$ over all sets $A\in \cB_k$, we get
\begin{equation}
\begin{split}
    & \sum_{A \in \cB_k} \bbP_\pi(T_A>\tsav - t_\Gamma^{**}) \\
    \leq  & \, \sum_{A\subset \Gamma \du |A| = k, \mindist(A) \leq \delta_\Gamma^*} \bbP_\pi(T_A>\tsavminuslittleoofone) + \sum_{A \in \cA^*(\Gamma, k) \du d(o,A) \leq \delta^*_\Gamma}\bbP_\pi(T_A>\tsavminuslittleoofone) \\
    \leq & \, S_\Gamma(\delta_\Gamma^*, k, s-o(1)) + o(n^k) n^{-k}.
\end{split}   
\end{equation}
Moreover, by Proposition \ref{P: bound S delta star k 0 combined} and \eqref{eq: good sequence delta implies Gumbel fluctuations inter 1}, we have  $S_\Gamma(\delta_\Gamma^*, k, s-o(1)) \lesssim S_\Gamma(\delta_\Gamma^*, k, 0) = o(1)$. We conclude that
\begin{equation}
     \sum_{A \in \cB_k} \bbP_o(T_A > \tsav(1+\varepsilon_\Gamma)) \leq n^{-6k} + S_\Gamma(\delta_\Gamma^*, k, s-o(1)) + o(n^k) n^{-k} = o(1). \qedhere
\end{equation}
\end{proof}

\begin{lemma}\label{lem: estimate probability uncovered set is A}
Let $(\Gamma)$ be a collection of finite connected vertex-transitive graphs of fixed degree satisfying \eqref{e:DC}. Let $k\geq 1$ and $s\in\bbR$. Let $(\varepsilon_\Gamma)_{\Gamma \in \cF}\in \bbR^\cF$ such that $\varepsilon_\Gamma = o(1/\log |\Gamma|)$ as $|\Gamma|\to \infty$.
As $n=|\Gamma| \to \infty$, uniformly over all $A \in \cA^*_o(\Gamma, k)$, and writing $\Tilde{t}_s = \tsav(1+\varepsilon_\Gamma)$, we have
 \begin{equation}
        \bbP_o(U(\Tilde{t}_s) = A) = \frac{e^{-e^{-s}}e^{-ks}}{n^k} (1+ o(1)).
    \end{equation}
\end{lemma}

\begin{proof}
We cannot directly apply here the moments method (or factorial moments method) as we did in Section \ref{S:main} but we note that there is a relatively simple way to use the work done in this section nevertheless, by exploiting instead the Bonferroni inequalities.

Let us fix a set $A$ as in the lemma, and observe that we can rewrite
\begin{equation}
\begin{split}
     \bbP_o(U(\Tilde{t}_s) = A) & = \bbP_o(A \text{ is not touched but all the points in } \Gamma\backslash A \text{ are}) \\
    & = \bbP_o\pg \ag T_A > \Tilde{t}_s\ad \cap \bigcap_{x\in \Gamma\backslash A}\ag T_x \leq \Tilde{t}_s\ad \pd.
\end{split}
\end{equation}
For $X \subset \Gamma\backslash A$, set $E_X := \ag T_{A\cup X}>\Tilde{t}_s\ad$, with the standard abuse of notation when $X$ is a singleton.
Then
\begin{equation}
\begin{split}
    \bbP_o(U(\Tilde{t}_s) \ne A) & = \bbP_o\pg \ag T_A \leq \Tilde{t}_s\ad \cup \bigcup_{x\in \Gamma\backslash A}\ag T_x > \Tilde{t}_s\ad \pd \\
    & = \bbP_o\pg T_A \leq \Tilde{t}_s \pd + \bbP_o\pg \ag T_A > \Tilde{t}_s\ad \cap \bigcup_{x\in \Gamma\backslash A}\ag T_x > \Tilde{t}_s\ad \pd \\
    & = \bbP_o\pg T_A \leq \Tilde{t}_s \pd + \bbP_o\pg \bigcup_{x\in \Gamma\backslash A}E_x \pd.
    \end{split}
\end{equation}
By Lemma \ref{lem: P o P pi equivalent k s}, uniformly over sets $A\in\cA^*_o(\Gamma, k)$,
\begin{equation}
\bbP_o(T_A > \Tilde{t}_s) =(1 + o(1)) n^{-k}e^{-ks}.
\end{equation}
Consequently, we have
\begin{equation}
    \bbP_o(U(\Tilde{t}_s) = A) = (1 + o(1)) \frac{e^{-ks}}{n^k}- \bbP_o\pg \bigcup_{x\in \Gamma\backslash A}E_x \pd.
\end{equation}
Observe that $E_X \cap E_Y = E_{X\cup Y}$.
Recall also from Proposition \ref{prop: lower bound capacities high mindist assumption c} that if $\mindist (A\cup X) \ge \delta_\Gamma^*$ and $|X|=j$, then $q_{A \cup X} = k + j - o(1/\log n)$. Therefore,
\begin{equation}\label{eq:convmacro}
n^k\sum_{X \subset \Gamma\backslash A \du |X| = j, A\cup X \in \cA^*_o(\Gamma, k+j) }\bbP_o\pg E_X\pd \to \frac{e^{-(k+j)s}}{j!}.
\end{equation}
Choose $J = J (k, \eps, s)$ such that
\begin{equation*}
\sum_{j=J+1}^\infty \frac{e^{(k+j)s}}{(k+j)!}\le \eps.
\end{equation*}
Having chosen $J$, proceeding as in the proof of Lemma \ref{L:non-macro sets negligible}, we can ignore sets $X$ of size $j$ such that $A\cup X \in \cB_{k+j}$, and write for every $1\leq j \leq J+1$,
\begin{equation}
n^k
 \sum_{X \subset \Gamma\backslash A\du |X| = j,  A\cup X \in \cB_{k+j} }\bbP\pg E_X\pd
 \leq \frac{\eps}{J}.
\end{equation}
We deduce from the Bonferroni inequalities and separating according to whether  $A\cup X$ is in $\cA^*_o(\Gamma, k+j)$ or in $\cB_{k+j}$ that for $n$ large enough, 
\begin{equation}
    \begin{split}
& n^k\bg\bbP_o\pg\bigcup_{x\in \Gamma\backslash A} E_x \pd - \sum_{j=1}^J (-1)^{j+1}\sum_{X \subset \Gamma\backslash A, |X| = j}\bbP_o\pg E_X \pd\bd \\ &\leq \sum_{X \subset \Gamma\backslash A: |X|=J+1 } \P_o(E_X) \le \frac{\eps}{J} + \frac{e^{(k+J+1)s}}{(k+J+1)!} + \eps \le 3\eps.        
    \end{split}
\end{equation}
Consequently,
\begin{equation}
     n^k\bg\bbP_o\pg \bigcup_{x\in \Gamma\backslash A}E_x \pd - \sum_{j=1}^J (-1)^{j+1}\sum_{X \subset \Gamma\backslash A, |X| = j, A\cup X \in \cA^*_o(\Gamma, k+j)}\bbP_o\pg E_X \pd\bd \leq 3\eps + J \frac{\eps}{J} = 4 \eps.
\end{equation}
By \eqref{eq:convmacro}, we may assume that $n$ is large enough that for $1\leq j \leq J$,
\begin{equation}
\bg n^k\pg \sum_{X \subset \Gamma\backslash A, |X| = j, A\cup X \in \cA^*_o(\Gamma, k+j)}\bbP_o\pg E_X \pd\pd -  \frac{e^{-(k+j)s}}{j!} \bd \leq \frac{\eps}{J}.
\end{equation}
Combining this with the definition of $J$, and since
\begin{equation}
e^{-ks}- \sum_{j=1}^\infty (-1)^{j+1} \frac{e^{-(k+j)s}}{j!} = e^{-ks}e^{-e^{-s}},
\end{equation}
another triangle inequality finally gives that for $n$ large enough we have
\begin{equation}
    \bg n^k\bbP_o(U(\Tilde{t}_s) = A) - e^{-ks}e^{-e^{-s}} \bd \leq 6\eps. \qedhere
\end{equation}
\end{proof}

We keep the following conditional statement as an intermediate result.

\begin{thm}\label{T:givensize}
Let $(\Gamma)$ be a collection of finite connected vertex-transitive graphs of fixed degree satisfying \eqref{e:DC}. Let $k\geq 0$ and $s\in\bbR$. Let $(\varepsilon_\Gamma)_{\Gamma \in \cF}\in \bbR^\cF$ such that $\varepsilon_\Gamma = o(1/\log |\Gamma|)$ as $|\Gamma|\to \infty$. Write $\Tilde{t}_s = \tsav(1+\varepsilon_\Gamma)$ and $\Tilde{Z}_s = |U(\Tilde{t}_s)|$.
As $n=|\Gamma| \to \infty$, we have \begin{equation}
\bbP_o(\Tilde{Z}_s = k) = \frac{e^{-e^{-s}}e^{-ks}}{k!} (1+ o(1)),    
\end{equation}
and uniformly over all $A \in \cA^*_o(\Gamma, k)$, we have
\begin{equation}
\bbP_o(U(\Tilde{t}_s) = A | \Tilde{Z}_s = k) =  \frac{k!}{n^k} (1+ o(1)).
\end{equation}
In particular as $|\Gamma|\to \infty$, if $\Unif_k$ denotes the uniform law on subsets of size $k$, and if $U(\Tilde{t}_s|k)$ denotes the law of $U(\Tilde{t}_s)$ conditionally given $\Tilde{Z}_s = k$ for the simple random walk on $\Gamma$ started at $o$, we have
\begin{equation}
    \dtv (U(\Tilde{t}_s|k), \Unif_k) \to 0.
\end{equation}
\end{thm}
\begin{proof}
    The case $k=0$ is the cover time result that we already established in Proposition \ref{prop: DC implies Gumbel fluctuations}. We may therefore assume that $k\geq 1$. First we have $|\cA_o^*(\Gamma, k)| = \frac{n^k}{k!}(1+o(1))$, that is a proportion $1-o(1)$ of subsets of $\Gamma$ of size $k$ are in $\cA_o^*(\Gamma, k)$. Therefore by Lemmas \ref{L:non-macro sets negligible} and \ref{lem: estimate probability uncovered set is A} we have
    \begin{equation}
        \bbP_o(\Tilde{Z}_s = k) = o(1) + \frac{n^k}{k!} \frac{e^{-e^{-s}}e^{-ks}}{n^k} (1+ o(1)) = \frac{e^{-e^{-s}}e^{-ks}}{k!} (1+ o(1)).
    \end{equation}
The second result and the convergence in total variation follow, since for $A \in \cA_o^*(\Gamma, k)$ we have 
\begin{equation}
    \bbP_o(U(\Tilde{t}_s) = A \mid \Tilde{Z}_s = k) = \frac{\bbP_o(U(\Tilde{t}_s) = A)}{\bbP_o( \Tilde{Z}_s = k)} = \frac{k!}{n^k}(1+o(1)). \qedhere
\end{equation}
\end{proof}

We are now ready to prove Proposition \ref{prop: theorem intro uncovered set first part: DC implies decorrelation}.

\begin{proof}[Proof of Proposition \ref{prop: theorem intro uncovered set first part: DC implies decorrelation}] 
Let $\varepsilon>0$, and let $K\geq 1$ such that 
\begin{equation}
    \sum_{k=K+1}^{\infty} \frac{e^{-e^{-s}}e^{-ks}}{k!} \leq \varepsilon.
\end{equation}
For each $k\geq 0$, we have by definition of $\mu_s^{\otimes \Gamma}$, uniformly over all $A\subset \Gamma$ of size $k$,
    \begin{equation}
        \mu_s^{\otimes \Gamma} (A) = \pg\frac{e^{-s}}{n}\pd^{|A|}\pg 1-\frac{e^{-s}}{n}\pd^{n-|A|} = (1+o(1)) \frac{e^{-ks}}{n^k} e^{-e^{-s}}.
    \end{equation}
Since $K$ is fixed, this bound, as well as that of Lemma \ref{lem: estimate probability uncovered set is A}, is uniform over all sets $A\in \cup_{k=0}^{K} \cA^*_o(\Gamma, k)$. Since a proportion $1-o(1)$ of subsets of size $\leq K$ of $\Gamma$ are in $\cup_{k=0}^{K} \cA^*_o(\Gamma, k)$, and the other sets do not contribute by Lemma \ref{L:non-macro sets negligible} (note that their contribution is also $o(1)$ for $\mu_s^{\otimes \Gamma}$ since this depends only on the number of sets), we obtain that under $\bbP_o$, writing again $\Tilde{t}_s = \tsav(1+\varepsilon_\Gamma)$,
\begin{equation}
    \dtv ( \cL(U(\Tilde{t}_s)), \mu_s^{\otimes \Gamma} ) \leq \varepsilon + o(1).
\end{equation}
Therefore $\dtv (\cL(U(\Tilde{t}_s)), \mu_s^{\otimes \Gamma} ) = o(1)$. By vertex-transitivity, we have $\dtv (\cL(U(\Tilde{t}_s)), \mu_s^{\otimes \Gamma} ) = o(1)$ under $\bbP_x$ for each $x\in \Gamma$, and by averaging we obtain the result under $\bbP_\pi$. Finally, the result for $\tsmax$ holds, since $\thit = \thitav(1+o(1/\log n))$ by Proposition \ref{P:thitmaxav}.
\end{proof}

\subsection{Convergence of the last $k$ points}
\label{SS:lastkpoints}

In this section we are interested in the first time at which the uncovered set has size $k$, (where $k\geq 1$ is fixed for this section), that is:
\begin{equation*}
\tau_k :=
\inf \{t\ge 0: | \{X_u, u \le t\}^c | = k\}.
\end{equation*}
We want to show that the distribution of $U(\tau_k)$ is close to the uniform distribution $\Unif_k$ over sets of size $k$.

\begin{thm}\label{T:convergence of the uncovered set}
Let $(\Gamma)$ be a collection of finite connected vertex-transitive graphs of fixed degree satisfying \eqref{e:DC}. Let $k\geq 1$. As $n \to \infty$, we have, under $\bbP_\pi$,
\begin{equation}
\dtv(U({\tau_k}), \Unif_k) \to 0.
\end{equation}
\end{thm}

We will need two ingredients.
Our first lemma improves on Theorem \ref{T:givensize} by showing that the position of the walk at time $\tsav$ is approximately independent of $U(\tsav)$ and is distributed (approximately again) according to $\pi$.

\begin{lemma}
  \label{L:indep}
  Let $(\Gamma)$ be a collection of finite connected vertex-transitive graphs of fixed degree satisfying \eqref{e:DC}. Let $s \in \bbR$. Recall that we denote the product over all the vertices of the graph of the Bernoulli law $\mu_s$ with parameter $e^{-s}/|\Gamma|$ by
$\mu_s^{\otimes \Gamma}$.
  As $|\Gamma|\to \infty$, we have, under $\bbP_\pi$,
  \begin{equation*}
  \dtv[ ( U(\tsav), X_{\tsav}) , ( \mu_s^{\otimes \Gamma} \otimes \pi) ] \to 0.
  \end{equation*}
\end{lemma}
\begin{proof}
We work under $\P = \P_{\pi}$. Recall that $Z_s=| U(\tsav)|$. Let $s' = s'(s,\Gamma) := s-1/\sqrt{f(\Gamma)}$.
Then
\begin{equation}\label{E:condition on s'}
D^2 \ll \tsav - \tsprimav \ll n = D^2f(\Gamma),
\end{equation}
and $s - s' = o(1)$ so that $\E( Z_{s'} - Z_s) \to 0$. As moreover $Z_s \le Z_{s'}$, we deduce that
  \begin{equation}
  \P ( Z_s \neq Z_{s'} ) = \bbP(Z_{s'} - Z_{s} \geq 1) \leq \E( Z_{s'} - Z_s) \to 0.\label{eq:nochange}
  \end{equation}
Furthermore, we already know from Theorem \ref{T:givensize} that $\dtv( U(\tsprimav), \mu_{s'}^{\otimes \Gamma} ) \to 0$. But since $|s- s'| \to 0$, we also have
  $\dtv ( U(\tsprimav), \mu_s^{\otimes \Gamma} ) \to 0$.

To conclude, let us compare $\P( X_{\tsav} = x, U(\tsav) = A) $ with our target $\pi(x) \lambda(A)$, where we set $\lambda (A) = e^{- e^{-s}} e^{-|A|s} / n^{|A|}$, i.e.\ $\lambda= \mu_s^{\otimes \Gamma}$. Then
  \begin{equation}\label{eq: convergence of the uncovered set first time k multiline inter eq}
  \begin{split}
    & \left| \P( X_{\tsav} =x , U(\tsav) =A) - \pi(x) \lambda(A) \right | \\ & \le \left| \P( X_{\tsav} =x , U(\tsav) =A) - \P (X_{\tsav} = x, U(\tsprimav) = A)\right| \\
    & \quad \quad + \left|  \P (X_{\tsav} = x, U(\tsprimav) = A) - \pi(x) \lambda(A)\right|\\
    & \le \P( X_{\tsav} = x, U(\tsprimav) = A, Z_{s} \neq Z_{s'} ) \\
    & \quad \quad + \left| \P (U(\tsprimav) = A) \P (X_{\tsav} = x | U(\tsprimav) =A) - \pi(x) \lambda(A) \right |.  
  \end{split}
  \end{equation}
  Summing over all $x, A$, for the first term in the right hand side of \eqref{eq: convergence of the uncovered set first time k multiline inter eq}, we get
  \begin{equation*}
  \sum_x \sum_A \P( X_{\tsav} = x, U(\tsprimav) = A, Z_{s} \neq Z_{s'} ) = \P( Z_s \neq Z_{s'} ),
  \end{equation*}
  which converges to zero by \eqref{eq:nochange}.
  The second term of the right hand side of \eqref{eq: convergence of the uncovered set first time k multiline inter eq}, on the other hand, can be written by Theorem \ref{T:givensize} as
  \begin{equation*}
  \left |\lambda (A) (1+ o(1))\P( X_{\tsav} = x | U(\tsprimav) =A) - \pi(x) \lambda(A)\right |
  \end{equation*}
  But observe that by choice of $s' $ compared to $s$, in the interval from $\tsprimav$ to $\tsav$ the walk has had time to mix. More specifically, the lower bound of \eqref{E:condition on s'}, $D^2 \ll \tsav - \tsprimav =: t$, implies, through \cite[Proposition 4.15]{LivreLevinPeres2019MarkovChainsAndMixingTimesSecondEdition}, \eqref{E:Poincaré}, and Corollary \ref{Cor:all bounds in one}, that uniformly over $x,y\in \Gamma$,
\begin{equation}
|p_t(x,y) - 1/n|\leq p_t(o,o) -1/n \leq \exp\pg-\frac{t-D^2}{dD^2}\pd p_{D^2}(o,o) = o\pg \frac{1}{n}\pd,
\end{equation}
i.e.\ that $p_t(x,y) = \pi(y)(1+o(1))$.

Thus, by averaging over $ X_{\tsprimav} $ we have $\P( X_{\tsav} = x | U(\tsprimav) =A) = \pi(x) (1+ o(1))$, and we conclude that this second term is equal to
  $
  \lambda (A) \pi(x) o(1)
  $. Summing over all $x $ and $A$, we see that the sum is $o(1)$ and we deduce that
   \begin{equation*}
   \sum_x \sum_A \left| \P( X_{\tsav} =x , U(\tsav) =A) - \pi(x) \lambda(A) \right | \to 0,
   \end{equation*}
   which is the desired result.
  \end{proof}

We also prove a small technical lemma that shows that hitting a set of points from far away takes about the same amount of time as hitting it from stationarity.
\begin{lemma}\label{L:hitunif_fixed prelim}
    Let $\cF = (\Gamma)$ be a collection of finite connected vertex-transitive graphs of fixed degree satisfying $\eqref{e:DC}$, let $(\delta_\Gamma)_{\Gamma \in \cF}$ be such that $1\leq \delta_\Gamma \to \infty$ as $n=|\Gamma|\to \infty$, and let $k\geq 1$. Then as $n\to\infty$, uniformly over all non-empty disjoint subsets $A,B$ of $\Gamma$ of size at most $k$ such that $\mindist(A\cup B) \geq \delta_\Gamma$, we have
    \begin{equation}
        \bbE_B T_A = \pg \frac{1}{|A|}+o(1)\pd\bbE_\pi T_o,
    \end{equation}
    where $\bbE_B = \frac{1}{|B|}\sum_{x\in B} \bbE_x T_A$.
\end{lemma}
\begin{proof}
    Up to averaging, it is enough to prove the result when $B = \ag x \ad$ is a singleton, which we assume for the rest of the proof. Denote by $(X_t)_{t\geq 0}$ the rate-1 simple random walk on $\Gamma$, started at $x$. Write also $t^*$ for $t^*_\Gamma$ as defined in \eqref{eq: definition of t Gamma star}.
    Since $t^* \gg D^2 \asymp \tmix$ (by Lemma \ref{P:tmix asymp diam squared}) we have $\dtv(X_{t^*}, \pi) = o(1)$. Moreover, since $t^*=o(n)$, we have $\max_{y\in \Gamma \du d(x,y) \geq \delta_\Gamma} \bbP_x (T_y\leq t^*) = o(1)$ by Corollary \ref{proppxygrand}. Hence by a union bound we have $\bbP_x(T_A\leq t^*) = o(1)$. 

    Denote by $\mu$ the law of $X_{t^*}$ conditionally on $\ag T_A >t^* \ad$.
    Since $\dtv(X_{t^*}, \pi) = o(1)$ and $\bbP_x(T_A \leq t^*) = o(1)$, we have $\dtv(\mu, \pi) =o(1)$. Moreover, 
    \begin{equation}
         \bbE_x T_A = \int_0^\infty \bbP_x(T_A>t)\mathrm{d}t = \int_0^{t^*} \bbP_x(T_A>t)\mathrm{d}t + \int_{t^*}^\infty \bbP_x(T_A>t)\mathrm{d}t.
    \end{equation}
But $\bbE_x T_A \asymp n$, $\int_0^\infty \bbP_x(T_A>t)\mathrm{d}t  =O(t^*) = o(n)$, $\bbP_x(T_A > t^*) = 1-o(1)$, and for $t\geq t^*$, $\bbP_x(T_A>t) = \bbP(T_A>t^*)\bbP_x(T_A>t\mid T_A>t^*)$ so
\begin{equation}
\begin{split}
     (1+o(1))\bbE_x T_A = \frac{\bbE_x T_A - \int_0^{t^*} \bbP_x(T_A>t)\mathrm{d}t}{\bbP_x(T_A > t^*)}  = \int_{t^*}^\infty \bbP_x(T_A>t\mid T_A>t^*)\mathrm{d}t.
\end{split}
\end{equation}
We obtain that
\begin{equation}
 (1+o(1))\bbE_x T_A = \int_{0}^\infty \bbP_\mu(T_A>t)\mathrm{d}t
     = \bbE_\mu T_A = (1-o(1))\bbE_\pi T_A + o(1)\cdot \max_{z\in \Gamma\backslash A} \bbE_z T_A,     
\end{equation}
so, since $\max_{z\in \Gamma\backslash A} \bbE_z T_A \leq \thit \lesssim n$, we obtain that
\begin{equation}\label{eq: inter E x T A E pi T A}
    \bbE_x T_A = (1+o(1))\bbE_\pi T_A.
\end{equation}
Finally, by Lemma \ref{lem: bound q A meso case} we have $q_A = |A|+o(1)$, that is $\frac{\bbE_\pi T_A}{\bbE_\pi T_o} = \frac{1}{|A|} + o(1)$. Plugging this into \eqref{eq: inter E x T A E pi T A} concludes the proof.
\end{proof}

\begin{lemma}\label{L:hitunif_fixed}
    Let $\cF = (\Gamma)$ be a collection of finite connected vertex-transitive graphs of fixed degree satisfying $\eqref{e:DC}$, and let $k\geq 1$. Uniformly over all $y\in \Gamma$, $A \in \cA^*(\Gamma, k)$ such that $d(y, A) \ge \delta^*$, and $x \in A$, 
the probability that $x$ is the first point of $A$ to be touched by the walk started at $y$ is $1/k + o(1)$:
\begin{equation}
\bbP_y(T_A=T_x) = \frac{1}{k} + o(1).
\end{equation}
\end{lemma}
\begin{proof}
    Denote $B = A\backslash\ag x\ad$. By \cite[Corollary 2.10]{AldousFill}, which we can apply with $(y,x,B)$ in place of $(i,j,\ell)$ up to collapsing the chain at $B$, we have
    \begin{equation}
        \bbP_y(T_A = T_x) = \bbP_y(T_x < T_B) = \frac{\bbE_y T_B + \bbE_B T_x - \bbE_y T_x}{\bbE_x T_B + \bbE_B T_x}.
    \end{equation}
Applying Lemma \ref{L:hitunif_fixed prelim} to all expectations, we conclude that 
 \begin{equation}
        \bbP_y(T_A = T_x) = \frac{1/(k-1)  + 1 - 1 + o(1)}{1/(k-1) + 1 + o(1)} = \frac{1}{k} + o(1). \qedhere
    \end{equation}
\end{proof}

\begin{lemma}\label{L:itérationkpoints}
 Let $\cF = (\Gamma)$ be a collection of finite connected vertex-transitive graphs of fixed degree satisfying $\eqref{e:DC}$, and let $k\geq 1$. Uniformly over all $A \in \cA^*(\Gamma, k)$ and $x \in A$, 
the probability that $x$ is the first point of $A$ to be touched by the walk starting from uniformity is $1/k + o(1)$:
\begin{equation}
\bbP_\pi(T_A=T_x) = \frac{1}{k} + o(1).
\end{equation}
\end{lemma}
\begin{proof}
    Since $\delta_\Gamma^* = o(\Diam(\Gamma))$ and $A$ is finite, we have $\pi(\ag y\in \Gamma \du \ag y\ad \cup A \in \cA^*(\Gamma, k) \ad) = 1-o(1)$. The result then follows from Lemma \ref{L:hitunif_fixed}.
\end{proof}

We can now prove the main result of this subsection.
\begin{proof}[Proof of Theorem \ref{T:convergence of the uncovered set}]
We work under $\P = \P_\pi$ again. Let $k \geq 2$ and $\eps > 0$. Let $s = s(k,\eps)\in \bbR$ be such that for $n$ large enough,
\begin{equation}
\bbP(\tau_k > \tsav) \geq 1 - \eps.
\end{equation}
Let also $K = K(s,k, \eps) \geq k+1$ be such that for all sufficiently large $n$ we have that,
\begin{equation}
\bbP(\tau_K \leq \tsav) \geq 1-\eps.
\end{equation}
In particular, with  probability at least $1-2\eps$, there are between $k+1$ and $K$ uncovered points at time $\tsav$:
\begin{equation}
\bbP(k+1\leq Z_s\leq K) = \bbP(\tau_K \leq \tsav < \tau_k)\geq 1-2\eps.
\end{equation}

Fix $j$ such that $k+1\leq j\leq K$, and condition on the event $Z_s = j$. Let us now describe the evolution of $U(\tsav)$, up to events of probability $o(1)$, in order to get an approximation of $U(\tau_k)$ in the total variation sense. Conditionally on $\{Z_s = j\}$, we know from Theorem \ref{T:givensize} that $U(\tsav)$ is a uniformly chosen set of size $j$, and we may assume it in $\cA^*(\Gamma,j)$. Furthermore, by Lemma \ref{L:indep}, the position of the walk at time $\tsav$ is uniformly distributed on $\Gamma$, independently from $U(\tsav)$.  (Again these descriptions refer in reality to approximation in the total variation sense.)
By Lemma \ref{L:itérationkpoints}, the next point that is removed from $U(\tsav)$ is therefore uniformly chosen among $U(\tsav)$. Applying next Lemma \ref{L:hitunif_fixed} $j-1-k$ times successively, from this point onwards, at each successive stage until time $\tau_k$, points are removed uniformly at random. Since $U(\tsav)$ was uniformly distributed among sets of size $j$ initially, it therefore follows that (still under the conditional law given $\{Z_s = j\}$) that the law of $U(\tau_k)$ is (close to, in the total variation sense) $\Unif_k$. Since this is true for every $k +1 \le j \le K$, we deduce
\begin{equation*}
\dtv (U(\tau_k), \Unif_k) \le 2 \eps + o(1).
\end{equation*}
The result follows.
\end{proof}

\subsection{Convergence in total variation of the cover time}
\label{SS:tv_conv}

We illustrate the results above by strengthening the mode of convergence for the rescaled cover time: namely the distribution $\nu_\Gamma$ of the random variable $Y_\Gamma = \frac{\taucov}{\bbE_\pi T_o} - \log n$, converges in total variation to a standard Gumbel distribution $\nu$.
\begin{thm}\label{T:tvconv}
Let $(\Gamma)$ be a collection of finite connected vertex-transitive graphs of fixed degree satisfying \eqref{e:DC}. As $n = |\Gamma| \to \infty$, we have
\begin{equation}
\dtv(\nu_n, \nu) \to 0.
\end{equation}
\end{thm}

\begin{proof}
We need several ingredients.
First, we saw in Theorem \ref{T:convergence of the uncovered set} that for each $k \geq 1$, $U(\tau_k)$ is approximately uniform, in particular, at time $\tau_1$, the walk is with probability $1-o(1)$ macroscopically far from the last point $x_0$ to be visited, and hence, with the same arguments as for the convergence of the $k$ last points, the walk with probability $1-o(1)$ gets mixed again before touching $x_0$.
Recalling moreover that the distance in total variation between the uniform distribution and the quasi-stationary distribution associated with $x_0$ is $o(1)$, we have that  
\begin{equation*}
\dtv\left( \frac{\tau_0 - \tau_1}{\E_\pi (T_o)}, \mathrm{Exp}(1) \right) \to 0.
\end{equation*}
Let $\nu_{\text{cont}}$ denote the law of $\frac{\tau_1}{\E_\pi (T_o)} + X  $, where $X$ is an independent exponential random variable with mean 1. Note that $\nu_{\text{cont}}$ is obtained by convolution with an exponential law and so has a 1-Lipschitz density with respect to Lebesgue measure. Furthermore, since $\taucov = \tau_0 = \tau_1 + (\tau_0 - \tau_1)$,
\begin{equation}\label{eq:weak}
\dtv\left( \frac{\taucov}{\E_\pi(T_o)}, \nu_{\text{cont}} \right) \to 0.
\end{equation}
Second, we have already shown that for each fixed $s\in\bbR$,
\begin{equation}
    F_\Gamma(s) := \bbP(Y_\Gamma \leq s) = \bbP(Z_s = 0) \to e^{-e^{-s}} =: F(s).
\end{equation}
Let now $\eps >0$ and let us fix $S = S(\eps)$ large enough such that for $n$ large enough, we have $F_\Gamma(S) \geq 1-\eps$  and $F_\Gamma(-S) \leq \eps$. In particular, for such $n$’s, we have
\begin{equation}
\nu_\Gamma(\cg -S,S \cd) = \bbP(Y_\Gamma\in \cg -S, S \cd) \geq 1- 2\eps.
\end{equation}
Let $f_Y$ be the density of $Y$ (which is just the density of the Gumbel law) and let $f_n$ denote the density of $\nu_{\text{cont}}$ after translating by $\log n$. Since $\dtv( \nu_\Gamma, f_n) \to 0$ it suffices to show that
\begin{equation}\label{goalL1}
\dtv(f_n, f_Y) = (1/2) \int_{s\in \bbR} |f_n(s) - f_Y(s)| \mathrm{d}s  \to 0.
\end{equation}
Observe that $f_{n}$ is 1-Lipschitz over $[-S, S]$, since it is a convolution of some given law with an exponential law. It is also pointwise bounded at, say $s = 0$ (indeed, since $f_n$ is Lipschitz, it cannot be large at any point without its integral being large, which is not possible by \eqref{eq:weak}). It is therefore uniformly equicontinuous, and by the Ascoli--Arz\'ela theorem has subsequential uniform limits. However, the limit can only be $f_Y$, again by \eqref{eq:weak}. Thus $f_n$ converges to $f_Y$ uniformly over $[-S, S]$, and hence also in the $L^1$ sense over $[-S, S]$. This proves \eqref{goalL1}.
\end{proof}

\section{When the diameter condition fails}\label{S:subcritical cases}

\report{In this section we complete the proofs of Theorems \ref{T:main}, \ref{T:uncovered_main}, and \ref{T:uncovered intro t hit av}.
To do so we start by extending the gradient inequality of Diaconis and Saloff-Coste \cite{DSC} from Cayley graphs to vertex-transitive graphs in Proposition \ref{P:gradient}, and then prove bounds on Green functions (defined in \eqref{eq: Green function def}). In particular we quantify the fact that occupation measures (and hence Green functions) are positively correlated at short (but still macroscopic) distances, and negatively correlated at larger distance, see Proposition \ref{P:greenfunctionmacroscopically continuous} and Proposition \ref{P:greenfunctiontwoballsmacrofar}.
 This allows us to prove that when \eqref{e:DC} fails, there is no product structure for the uncovered set.}

\subsection{Gradient inequality and Green function estimates}

In this section $\Gamma = (V,E)$ is a vertex-transitive graph of degree $d$. Denote the transition matrix of simple random walk on $\Gamma$ by $P$ and the time $t$ transition probabilities for the rate 1 continuous-time simple random walk on $\Gamma$ by $P_t=e^{-t(I-P)}$ (so that we have $P_t(x,y) = p_t(x,y)$). 
Recall that for $x,y\in V$ and $t\geq 0$,
\begin{equation}
    h_t(x,y) = p_t(x,y) - \frac{1}{n}.
\end{equation} 
We start by adapting a gradient inequality, due to Diaconis and Saloff-Coste \cite{DSC} for Cayley graphs, to vertex-transitive graphs. The main difference is that we need to use the mass transport principle.

\begin{prop}
\label{P:gradient}
Let $\Gamma=(V,E)$ be a finite vertex-transitive graph of degree $d$. Then for all $s,t \ge 0$ and $o,y,z \in V$ we have that
\begin{equation}
\label{e:DSC1}
|h_{t+s}(o,y)-h_{t+s}(o,z)| \le d (y,z) \sqrt{\frac{d}{2es}} h_t(o,o).  \end{equation}
\end{prop}
\begin{proof}  
By  the triangle inequality it suffices to prove the inequality when $y$ and $z$ are neighbours. Let then $y,z \in V$ such that $y \sim z$. First observe that
\begin{align*}
    h_{t+s}(o,y)-h_{t+s}(o,z) & = p_{t+s}(o,y)-p_{t+s}(o,z) \\
    & = \sum_{w\in V} p_{t/2}(o,w)(p_{t/2+s}(w,y)-p_{t/2+s}(w,z)) \\
    & = \sum_{w\in V} h_{t/2}(o,w)(h_{t/2+s}(w,y)-h_{t/2+s}(w,z))
\end{align*}
because, by reversibility, $\sum_{w\in V} (h_{t/2+s}(w,y)-h_{t/2+s}(w,z)) = 0$.
By reversibility again, we have 
\begin{equation}\label{e:hproduitscalaire}
    \sum_{w\in V} h_{t/2}(o,w)^2 = \pg \sum_{w\in V} p_{t/2}(o,w)p_{t/2}(w,o)\pd - \frac{1}{n} = p_t(o,o) - \frac{1}{n} = h_t(o,o).
\end{equation}
It follows from the triangle inequality and the Cauchy--Schwarz inequality that
\begin{equation}
\begin{split}
\label{e:DSC2}
 |h_{t+s}(o,y)-h_{t+s}(o,z)|^2  & \le \pg\sum_{w} h_{t/2}(o,w)^2\pd\pg \sum_w |h_{t/2+s}(w,y)-h_{t/2+s}(w,z)|^{2}\pd \\ & =  h_{t}(o,o) \sum_w |p_{t/2+s}(y,w)-p_{t/2+s}(z,w)|^{2} \\ & \leq h_{t}(o,o) \sum_w \sum_{y'\du y'\sim y}|p_{t/2+s}(y,w)-p_{t/2+s}(y',w)|^{2} \\
 & = h_t(o,o) \sum_w F(y,w),\end{split}
\end{equation} 
where
\begin{equation*}
F \du(a,b) \mapsto \sum_{x \du x\sim a} \bg p_{t/2 + s}(a,b) -p_{t/2 + s}(x,b)\bd^2.
\end{equation*}

Observe that for any $\gamma \in \Aut(\Gamma)$, and $a,b \in V$, we have $F(\gamma(a), \gamma(b)) = F(a,b)$. Moreover $\Aut(\Gamma)$ is a discrete group of automorphisms, so by \cite[Corollary 8.9]{LyonsPeresBook} is unimodular. Since $\Gamma$ is vertex-transitive, $\Aut(\Gamma)$ is (by definition) transitive, so we can apply the mass transport principle. By \cite[Equation (8.4)]{LyonsPeresBook}, we hence have, for all $y \in V$, that
\begin{equation}\label{e:DSC3}
    \sum_{w\in V} F(y,w) = \sum_{w\in V} F(w,y)  = \sum_{w\in V} \sum_{w'\du w' \sim w} |p_{t/2+s}(w,y)-p_{t/2+s}(w',y)|^{2}. 
\end{equation}
Denote $P_{\mathrm{L}}:=\frac 12(I+P)$. 
For $f,g:V \to \mathbb{R}$ denote $\langle f,g\rangle=\pi(o) \sum_{v \in V}f(v)g(v)$ and $\|f\|_2^2:=\langle f,f\rangle$.
Since $I-P_{\mathrm{L}}$ is self-adjoint and $\langle (I-P_{\mathrm{L}}) f,f\rangle \ge 0 $ for all $f: V \to \mathbb{R}$, we may consider $K:=\sqrt{I-P_{\mathrm{L}}}$ which is a self-adjoint operator satisfying $K^2=I-P_{\mathrm{L}}$, and $\langle (I-P_{\mathrm{L}})f,g\rangle=\langle Kf,Kg\rangle$ for all $f,g:V \to \bbR$.
Noting that $P_{s}=e^{-2s(I-P_{\mathrm{L}})}$ we have that $KP_{s}=q(I-P_{\mathrm{L}})$, where $q:[0,1] \to [0,1]$ is defined by $q(u)=\sqrt{u}e^{-2su}$, notation which we also extend (slightly abusing notation) to matrices.
Since the spectral measure of  $I-P_{\mathrm{L}}$ is supported on $[0,1]$,   it is standard that $ \|KP_s \|_2^2 =\|q(I-P_{\mathrm{L}}) \|_{2}^2 \le \max_{u \in [0,1]}q(u)^2=\frac{1}{4es}$, where  $  \|KP_s \|_2 :=\sup_{f \in \mathbb{R}^V:\|f\|_2=1 } \|KP_sf \|_2  $ is the operator norm of $KP_s$. 

For $r \ge 0$ and $w\in V$, denote $f_{r}(w)= h_r(w,y)$. By reversibility, the sum on the right hand side of \eqref{e:DSC3} equals $2d$ times
(see e.g., \cite[Lemma 13.6]{LivreLevinPeres2019MarkovChainsAndMixingTimesSecondEdition})
\begin{equation*}
\begin{split}
    \frac{1}{\pi(o)}\langle\left(I-P_{\mathrm{L}}\right)f_{t/2+s},f_{t/2+s}\rangle =\frac{1}{\pi(o)}\langle K^{2}P_{s}f_{t/2},P_{s}f_{t/2}\rangle \leq\frac{1}{\pi(o)}\|KP_s \|_{2}^2 \|f_{t/2}\|_2^2 \le \frac{h_{t}(o,o)}{4es}, 
\end{split}
\end{equation*}
where in the last inequality we used that $\|KP_s \|_2^2\le\frac{1}{4es}$ and \eqref{e:hproduitscalaire}.

It follows that the right hand side of \eqref{e:DSC3} is less than $\frac{d h_{t}(o,o)}{2es}$, which in conjunction with \eqref{e:DSC2} concludes the proof.
  \end{proof}

We now prove that at times $t$ proportional to $D^2$ the distribution of the walk is still far from being uniform.

\begin{prop}\label{P:lowerbounddiagonalheatkernel}
Assume that $n \leq D^3$. There exists a constant $0<a = a(d)<1$ such that for all $ t \leq aD^2$, 
\begin{equation}
nh_t(o,o)\geq 1/2.
\end{equation}
\end{prop}
\begin{proof}
Recall the off-diagonal heat kernel bound \eqref{Folz}. Note that the constant $C$ in \eqref{g}, as well as the constant $c$ and the implicit constant in \eqref{Folz} (call it $C'$ for this proof), depend only on the graph degree, and the growth of the function $g$ defined in \eqref{g}. We can hence take the same constants for all (connected) finite vertex-transitive graphs of degree $d$ such that $D^3\geq n$.

First assume that $D^2 \leq n < D^3$, and let $x$ such that $d(o,x) = D$. We have for every $D\leq t \leq D^2$,
\begin{equation}
    p_t(o,x) \leq C' \frac{1}{g(t)}\exp\pg-c\frac{D^2}{t}\pd = C' \frac{C}{\min(f(\Gamma)t, t^{3/2})} \exp\pg-c\frac{D^2}{t}\pd.
\end{equation}
Let $0<a<1$. At time $t = aD^2$ we have \begin{equation}
    \min(f(\Gamma)t, t^{3/2}) = \min(a n, a^{3/2}D^3) \geq a^{3/2}n.
\end{equation}
It follows that
\begin{equation}\label{e:bound on pt(o,x)}
    p_t(o,x) \leq \frac{CC'}{a^{3/2}} \frac{1}{n} \exp\pg-c\frac{1}{a}\pd.
\end{equation}
Fixing $a = a(c,C,C')$ small enough, we hence have at time $ t = aD^2$ that $np_t(o,x) \leq 1/2$, i.e.\
\begin{equation}\label{e:heatkernelatmaxdistancesmall}
    nh_t(o,x) \leq -\frac{1}{2}.
\end{equation}
Finally, we have
$h_t(o,o) = \max_{z,w \in \Gamma} \bg h_t(z,w)\bd$ by \cite[Proposition 4.15]{LivreLevinPeres2019MarkovChainsAndMixingTimesSecondEdition} and vertex-transitivity.
We conclude that
\begin{equation}
    nh_t(o,o) \geq \frac{1}{2}.
\end{equation}
The proof for $D\leq n < D^2$ is similar.
\end{proof}
   
We define the Green function between two points $x,y\in V$ by
\begin{equation}\label{eq: Green function def}
G(x,y)=\int_{0}^{\infty}h_t(x,y)\mathrm{d}t.    
\end{equation}
The following proposition collects some useful identities involving Green functions.
\begin{prop}\label{P:Green function identities}
Let $\Gamma$ be any finite connected vertex-transitive graph. We have the following identities.
\begin{enumerate}
    \item Fix $o\in \Gamma$. We have
    \begin{equation}\label{e:identity thitav green function}
    \bbE_\pi T_o = nG(o,o) \geq n - (1 + \log n).
\end{equation}
\item For every $x,y \in \Gamma$,
\begin{equation}\label{e:identity ExTy green function}
    \bbE_x T_y = n\pg G(o,o) - G(x,y) \pd.
\end{equation}
\item For any subset $A = \ag x,y \ad$ of size 2 of $\Gamma$, recalling that $q_A = \frac{\E_\pi[T_o]}{ \E_\pi [T_A]}$ by definition, we have
\begin{equation}\label{e:identité capacité}
    q_A = \frac{2}{1+\frac{G(x,y)}{G(o,o)}}.
\end{equation}
\end{enumerate}
\end{prop}
\begin{proof}
The equalities in the first two points are stated in \cite{AldousFill} in a more general framework, as Lemma 2.11 and Lemma 2.12, respectively, and Section 2.2.3 explains why they also hold in continuous time.
We moreover have, lower bounding the probability to be at $o$ by the probability to have never jumped:
\begin{equation}
    G(o,o)\geq \int_0^{\log n} \cg e^{-t} - \frac{1}{n}\cd \mathrm{d}t= 1 - \frac{1 + \log n}{n},
\end{equation}
proving the inequality from the first point.

For the third point, let us denote with tildes the collapsed chain where $A$ is reduced to a point which we denote $\Tilde{A}$ (which, if $x$ and $y$ are neighbours, has a loop at $\Tilde{A}$), and jumps at rate $1$. Then by \cite{AldousFill}, Lemma 2.11, and transitivity of the original graph $\Gamma$ (which implies in particular that $p_t(x,x) +p_t(x,y) = p_t(y,y) + p_t(y,x)$), we have
\begin{equation}
    \bbE_\pi T_A = \bbE_{\Tilde{\pi}} T_{\Tilde{A}} = \frac{1}{\Tilde{\pi}(\Tilde{A})}\Tilde{G}(\Tilde{A}, \Tilde{A}) = \frac{n}{2}\int_0^\infty p_t(x,x) +p_t(x,y) \mathrm{d}t = \frac{n}{2}(G(o,o) + G(x,y)).
\end{equation}
We conclude that
\begin{equation}
    q_A = \frac{\bbE_\pi T_o}{\bbE_\pi T_A} =\frac{nG(o,o)}{\frac{n}{2}(G(o,o) + G(x,y))} = \frac{2}{ 1+\frac{G(x,y)}{G(o,o)}}. \qedhere
\end{equation}
\end{proof}
We first show that the tail of the integral defining the Green function above decreases exponentially fast.

\begin{lemma}\label{L:poincaré tail integral Green function}
Let $d\geq 2$ be an integer. There exist constants $c_1,c_2$ depending on $d$ such that the following holds. Let $b\in [1,\infty)$ and $\Gamma$ be a finite (connected) vertex-transitive graph of degree $d$. Then, writing again $D = \Diam(\Gamma)$ and $n = |\Gamma|$, we have
\begin{equation}
    \int_{bD^2}^\infty h_t(o,o) \mathrm{d}t\leq c_1 e^{-c_2b} D^2 h_{D^2}(o,o).
\end{equation}
Moreover, if $n\leq D^5$, we have
\begin{equation}
    \int_{bD^2}^\infty h_t(o,o) \mathrm{d}t\leq c_1 e^{-c_2b}\frac{D^2}{n}.
\end{equation}
\end{lemma}
\begin{proof}
By spectral estimates, we have for all $t,s\geq 0$ that
\begin{equation}
    h_{t+s}(o,o) \leq e^{-s/\trel}h_t(o,o).
\end{equation}
Recalling from \eqref{treldiff} that $\trel \leq dD^2$ and using that $t\mapsto h_t(o,o)$ is decreasing, we deduce that for any $j\geq 0$,
\begin{equation}
    \int_{(j+1)D^2}^{(j+2)D^2}h_t(o,o) \mathrm{d}t
    \leq \pg \int_{D^2}^{2D^2} h_{t}(o,o) \mathrm{d}t\pd e^{-jD^2/\trel} \leq  D^2 h_{D^2}(o,o) e^{-j/d}.
\end{equation}
Assume that $b$ is an integer. Splitting the integral into parts of length $D^2$ we get
\begin{equation}
\begin{split}
     \int_{bD^2}^{\infty} h_t(o,o) \mathrm{d}t = \sum_{j = b-1}^{\infty} \int_{(j+1)D^2}^{(j+2)D^2}h_t(o,o) \mathrm{d}t & \leq  D^2 h_{D^2}(o,o) \sum_{j = b-1}^{\infty} e^{-j/d} \\
     & = \frac{ e^{-(b-1)/d}}{1-e^{-1/d}}D^2 h_{D^2}(o,o).
\end{split}
\end{equation}
This proves the first point when $b$ is an integer, and the claim for real $b$ follows immediately, since $b\mapsto \int_{bD^2}^{\infty}h_t(o,o) \mathrm{d}t$ is decreasing. The second point follows, since if $n\leq D^5$ we have $h_{D^2}(o,o) \leq p_{D^2}(o,o) \lesssim 1/n$ by Proposition \ref{P: bounds on return probabilities}.
\end{proof}

\begin{prop}\label{P:lowerboundgreenfunctionsmallmacro}
Let $(\Gamma)$ be a sequence of finite (connected) vertex-transitive graphs of fixed degree $d$, such that $|\Gamma| = n\leq D^{3}$. There exist positive constants $\eta = \eta(d), C = C(d)$ such that for all $x \in B(o,\eta D)$,
\begin{equation}
G(o,x) \geq C \frac{D^2}{n}.   
\end{equation}
\end{prop}
\begin{proof}
Let $\varepsilon>0$, to be chosen later.
By Lemma \ref{L:poincaré tail integral Green function}, we can choose $b = b(d,\varepsilon)\geq 1$ large enough such that
\begin{equation}\label{e:boundtailgreenfunction}
    \int_{bD^2}^\infty h_t(o,o) \mathrm{d}t\leq \varepsilon\frac{D^2}{n}.
\end{equation}
We hence have
\begin{equation*}
    G(o,x) = \int_{0}^{\varepsilon D^2} h_t(o,x) \mathrm{d}t+ \int_{\varepsilon D^2}^{b D^2} h_t(o,x) \mathrm{d}t+ \int_{b D^2}^{\infty} h_t(o,x) \mathrm{d}t.
\end{equation*}
Since for every $t$, $h_t(o,x) \geq -1/n$, and $h_t(o,x) \geq - h_t(o,o)$, we deduce that
\begin{align*}
    G(o,x) & \geq \int_{\varepsilon D^2}^{b D^2} h_t(o,x) \mathrm{d}t-2\varepsilon \frac{D^2}{n}.
\end{align*}
By Proposition \ref{P:gradient}, we have
\begin{equation}
\int_{\varepsilon D^2}^{b D^2} h_t(o,x) \mathrm{d}t\geq \int_{\varepsilon D^2}^{b D^2} h_t(o,o) \mathrm{d}t- d(o,x)\int_{\varepsilon D^2}^{b D^2} \sqrt{\frac{d}{et}} h_{t/2}(o,o) \mathrm{d}t.
\end{equation}
By Proposition \ref{P:lowerbounddiagonalheatkernel}, we have, assuming without loss of generality that $\varepsilon\leq a/2$, that
\begin{equation}
    \int_{\varepsilon D^2}^{b D^2} h_t(o,o) \mathrm{d}t\geq \int_{aD^2/2}^{a D^2} h_t(o,o) \mathrm{d}t\geq \frac{a}{4}\frac{D^2}{n}.
\end{equation}
Finally, by Proposition \ref{P: bounds on return probabilities}, there exists a constant $C = C(d, \varepsilon)$ such that $h_{\varepsilon D^2/2} \leq C(d,\varepsilon)/n$, so we have, setting $C' = Cb\sqrt{\tfrac{d}{e\varepsilon}}$, for every $x\in B(o,\tfrac{\varepsilon}{C'}D)$, that 
\begin{equation}
\begin{split}
    d(o,x)\int_{\varepsilon D^2}^{b D^2} \sqrt{\frac{d}{et}} h_{t/2}(o,o) \mathrm{d}t& \leq d(o,x)(b-\varepsilon)D^2\sqrt{\frac{d}{e\varepsilon D^2}} h_{\varepsilon D^2/2}(o,o) \\ & \leq C' \frac{d(o,x)}{D} \frac{D^2}{n} \leq \varepsilon \frac{D^2}{n}.
\end{split}
\end{equation}

The arguments above show that taking $\varepsilon = a/24$, we have, for some $\eta >0$ (depending only on $d$), that for every $x \in B(o,\eta D)$,
\begin{equation}
    G(o,x) \geq \pg\frac{a}{4}-3\varepsilon\pd\frac{D^2}{n} \geq \frac{a}{8}\frac{D^2}{n},
\end{equation}
which concludes the proof. 
\end{proof}

\begin{lemma}\label{lem:Greenmacro1}
    Let $(\Gamma)$ be a sequence of finite (connected) vertex-transitive graphs of fixed degree $d$, such that $|\Gamma|=n\leq D^{3}$. There exists a constant $C = C(d) >0$ such that for all $x \in \Gamma$,
\begin{equation}
G(o,x) \geq - C \frac{D^2}{n}.   
\end{equation}
\end{lemma}
\begin{proof}
We have
\begin{equation}
    G(o,y) = \int_0^\infty h_t(o,x) \mathrm{d}t\geq \int_{0}^{D^2}\pg - \frac{1}{n}\pd \mathrm{d}t - \int_{D^2}^\infty h_t(o,o)\mathrm{d}t.
\end{equation}
Applying Lemma \ref{L:poincaré tail integral Green function} to the second integral on the right hand side concludes the proof.
\end{proof}
Let us set, for every $x\in \Gamma$ and $c>0$,
\begin{equation}
    S_{o,c} := \ag z \in \Gamma\du G(o,z) \leq - c \frac{D^2}{n} \ad.
\end{equation}

\begin{corollary}\label{cor:manypointshaveasmallgreenfunction}
There exist constants $c = c(d) > 0$ and $c' = c'(d)>0$ such that for all finite (connected) vertex-transitive graphs of degree $d$ such that $n = |\Gamma|\leq D^3$, we have
\begin{equation}\label{e:macrosetssmallgreenfunction1}
    |S_{o,c}| \geq c' n.
\end{equation}
\end{corollary}
\begin{proof}
Let us call for this proof $C_1$ the constant from Proposition \ref{P:lowerboundgreenfunctionsmallmacro} and $C_2$ the constant from Lemma \ref{lem:Greenmacro1}. Also, by \cite[Proposition 6.1]{TesseraTointonIsop}, there is a constant $C_3 = C_3(d, \eta)$ such that $V(\eta D) \geq C_3 n$.  We hence have
\begin{equation*}
\begin{split}
    0  = \sum_{x\in\Gamma}G(o,x) & = \sum_{\underset{G(o,x)\geq 0}{x\in \Gamma}}G(o,x) + \sum_{\underset{0>G(o,x)> -\tfrac{C_1C_3}{2}\tfrac{D^2}{n}}{x\in \Gamma}}G(o,x) + \sum_{\underset{-\tfrac{C_1C_3}{2}\tfrac{D^2}{n} \geq G(o,x)}{x\in \Gamma}}G(o,x) \\
    & =: A_1 + A_2 + A_3.
\end{split}
\end{equation*}
The first term can be lower bounded using Proposition \ref{P:lowerboundgreenfunctionsmallmacro}:
\begin{equation}
    A_1 \geq \sum_{x\in B(o,\eta D)}G(o,x) \geq C_1 \frac{D^2}{n}V(\eta D) \geq C_1C_3 D^2.
\end{equation}
Since the second sum is over at most $n$ terms, we have the raw bound
\begin{equation}
    A_2 \geq -\frac{C_1C_3}{2}D^2.
\end{equation}
Finally, since by Lemma \ref{lem:Greenmacro1} for every $x$, $G(o,x)\geq - C_2\tfrac{D^2}{n}$,
\begin{equation}
    A_3 \geq -\bg S_{o,c} \bd C_2 \frac{D^2}{n}
\end{equation}
where $c = C_1 C_3 / 2$.

All in all, dividing by $D^2$, we have proved that 
\begin{equation}
   0 \geq C_1C_3 -\frac{C_1C_3}{2} - \frac{C_2}{n}|S_{o,c}| = c - \frac{C_2}{n}|S_{o,c}|,
\end{equation}
i.e., setting $c := \frac{C_1C_3}{2}$ and $c':= c/C_2 = \frac{C_1C_3}{2C_2}$, that
\begin{equation}
    |S_{o,c}| \geq c'n. \qedhere
\end{equation}
\end{proof}

Roughly speaking, one should think of the set $S_{o,c}$ as the set of points which are ``far away'' from $o$. 
However, it turns out that the variations of the Green function are \textit{macroscopically continuous}, and hence that the set $S_{o,c}$ contains a ball of size $\asymp n$. This is the content of Proposition \ref{P:greenfunctionmacroscopically continuous} and Proposition \ref{P:greenfunctiontwoballsmacrofar}.

\begin{prop}\label{P:greenfunctionmacroscopically continuous}
Let $(\Gamma)$ be a sequence of finite (connected) vertex-transitive graphs of fixed degree $d$, such that $|\Gamma| = n\leq D^{3}$. Let $0<\rho<1$. For every $\varepsilon>0$, there exists a constant $\delta = \delta(\rho, \varepsilon, d)>0$ such that for every $x,y,z\in\Gamma $ such that $d(o,x) \geq \rho D$, $d(x,y) \leq \delta D$, and $d(o,z) \leq \delta D$, we have
\begin{equation}
|G(o,x)-G(z,y)| \leq \varepsilon \frac{D^2}{n}.   
\end{equation}
\end{prop}

\begin{proof}
    Let $\varepsilon>0$. We first assume that $z = o$ and $D^2 \leq n \leq D^3$. Let $a>0$ be such that $\frac{e^{-c\rho^2/a}}{a^{1/2}}\leq \varepsilon$, and let $b>0$ be such that \eqref{e:boundtailgreenfunction} holds. Let $x,y\in B(o,\rho D)^c$. Denote for this proof, for $0<t_1<t_2\leq \infty$,  $I(t_1,t_2) := \int_{t_1}^{t_2} |h_t(o,x)-h_t(o,y)| \mathrm{d}t$. By the triangle inequality,
    \begin{equation}
        |G(o,x)-G(o,y)|  \leq I(o,D) +I(D,aD^2) + I(aD^2, b D^2) + I(bD^2, \infty).
    \end{equation}
By the Carne--Varopoulos inequality, we have $I(o,D) = o(D^2/n)$. Now observe that for all $t$, 
\begin{equation*}
|h_t(o,x) - h_t(o,y)| = |p_t(o,x) - p_t(o,y)| \leq \max(p_t(o,x), p_t(o,y)).    
\end{equation*}
To bound $p_t(o,x)$ (and $p_t(o,y)$) we proceed as in the proof of Proposition \ref{P:lowerbounddiagonalheatkernel}  to get \eqref{e:bound on pt(o,x)}. (Now $d(o,x) \geq \rho D$ instead of having $d(o,x) = D$.)
Let $c$ be as in \eqref{e:bound on pt(o,x)}. Observing that the bound on $p_t(o,x)$ in the integral is maximised at $t = aD^2$, we have 
\begin{equation}
    I(D,aD^2) \lesssim \int_D^{aD^2} \frac{1}{a^{3/2}} \frac{1}{n} \exp\pg-c\frac{\rho^2}{a}\pd \mathrm{d}t\leq \frac{e^{-c\rho^2/a}}{a^{1/2}}\frac{D^2}{n} \leq \varepsilon\frac{D^2}{n}.
\end{equation}
Then by the triangle inequality (and since for all $t$, $h_t(o,o) = \max_x |h_t(o,x)|$), we have $I(bD^2, \infty) \leq 2\varepsilon D^2/n$. Finally, by Proposition \ref{P:gradient} with $s = aD^2/2$ and making a change of variables, we have
\begin{equation}
    I(aD^2,bD^2) \leq \sqrt{d/e}\frac{d(x,y)}{D}\int_{aD^2/2}^{(b-a/2)D^2}h_t(o,o)\mathrm{d}t\lesssim_{a,d} \frac{d(x,y)}{D} \frac{D^2}{n}, 
\end{equation}
so for $\delta>0$ fixed small enough, if $d(x,y)\leq \delta D$, we also have $I(aD^2,bD^2) \leq \varepsilon D^2/n$.
Putting everything together, we have proved that for $\delta>0$ fixed as above and $x,y$ such that $d(x,y)\leq \delta D$, we have $|G(o,x)-G(o,y)| \leq (4\varepsilon +o(1)) \frac{D^2}{n}$. We can bound similarly $|G(o,y) - G(z,y)|$, and hence by the triangle inequality we have $|G(o,x)-G(z,y)| \leq (8\varepsilon +o(1)) \frac{D^2}{n}$, concluding the proof when $D^2\leq n\leq D^3$. The proof for $D\leq n\leq D^2$ is analogous.
\end{proof}

\begin{prop}\label{P:greenfunctiontwoballsmacrofar}
There exist constants $c, \rho > 0$ (depending on $d$) such that for all finite (connected) vertex-transitive graphs of degree $d$ such that $n = |\Gamma|\leq D^3$, there exist two (disjoint) balls $S_1$ and $S_2$ of radius $\rho D$ such that for every $x\in S_1$ and $y \in S_2$,  
\begin{equation}
    G(x,y) \leq -c \frac{D^2}{n}.
\end{equation}
In particular, for such $c,\rho$, $S_{o,c}$ contains a ball of radius $\rho D$.
\end{prop}
\begin{proof}
    By Corollary \ref{cor:manypointshaveasmallgreenfunction}, there exists $c>0$ and $x\in \Gamma$ such that $G(o,x) \leq -cD^2/n$. Therefore, $d(o,x) \geq \eta D$ where $\eta$ is as in Proposition \ref{P:lowerboundgreenfunctionsmallmacro}. Therefore, by Proposition \ref{P:greenfunctionmacroscopically continuous} (where the role of $\rho$ there is played by $\eta$ here) we can find $\delta>0$ such that for every $y\in B(x,\delta D)$, and $z\in B(o, \delta D)$
    \begin{equation}
        |G(z,y) - G(o,x)| \leq \frac{c}{2} \frac{D^2}{n}.
    \end{equation}
It follows that for such $y,z$, we have
$G(z,y)\leq -\frac{c}{2}\frac{D^2}{n}$, concluding the proof.
\end{proof}
We will also need the following bound on $\thit$.
\begin{prop}\label{prop:subcritical bound diagonal green function}
Let $C>0$. Let $(\Gamma)$ be a collection of finite (connected) vertex-transitive graphs of fixed degree $d$, such that $n=|\Gamma| \leq CD^2\log n$. Then 
\begin{equation}
    \thitav \lesssim_{d,C} D^2 \log n.
\end{equation}
\end{prop}
\begin{proof}
    Recall that we have $\thitav = nG(o,o)$ from Lemma 2.11 in \cite{AldousFill}. We split the proof into two cases, depending on the diameter of $\Gamma$.
    
First, if $D^2 \leq n \leq CD^2\log n$, we write $n = D^2R$. Integrating the bound on $p_t(o,o)$ from Proposition \ref{P: bounds on return probabilities}, and bounding the tail of the integral with Proposition \ref{L:poincaré tail integral Green function}, we have 
\begin{equation}
    G(o,o) = \int_0^{\infty} h_t(o,o) \mathrm{d}t\lesssim_d 1+ \frac{1}{R} + \frac{\log(D/R)}{R} + \frac{D^2}{n} \lesssim_{d,C} \frac{\log n}{R} = \frac{1}{n} D^2\log n, 
\end{equation}
as desired.

Suppose now that $D \leq n \leq D^2$, and let $H$ such that $n = DH$. Proceeding as above, we get
\begin{equation}
    G(o,o) \lesssim_d 1 + \log H + \frac{D}{H} + \frac{D^2}{n} \lesssim \log n + \frac{D}{H}.
\end{equation}
We hence have
\begin{equation}
    nG(o,o) \lesssim DH\log n + D^2 \lesssim D^2 \log n,
\end{equation}
which concludes the proof.
\end{proof}

\subsection{Proof that the diameter condition is necessary in Theorem \ref{T:main} and Theorem \ref{T:uncovered_main}}
In this section, we prove that when the diameter condition \eqref{e:DC} does not hold, then the renormalised cover time $\tfrac{\taucov}{\thit} - \log |\Gamma|$ does not have asymptotic Gumbel fluctuations. 
We first prove that if $n\leq D^5$, Proposition \ref{prop:comparison average and max hitting times} still holds when the average hitting time is replaced by the quasi-stationary hitting time.

\begin{prop}\label{prop:comparison QS and max hitting times} Let $d\geq 2$. Uniformly over all finite (connected) vertex-transitive graphs of degree $d$ such that $n = |\Gamma| \leq D^5$, we have
\begin{equation}
    \thit - \bbE_\alpha T_o \asymp_d D^2,
\end{equation}
where $\alpha$ is the quasi-stationary distribution associated with $o$.
\end{prop}

To prove Proposition \ref{prop:comparison QS and max hitting times} we will need a few lemmas. For $1\leq r \leq D$, set $B_r := B(o,r)$,  and $\beta(r) = 1/\bbE_{\alpha_{B_r^c}}T_{B_r^c}$. 

\begin{lemma}
\label{lem:ondiaglowervolumedubling}  For every $2 \le \ell < D$ and $t \ge 0$, we have, writing $r:=\lfloor \ell /2 \rfloor$, that
\begin{equation}
\label{e:tailofexitball}
\P_o[T_{B(o,\ell)^c}>t] \ge \max_{v \in B(o,r)}\P_v[T_{B(o,r)^c}>t] \ge e^{-t\beta(r) },
\end{equation}
\begin{equation}
\label{e:tailofexitball2}
P_t(o,o) \ge \frac{\P_o[X_t \in B(o,\ell) ]}{V(\ell)} \ge \frac{\P_o[T_{B(o,\ell)^c}>t] }{V(\ell)}\ge\frac{e^{-t\beta(r) }}{V(\ell)},
\end{equation}
\begin{equation}
\label{e:tailofexitball3}
 \beta(\ell)\le \frac{4V(\ell )}{V(r) }\ell^{-2}. 
 \end{equation}
\end{lemma}
\begin{proof} Let $v \in B(o,r)$. The first inequality in \eqref{e:tailofexitball} follows by noting that if $X_0=v  $ then $d(X_t,X_0) \le 2r \le \ell $ for all $t< T_{B(o,r)^c}$ and so $T_{B(o,r)^c} \le T_{B(v,\ell)^c}$. It follows by transitivity that $\P_v[T_{B(o,r)^c}>t] \le \P_v[T_{B(v,\ell)^c}>t]=\P_o[T_{B(o,\ell)^c}>t]$ for all $t \ge 0$, as desired. 

Let $A = B(o,r)^c$, and $\alpha = \alpha_{A}$ be the quasi-stationary distribution associated with $A$. The second inequality in \eqref{e:tailofexitball} follows from the fact that
\begin{equation}
    \max_{v \in B(o,r)}\P_v[T_{B(o,r)^c}>t] \ge\P_{\alpha_{r}}[T_{B(o,r)^c}>t] =e^{-t\beta(r) }.
\end{equation}
The first inequality in \eqref{e:tailofexitball2} follows from the fact that by transitivity $p_t(o,o)=\max_{x\in\Gamma} p_t(o,x)$, while the second is trivial and the third is exactly  \eqref{e:tailofexitball}.

We now prove \eqref{e:tailofexitball3}. For a distribution $\nu$ on $V$ and $g,g':V \to \bbR$ we write $\langle g,g' \rangle_{\nu}:=\mathbb{E}_{\nu}[gg']=\sum_{x\in V}\nu(x)g(x)g'(x) $.  Let $\mathcal{C}_{\ell}$ be the collection of all $g:V \to \mathbb{R}$ whose support $\{x \in V:g(x) \neq 0\}$ is non-empty and contained in $B(o,\ell)$, and denote by $\pi_{\ell}$ the uniform distribution on   $B(o,\ell)$. Observe that for all $g \in \mathcal{C}_{\ell}$ we have that
\begin{equation}\label{e:scalarproductcourantfischerproof}
\begin{split}
\langle(I-P_{B(o,\ell)})g,g \rangle_{\pi_{\ell}} & =\langle(I-P)g,g \rangle_{\pi_{\ell}}=\langle(I-P)g,g \rangle_{\pi}/\pi(B(o,\ell))
\\ &=\frac{1}{2\pi(B(o,\ell))} \mathbb{E}_{\pi}\left[\left(g(X_{0})-g(X_{1}) \right)^2  \right]
\\ &=\mathbb{E}_{\pi_{\ell}}\left[\left(g(X_{0})-g(X_{1}) \right)^2 \frac {1+\mathbf{1}\{X_1 \notin B(o,\ell)\}}{2} \right]  
\end{split}
\end{equation}
(c.f.\ \cite[Lemma 13.6]{LivreLevinPeres2019MarkovChainsAndMixingTimesSecondEdition} for the third equality). Since $\beta(\ell)$ is the spectral gap of the chain killed at $B(o,\ell)^c$, we have by the Courant--Fischer characterization of $\beta(\ell)$ that
\begin{equation}
    \beta(\ell) = \min_{g \in \mathcal{\mathcal{C}_{\ell}}} \left\{ \frac{\langle(I-P_{B(o,\ell)})g,g \rangle_{\pi_{\ell}}}{ \mathbb{E}_{\pi_{\ell}}[g^2]}\right\}. 
\end{equation}
  Using \eqref{e:scalarproductcourantfischerproof}, and since the test function $f \du x \mapsto d(x,B(o,\ell)^c)$ is 1-Lipschitz we get that $\langle(I-P_{B(o,\ell)})f,f \rangle_{\pi_{\ell}}\le 1$. Since moreover $f \in \mathcal{C}_{\ell} $, we have
  \begin{equation*}
    \frac{1}{\beta(\ell)} \ge \frac{\mathbb{E}_{\pi_{\ell}}[f^2] }{\langle(I-P_{B(o,\ell)})f,f \rangle_{\pi_{\ell}}}  \geq \mathbb{E}_{\pi_{\ell}}[f^2] \ge (\ell - r)^2 \pi_{\ell}(B(o,r)) \ge \frac{\ell^2V(r)}{4 V(\ell)},  
  \end{equation*}
where we used the fact that $f(x) \ge \ell-r$ for all $x \in B(o,r)$. This concludes the proof.
\end{proof}

\begin{lemma}\label{lem:positive proba to stay in macroscopic balls}
        Let $(\Gamma)$ be a sequence of finite (connected) vertex-transitive graphs of degree $d$ such that $n = |\Gamma| \leq D^5$, and let $a, b>0$. There exists $C = C(a,b,d) >0$ such that
  \begin{equation}
    \bbP_o(T_{B_{aD}^c} > bD^2) \geq C.
\end{equation}  
\end{lemma}
\begin{proof} 
Let $\ell := 2\left\lfloor aD/2 \right \rfloor$, $r = \ell/2$, and $t = bD^2$, so that $\ell$ is even and $\bbP_o(T_{B_{aD}^c} > bD^2) \geq  \bbP_o(T_{B_{\ell}^c} > t)$. By \cite[Corollary 1.5]{TesseraTointonfinitary}, we have $V(r) \gtrsim n$. Plugging this into \eqref{e:tailofexitball3}, we deduce that $\beta(r) \lesssim r^{-2}$, and therefore that $e^{-\beta(r)t} \gtrsim 1$. Plugging this into \eqref{e:tailofexitball} concludes the proof.
\end{proof}

\begin{lemma}\label{lem:QSdistrib positive proba on macroscopic balls}
    Let $(\Gamma)$ be a sequence of finite (connected) vertex-transitive graphs of degree $d$ such that $n = |\Gamma| \leq D^5$, and let $\varepsilon>0$. There exists $c>0$ such that
\begin{equation}\label{e:quasistat distrib of a small ball positive}
    \alpha_o(B_{\varepsilon D}) \geq c(\varepsilon).
\end{equation}
\end{lemma}
\begin{proof}
The intuition of this proof is that if a walk starts more than $\eps D$ away from $o$, there should be a positive probability that, after time of order $D^2$, it finds itself in the ball $B_{\varepsilon D}$ (indeed, mixing times are hitting times of large sets, see \cite{PeresSousi2015}) but has not touched $o$ yet (since starting from the quasi-stationary distribution $\alpha_o$ and conditioning on hitting $o$, the distribution at any given later time remains $\alpha_o$). 

Suppose that $x \in B_{\eps D}^c$. Then 
\begin{equation}
    \bbP_x(T_{B_{\varepsilon D/2}} \leq \tmixinfty(1/2)) \geq \bbP_x(X_{\tmixinfty(1/2)} \in B_{\varepsilon D/2}) \geq \frac{1}{2}\frac{V(\varepsilon D/2)}{n}\gtrsim 1. 
\end{equation}
Furthermore, by Lemma \ref{lem:positive proba to stay in macroscopic balls} (recalling that $\tmixinfty(1/2) \asymp D^2$ by Proposition \ref{P:tmix asymp diam squared}), 
given that the walk enters $B_{\eps D/2}$, the conditional probability that it remains in the annulus $B_{3\varepsilon D/4}\backslash B_{\varepsilon D/4}$ until time $\tmixinfty(1/2)$ is at least $c$ for some constant $c$. Therefore, setting $t = \tmixinfty (1/2)$,  
\begin{align*}
\alpha ( B_{\eps D} )  = \P_\alpha ( X_t \in B_{\eps D} | T_o \ge t ) 
& \ge \alpha (B_{\eps D}^c)  \E_{\alpha |_{B_{\eps D}^c}} ( X_t \in B_{\eps D} | T_o \ge t )  \\
& \ge \alpha ( B_{\eps D}^c) \inf_{x \in B_{\eps D}^c} \P_x ( X_t \in B_{\eps D}, T_o \ge t)\\
& \gtrsim  \alpha ( B_{\eps D}^c).
\end{align*}
Thus $\alpha (B_{\eps D})$ is bounded away from zero, as desired.
\end{proof}

\begin{proof}[Proof of Proposition \ref{prop:comparison QS and max hitting times}]
By Proposition \ref{P:Green function identities}, for $x,y \in \Gamma$, we have $\bbE_x T_y = \bbE_\pi T_o - nG(x,y)$, so $\bbE_x T_y \leq \bbE_\pi T_o$  if and only if $G(x,y) \geq 0$. Combining this with Proposition \ref{P:lowerboundgreenfunctionsmallmacro}, there exists $\eta>0$ such that for every $x\in B(o,\eta D)$,
\begin{equation}
    \bbE_x T_o \leq \bbE_\pi T_o.
\end{equation}
Therefore, since $\alpha(B_{\eta D})\gtrsim 1$ by Lemma \ref{lem:QSdistrib positive proba on macroscopic balls}, and $\thit - \bbE_\pi T_o \asymp D^2$ by Propositions \ref{prop:comparison average and max hitting times} and \ref{P:tmix asymp diam squared}, we conclude that
\begin{equation}
    \thit - \bbE_\alpha T_o \geq \sum_{x\in B_{\eta D}} \alpha(x)\pg \thit - \bbE_\pi T_o \pd =  \alpha(B_{\eta D})\pg \thit - \bbE_\pi T_o\pd \gtrsim D^2. \qedhere
\end{equation}
\end{proof}

\begin{prop}\label{P:proof no Gumbel fluct max hit time}
Let $(\Gamma)$ be a collection of finite (connected) vertex-transitive graphs of fixed degree $d$, and assume that as $n =|\Gamma| \to \infty$, 
\begin{equation}\label{e:lesssim sous-critique}
D^2 \log n \gtrsim n.
\end{equation}
Then there exists a constant $\kappa <1$ such that for $s\in\bbR$ large enough, we have at time $\tsmax = \thit((\log n)+ s)$ 
\begin{equation}
    \bbP(\taucov > \tsmax) \leq \kappa \pg 1-e^{-e^{-s}}\pd.
\end{equation}
In particular $\tfrac{\taucov}{\thit} - \log |\Gamma|$ asymptotically does not have Gumbel fluctuations, and for $s$ large enough, $\dtv ( \cL(U(\tsmax)), \mu_s^{\otimes \Gamma} )$ does not converge to $0$.
\end{prop}

\begin{proof}
From Proposition \ref{prop:comparison QS and max hitting times}, there exists a constant $c>0$ such that
    \begin{equation}
        \frac{\thit}{\bbE_\alpha T_o} \geq 1 + c\frac{D^2}{\bbE_\alpha T_o}.
    \end{equation}
It follows from Proposition \ref{prop:subcritical bound diagonal green function} that there exists a constant $C>0$ (depending on $d$, on $c$ and the implicit constant in \eqref{e:lesssim sous-critique}), such that
 \begin{equation}\label{e:thitthitavratiolowerbound}
     \frac{\thit}{\bbE_\alpha T_o} \geq 1 + \frac{C}{\log n}.
 \end{equation}
From a simple union bound and \eqref{e:thitthitavratiolowerbound}, we deduce that
\begin{equation}
    \bbP(\taucov > \tsmax) \leq n\bbP_\pi(T_o > \tsmax)  \leq n \bbP_\alpha(T_o > \tsmax)   = n\exp\pg - \frac{\thit}{\bbE_\alpha T_o}\pg (\log n) + s\pd\pd  \leq e^{-s}e^{-C}. 
\end{equation}
Hence, for $s$ larger than some $s_0 = s_0(C)$, we have
\begin{equation}
    \bbP\pg\frac{\taucov}{\thit} - \log |\Gamma| > s\pd = \bbP(\taucov > \tsmax)  \leq e^{-C}\bbP\pg \chi > s \pd  \leq \kappa\pg 1-e^{-e^{-s}}\pd,
\end{equation}
where $\kappa :=\tfrac{1+e^{-C}}{2}< 1$. This proves that $\tfrac{\taucov}{\thit} - \log |\Gamma|$ does not have Gumbel fluctuations. Since $\bbP(\taucov > \tsmax) = \bbP(U(\tsmax) \ne \eset)$, it also shows that for $s\geq s_0$, we have 
\begin{equation}
    \dtv ( \cL(U(\tsmax)), \mu_s^{\otimes \Gamma}) \geq (1-\kappa )\pg 1-e^{-e^{-s}}\pd. \qedhere
\end{equation}
\end{proof}

\report{
We can now complete the proof of Theorems \ref{T:main} and \ref{T:uncovered_main}.
\begin{proof}[Proof of Theorems \ref{T:main} and \ref{T:uncovered_main}]
    The facts that the diameter condition implies Gumbel fluctuations and even a product structure for the uncovered set at the desired times were proved in Propositions \ref{prop: DC implies Gumbel fluctuations}  and \ref{prop: theorem intro uncovered set first part: DC implies decorrelation}, and the reverse directions were proved above in Proposition \ref{P:proof no Gumbel fluct max hit time}. 
\end{proof}}

\subsection{Proof that the diameter condition is necessary in Theorem \ref{T:uncovered intro t hit av}}

\label{SS:diam_nec_uncovered_av}
In this section, we assume that $\trel = o( \thit)$. We want to show that the law of the uncovered set at time $\tsexp$ is far, in the total variation sense, from the product measure $\mu_s^{\otimes \Gamma}$. Note that this immediately implies Theorem \ref{T:uncovered intro t hit av}: indeed, since $\trel \le d D^2$ by \eqref{lem:trel asymp diam squared} and $\thit \ge \thitav = n G(o,o) \ge n(1 - o(1/\log n))$ by Proposition \ref{P:Green function identities}, the assumption $D^2 = o(n)$ implies that $\trel =o(\thit)$. Therefore, it suffices to prove the result under the sole assumption $\trel = o(\thit)$.

To prove that the law of the uncovered set is far from a product measure when the diameter condition fails, we might initially be tempted to show that the uncovered set is ``too'' clustered, i.e.\ the probability that two relatively nearby points are uncovered is larger than it should be under the independent scenario. However, this turns out to be very difficult to control as the contribution to the moments of order $k$ of the size of the uncovered set coming from nearby points start exploding when the diameter condition fails. We cannot translate this into 
estimates about events of positive probability for the uncovered set; roughly speaking, either the Bonferroni inequality goes in the wrong direction, or one would need to keep track of moments of higher order and compare how they blow up. 

Instead, we show that points that are sufficiently far apart are negatively correlated. It turns out that this enables us to use the Bonferroni inequality (i.e., union bound) as this is an upper bound. 

\begin{prop}\label{P:no DC and trel <<< thit implies no product measure}
 Let $\gamma >0$, and $(\Gamma)$ be a collection of finite (connected) vertex-transitive graphs of fixed degree $d$, such that $n = |\Gamma| \leq \gamma D^2\log n$ and $\trel = o(\thit)$, and let (for every $s\in\bbR$)
$\mu_s^{\otimes \Gamma}$ denote the product over all the vertices of the graph of the Bernoulli law $\mu_s$ with parameter $e^{-s}/|\Gamma|$. Then for every fixed $s\in \bbR$, there exists a constant $c^* = c^*(s,d, \gamma)$ such that as $|\Gamma|\to\infty$, for $t\in\ag \tsav, \tsexp, \tsmax\ad$, and under $\P_\pi$ (and thus also under $\P_x$ for any $x\in \Gamma$),
  \begin{equation*}
  \dtv ( \cL(U(t)), \mu_s^{\otimes \Gamma} ) \geq c^* + o(1).
  \end{equation*}
\end{prop}

We start by comparing $\tsav$ and $\tsexp$. Recall that $\tsav = \thitav((\log n) + s)$ and $\tsexp$ is such that $\bbE|U(\tsexp)| = e^{-s}$.
\begin{lemma}\label{lem:comparasion ts et tdoubleprime s}
    Let $s\in \bbR$. For $n$ large enough, we have $\bbE |U(\tsav)|\geq e^{-s}$, or in other words $\tsav \leq \tsexp$.  
\end{lemma}

(Note that by transitivity there is no need to specify the starting point in the above statement.)

\begin{proof}
Let $\alpha$ be the quasi-stationary distribution associated with the set $A=\ag o \ad$, and let $s\in\bbR$. We have
\begin{equation}
\begin{split}
 \bbE|U(\tsav)| = n\bbP_\pi(T_o > \tsav) & = n\frac{\bbP_\pi(T_o > \tsav)}{\bbP_\alpha(T_o > \tsav)}e^{-\tsav/\bbE_\alpha T_o} \\ & = n\frac{\bbP_\pi(T_o > \tsav)}{\bbP_\alpha(T_o > \tsav)}e^{-\tfrac{\bbE_\pi T_o}{\bbE_\alpha T_o}((\log n)+ s)}.   
\end{split}
\end{equation}
Here we use a further refinement of the Aldous--Brown approximations. Let $W = \sum_{x\in \Gamma} \frac{\alpha(x)^2}{\pi(x)}$ and $\theta = \frac{1-W}{W}$. The quantity $\theta$ is the quasi-stationary default of stationarity for the set $B=\Gamma\backslash\ag o\ad$, denoted $R_B$ in \cite{BerestyckiHermonTeyssier2025AB}. By \cite[Lemma 3.2 and Lemma 3.3]{BerestyckiHermonTeyssier2025AB}, we have
\begin{equation}
    \frac{\bbP_\pi(T_o > \tsav)}{\bbP_\alpha(T_o > \tsav)} \geq 1 - \theta \quad \quad \text{ and } \quad \quad  \frac{\bbE_\pi T_o}{\bbE_\alpha T_o} \leq 1 - \theta + \theta \frac{\trel}{\bbE_\alpha T_o}.
\end{equation}
Since $\trel = o(\thit)$ by hypothesis, we have $\tfrac{\trel}{\bbE_\alpha T_o} \to 0$. Since $\frac{\bbE_\pi T_o}{\bbE_\alpha T_o} = 1+o(1)$ (for instance by \cite{AldousBrown1992} since $\trel=o(\thit)$), we also know that $\theta = o(1)$, and since $\frac{\bbE_\pi T_o}{\bbE_\alpha T_o} <1$, we have $\theta>0$. Therefore, using the inequality $e^x \geq 1+x$, valid for all $x \in \bbR$, we have 
\begin{equation}
    \bbE|U(\tsav)| \geq e^{-s}(1-\theta)(1+\theta\log n + o(\theta \log n)) = e^{-s}(1 + \theta \log n + o(\theta \log n)),
\end{equation}
which is $>e^{-s}$ for $|\Gamma|$ large enough. This concludes the proof.
\end{proof}
We now prove a technical lemma.
\begin{lemma}\label{L:bound hitting proba when no DC nor trel <<< thit}
    Let $(\Gamma)$ as in Proposition \ref{P:no DC and trel <<< thit implies no product measure}, and $c,S_1,S_2$ as in Proposition \ref{P:greenfunctiontwoballsmacrofar}.
There exists a constant $C>0$ such that for every $(x,y)\in S_1 \times S_2$ and $s\in \bbR$, we have (as $|\Gamma|\to\infty$), setting $A = \{ x, y\}$,
\begin{equation}
    \bbP_\pi(T_A > \tsav) \leq \frac{e^{-2s}}{n^2}e^{-C + o(1)}.
\end{equation}
\end{lemma}
\begin{proof}
Let us fix $s\in\bbR$. For every $(x,y)\in S_1 \times S_2$, we have by \eqref{e:identité capacité}
\begin{equation}
    q_A = \frac{2}{1 + \frac{G(x,y)}{G(o,o)}} \geq \frac{2}{1 - c \frac{D^2}{nG(o,o)}} \geq 2\pg 1 + c \frac{D^2}{nG(o,o)} \pd.
\end{equation}
Moreover, by Theorem \ref{thm: improvement of AB for finite sets under DC}, we have uniformly over all $A\subset \Gamma$ of size 2, 
\begin{equation}
       \bbP_\pi(T_A > \tsav) \leq \bbP_{\alpha_A}(T_A > \tsav) =  \exp\pg-q_A \tsav \pg 1 + O \pg\frac{\trel}{\thit}\pd^2\pd \pd.
\end{equation}
It follows, recalling that $\tfrac{D^2}{nG(o,o)}\asymp \tfrac{\trel}{\thit}$, that there exists a constant $c'>0$ such that for every $(x,y)\in S_1 \times S_2$,
\begin{equation}
     \bbP_\pi(T_A > \tsav) \leq \frac{e^{-2s}}{n^2}\exp\pg- 2c'\frac{\trel}{\thit}(\log n)  \pg 1 + \frac{s}{\log n} + O \pg\frac{\trel}{\thit}\pd\pd \pd.
\end{equation}
From Proposition \ref{prop:subcritical bound diagonal green function}, there exists a constant $C$ such that (recalling that $\trel \asymp D^2$)
\begin{equation}
    2c'\frac{\trel \log n}{\thit} \geq C,
\end{equation}
which allows us to conclude, using that $\trel =o(\thit)$, that
\begin{equation}
    \bbP_\pi(T_A > \tsav) \leq \frac{e^{-2s}}{n^2}e^{-C + o(1)}. \qedhere
\end{equation}
\end{proof}
\begin{proof}[Proof of Proposition \ref{P:no DC and trel <<< thit implies no product measure}]
Let $s\in \bbR$. By definition of the total variation distance, it is enough to find a subset $B$ of $\cP(\Gamma)$ and a constant $c^*$ such that for $t\in\ag \tsav, \tsexp, \tsmax\ad$,
\begin{equation}
 \mu_s^{\otimes \Gamma}(B) - \bbP(U(t) \in B) > c^*.
\end{equation}

Let $c, S_1, S_2$ as in Proposition \ref{P:greenfunctiontwoballsmacrofar}, and $C>0$ from Lemma \ref{L:bound hitting proba when no DC nor trel <<< thit}. Set $b = b(s) = e^{-s}$ and let  $a>0$ be small enough such that $|S_1|\geq an$ (for all $\Gamma$ as in the statement of the proposition), and $1 - 2e^{-ab} + e^{-2ab} \geq (ab)^2(1+e^{-C})/2$. (Such a choice for $a$ is possible since as $a\to 0$, we have $ 1 - 2e^{-ab} + e^{-2ab} = (ab)^2 + O(a^3)$.) Let $S'_1$ and $S'_2$ be sub-balls of (respectively) $S_1$ and $S_2$ such that $|S'_1| = |S'_2| = (a+o(1))n$, and set
\begin{equation}
    B := \ag A\subset \Gamma\du A\cap S'_1 \ne \eset \text{ and } A\cap S'_2 \ne \eset \ad.
\end{equation}
  
By a union bound and Lemma \ref{L:bound hitting proba when no DC nor trel <<< thit}, we get the following upper bound on $\bbP(U(\tsav) \in B)$
\begin{equation}\label{e:bound event two macro balls}
\begin{split}
    \bbP(U(\tsav) \in B) & \leq \sum_{(x,y)\in S'_1\times S'_2}\bbP_\pi \pg \ag x,y \ad \subset U(\tsav) \pd
   \\ & \leq \frac{|S'_1||S'_2|}{n^2} e^{-2s}\pg e^{-C} + o(1)\pd = (ab)^2e^{-C} + o(1). 
\end{split}
\end{equation}
 
Let us now lower bound $\mu_s^{\otimes \Gamma}(B)$.
Let $B_i = \ag A \subset \Gamma\du A\cap S'_i \ne \eset\ad$ for $i\in \ag 1,2 \ad$. Since $B = B_1\cap B_2$, we have
\begin{equation}
     \mu_s^{\otimes \Gamma}(B) = 1- \mu_s^{\otimes \Gamma}(B_1^c\cup B_2^c) = 1 - \pg \mu_s^{\otimes \Gamma}(B_1^c) + \mu_s^{\otimes \Gamma}(B_2^c) - \mu_s^{\otimes \Gamma}(B_1^c\cap B_2^c)\pd. 
\end{equation}
Moreover,
\begin{equation}
    \mu_s^{\otimes \Gamma}(B_1^c) = \mu_s^{\otimes \Gamma}(B_2^c) = \pg 1 - \frac{b}{n}\pd^{(a+ o(1))n} = e^{-ab + o(1)},
\end{equation}
and we have similarly $\mu_s^{\otimes \Gamma}(B_1^c\cap B_2^c) = e^{-2ab + o(1)}$. If follows that
\begin{equation}
    \mu_s^{\otimes \Gamma}(B) \geq 1- 2e^{-ab + o(1)} + e^{-2ab + o(1)},
\end{equation}
and hence, by definition of $a$, that
\begin{equation}\label{e:conclusion proof two macro balls}
\begin{split}
    \mu_s^{\otimes \Gamma}(B) - \bbP(U(\tsav) \in B) & \geq (ab)^2\frac{1+e^{-C}}{2} - (ab)^2e^{-C} + o(1) \\ & = (ab)^2\frac{1-e^{-C}}{2} + o(1).
\end{split}
\end{equation}
Recall that $\tsexp \geq \tsav$ (for $n$ large enough) and $\tsmax\geq \tsav$ (for all $n$), so since $t \mapsto \bbP_\pi(T_A > t)$ is decreasing for every $A$, \eqref{e:bound event two macro balls} also holds with $\bbP(U(\tsexp) \in B)$ or $\bbP(U(\tsmax) \in B)$ on the left hand side. Therefore, \eqref{e:conclusion proof two macro balls} also holds with $\tsexp$ or $\tsmax$ instead of $\tsav$, and the proof is complete.
\end{proof}

\report{We can now prove Theorem \ref{T:uncovered intro t hit av}.
\begin{proof}[Proof of Theorem \ref{T:uncovered intro t hit av}]
The (contraposition of the) implication “(i) $\implies$ (ii)” was proved in Proposition \ref{P:no DC and trel <<< thit implies no product measure}.
Now, by Proposition \ref{P:tdoubleprime s} we have, for each $s\in \bbR$ and under \eqref{e:DC}, that $\tsexp = \tsav(1+o(1/\log n))$. It follows from Proposition \ref{prop: theorem intro uncovered set first part: DC implies decorrelation} that the implication “(ii) $\implies$ (i)” also holds, which concludes the proof.
\end{proof}}

\section{A strongly uniformly transient graph without Gumbel fluctuations}
\label{s:examples}

In this section, we construct an example of a sequence of vertex-transitive graphs satisfying the strong uniform transience (SUT) of \eqref{e:SUT}, but not the diameter condition \eqref{e:DC}. 

The idea is to consider the following graph: for $m \ge 2$ even, let $\Gamma = C_m \times G$ be the Cartesian {product} of a cycle $C_m$ of length $m$, together with an expander Cayley graph $G$  (in fact it will be convenient to take a Ramanujan graph of degree greater than three, see Morgenstern \cite{Morgenstern1994} for an explicit construction generalising the famous construction of Lubotzky, Phillips and Sarnak \cite{LubotzkyPhillipsSarnak1988}) of size $m h$, where $h = h_m \ge 1 $, and say $ h_m \le m$, say. (The interesting examples will be those for which $1\ll h \lesssim \log m$.)  

\report{Intuitively, the graph $\Gamma$ is locally ``very recurrent'' in the $C_m$-direction, but locally ``as transient as can be'' in the $G$-direction (indeed, Ramanujan graphs are locally tree-like and have the fastest possible expansion among graphs of a given size and degree). Below we show that the transient behaviour in the $G$-direction is sufficient to guarantee the SUT property, as soon as $|G| \gg m$, but  the diameter condition is equivalent to $|G| \gg m \log m$. Therefore there is some space between the two conditions and by playing with the sizes of the two factors, one can ensure that SUT holds and simultaneously that the diameter condition is not fulfilled.}

\begin{prop}\label{T:counter}
Let $(\Gamma)$ be as above. Assume that $1\ll h_m \lesssim \log m$ as $m\to \infty$. Then $(\Gamma)$ satisfies the strong uniform transience (SUT) condition but not the diameter condition \eqref{e:DC}.
\end{prop}

\begin{proof}
For a vertex $x \in \Gamma$, let us write $x = (\hat x, \check x ) $ with $\hat x \in C_m$ and $\check x \in G$ so that $\hat x$ denotes the cycle coordinate and $\check x$ the Ramanujan coordinate. Since both $C_m$ and $G$ are vertex-transitive, 
for a continuous-time random walk $X_t = ( \hat X_t, \check X_t)\in \Gamma$, the coordinates $\hat X$ and $\check X$ are in fact independent continuous-time random walks with rates $2/d$ and $1- 2/d$ respectively (where $d = \deg (x)$ is the total degree on $\Gamma$.) 

Since the coordinates are independent, we have  that for every $t\ge 0$, 
\begin{equation*}
p_t(x,y) = \hat p_t(\hat x,\hat y) \check p_t(\check x, \check y),
\end{equation*}
where $\hat p_t$ denotes the transition probabilities of the random walk on the cycle $C_m$ (with rate $2/d$), and $\check p_t$ denote those on the Ramanujan graph $G$ (with rate $1-2/d$).

Furthermore, by Proposition \ref{P:tmix asymp diam squared}, 
$\tmix = \tmix(1/4) \lesssim m^2 $. It is therefore easy to estimate both $\hat p_t$ and $\check p_t$ for $t \le m^2$: namely, on the cycle we know (for instance from Proposition \ref{P: bounds on return probabilities}) that 
\begin{equation*}
\hat p_t (\hat o, \hat o) \lesssim \frac1{\sqrt{t+1}};
\end{equation*}
and on the Ramanujan component, by \eqref{Eq:spec},
\begin{equation*}
|\check p_t(\check o, \check o) - \check \pi (\check o) | \le e^{- \lambda t}
\end{equation*}
where $\lambda>0$ is the spectral gap on the Ramanujan component (which, by definition, is bounded away from zero), and $\check \pi$ is the uniform distribution on $G$. 
Since $\check \pi (\check o) = 1/ (mh)$, we deduce that $\check p_t(\check o, \check o) \le 1/ (mh) + e^{-\lambda t}$.

For each fixed $s >0$
\begin{align*}
\E_o [ L_o(\tmix ) - L_o(s)]  = \int_s^{\tmix} p_t(o, o) \mathrm{d}t
& \lesssim \int_s^{\tmix} \frac1{\sqrt{t+1}} \pg \frac1{mh} + e^{-\lambda t } \pd \mathrm{d}t\\
& \lesssim \frac{\sqrt{\tmix}}{mh} + \int_s^\infty e^{ - \lambda t} \frac1{\sqrt{t+1}} \mathrm{d}t
\end{align*}
Since $\tmix \lesssim m^2$ and $h = h_m \to \infty$, the first term tends to zero. 
The second term does not depend on $m$, and is the integral of a function which is clearly integrable. Consequently, the limsup as $s \to \infty$ is zero. 
Thus
\begin{equation*}
   \limsup_{s\to \infty} \limsup_{n\to \infty} \E_o(L_o(\tmix ) - L_o(s) ) = 0, 
\end{equation*}
which shows \eqref{e:SUT}. Finally, if $h_m \lesssim \log m$, then the diameter condition \eqref{e:DC} fails by definition, since $\Diam(\Gamma) \asymp m$ and $|\Gamma| = m^2h_m$. This concludes the proof.
\end{proof}

\bibliography{bibliographyCTU}
\bibliographystyle{alpha}
\end{document}